%% file: multilevel.tex
\documentclass[final]{siamltex}

\usepackage[margin=1in]{geometry}

\usepackage{amsmath,amssymb,amscd}

\usepackage{hyperref}
\usepackage{xcolor}
\usepackage{graphicx}
\usepackage{subfigure}
\usepackage[title]{appendix}
\usepackage{tikz}
\usepackage{algorithm}
\usepackage{algpseudocode}
\usepackage{booktabs}
\usepackage{tikz}
\usepackage{xparse}
\usepackage{enumitem}

\usetikzlibrary{matrix,backgrounds,decorations.pathreplacing}
\pgfdeclarelayer{myback}
\pgfsetlayers{myback,background,main}

\tikzset{mycolor/.style = {line width=1bp,color=#1}}%
\tikzset{myfillcolor/.style = {draw,color=#1,fill=#1}}%

\NewDocumentCommand{\highlight}{O{blue!40} m m}{%
\draw[mycolor=#1] (#2.north west)rectangle (#3.south east);
}

\NewDocumentCommand{\fhighlight}{O{blue!40} m m}{%
\draw[myfillcolor=#1] (#2.north west)rectangle (#3.south east);
}

\newtheorem{remark}{Remark}[section]

\def\bU{{\mathbf U}}
\def\bSigma{{\boldsymbol \Sigma}}
\def\bu{{\mathbf u}}
\def\bp{{\mathbf p}}
\def\bq{{\mathbf q}}
\def\bv{{\mathbf v}}
\def\bx{{\mathbf x}}
\def\by{{\mathbf y}}
\def\bb{{\mathbf b}}
\def\bff{{\mathbf f}}
\def\bsigma{{\boldsymbol \sigma}}

\def\bone{{\boldsymbol 1}}
\def\bzero{{\boldsymbol 0}}

\newcommand{\vertexspace}{\bU}
\newcommand{\edgespace}{\bSigma}

\newcommand{\graph}{G}
\newcommand{\vertices}{V}
\newcommand{\edges}{E}
\newcommand{\agg}{{\mathcal A}}
\newcommand{\bgg}{{\mathcal B}}

\newcommand{\face}{{\mathcal F}}

\newcommand{\Lapl}{{\mathcal L}}

\newcommand{\Span}[1]{\text{Span}\left\{\; #1 \;\right\}}
\newcommand{\Null}[1]{\text{Null}\left( #1 \right)}
\renewcommand{\Range}[1]{\text{Range}\left( #1 \right)}

\newcommand{\bdr}{\partial}

\newcommand{\corr}{\gamma} 

\graphicspath{}

\title{Multilevel Spectral Coarsening for Graph Laplacian Problems with Application to Reservoir
Simulation\thanks{Performed under the auspices of the U.S. Department of Energy under Contract DE-AC52-07NA27344 (LLNL-JRNL-795643).}}

\author{Andrew~T.~Barker\thanks{Center for Applied Scientific Computing, Lawrence Livermore National Laboratory, Livermore, CA 94550
({\tt atb@llnl.gov, cslee@llnl.gov, osborn9@llnl.gov, panayot@llnl.gov})} \and Stephan~V.~Gelever\thanks{The Fariborz Maseeh Department of Mathematics and Statistics, Portland State University, Portland, OR 97201
({\tt stephan.gelever@gmail.com, panayot@pdx.edu})} \and Chak~S.~Lee\protect\footnotemark[2] \and Sarah~V.~Osborn\protect\footnotemark[2] \and Panayot~S.~Vassilevski\footnotemark[2]\hspace{1.2mm}\textsuperscript{,}\footnotemark[3]}

\begin{document}
\maketitle

\begin{abstract}
We extend previously developed two-level coarsening procedures for graph Laplacian problems written in a mixed saddle point form to the fully recursive multilevel case. The resulting hierarchy of
discretizations gives rise to a hierarchy of upscaled models, in the sense that they provide approximation in the natural norms (in the mixed setting). This property enables us to utilize them in
three applications: (i) as an accurate reduced model, (ii) as a tool in multilevel Monte Carlo simulations (in application to finite volume discretizations), and (iii) for providing a sequence of nonlinear operators in FAS
(full approximation scheme) for solving nonlinear pressure equations discretized by the conservative
two-point flux approximation. We illustrate the potential of the proposed multilevel technique in all
three applications on a number of popular benchmark problems used in reservoir simulation.
\end{abstract}


\begin{AMS}{65N55, 65N08, 65F15}\end{AMS}

\section{Introduction}
The massive amount of data available today makes it typically impractical to take into account all of the data for decision making, so that it is becoming increasingly important to extract useful information from raw data effectively and efficiently.
For example, in reservoir simulations, geological models are often very detailed so that direct simulations on these fine scale models are prohibitively expensive. Traditionally, rock properties (e.g.~porosity, permeability) in fine scale geological models are upscaled or homogenized to a coarser resolution that is suitable for simulations \cite{wg96, wdle00,hornung97,cps07}. Other than upscaling and homogenization, several sophisticated methods for constructing coarse models were proposed in the past two decades, including (generalized) multiscale finite element methods \cite{hw97, EfendievHouBook, egh12, ArPeWY07}, multiscale finite volume methods \cite{jennylt03, jlt06, tjl07}, heterogeneous multiscale methods \cite{ee2003, aeev2012}, and (algebraic) multigrid \cite{mm06, lv08, vassilevski-upscaling, kalchev-lee-upscaling-mixed, blv17}. Among these methodologies, generalized multiscale finite element methods \cite{egh12} and spectral AMGe \cite{cfhjmmrv03, kalchev-lee-upscaling-mixed, blv17} share the idea of constructing coarse spaces via local spectral problems, which allows the user to adjust approximation quality of the coarse space by altering the number of local basis functions (without changing the coarse grid). For exceedingly large problems, multilevel coarsening is desirable to provide the user a range of accuracy levels. In addition, for multilevel Monte Carlo simulations and multilevel solvers for nonlinear problems, nested hierarchies of various resolutions are indispensable. To produce nested hierarchical models, multigrid methods are the obvious choice due to their intrinsic multilevel nature.

In \cite{blv17}, a finite volume discretization using two-point flux approximation (TPFA) for Darcy's flow problems is interpreted as a graph Laplacian, and coarse reservoir models are constructed by algebraically coarsening the Laplacian. The input required by this algebraic method is essentially the linear system only, which makes the method very easy to apply. Moreover, it enables flexible local enrichment of the coarse approximation space, and hence produces coarse models with a range of accuracies. In this paper, our goal is to extend the two-level method developed in \cite{blv17} to a multilevel method, and demonstrate the applications of the method in numerical upscaling, multilevel Monte Carlo simulations (of finite volume problems), and nonlinear multigrid solvers (full approximation scheme). It is shown in \cite{blv17} that the coarse spaces have approximation properties in the natural norms (in mixed setting), which is the key for the targeted applications to work effectively.

The coarse system produced by the two-level spectral coarsening method in \cite{blv17} is not a graph Laplacian, so it seems that the two-level coarsening method cannot be applied directly to coarsen the coarse system further. But if we associate degrees of freedom to vertices/edges of the graph, and think of the coarsening with respect to degrees of freedom rather than with respect to the graph topology, then the coarsening method can actually be applied recursively in a natural way. The main idea is to aggregate degrees of freedom based on their associated vertices/edges of the graph so that the local problems are well-defined. In the finest level, there is a one-to-one correspondence between the vertex/edge degrees of freedom and vertices/edges of the original graph; in coarse levels, the relation between coarse degrees of freedom and vertices/edges of the coarse graph is available from the coarsening procedure. Hence, the required input of the overall multilevel coarsening procedure is still the original fine scale system only.

Since our spectral coarsening method involves solving local eigen and extension problems, the setup of the hierarchy is quite expensive (compared with the classic algebraic multigrid \cite{classicAMG}). Therefore, in multilevel Monte Carlo simulations and FAS, it is crucial that the coarse spaces can be reused even when the problem parameters (e.g. permeability coefficient, source term) are changed. This is achieved by first storing the local coarse matrices, and when the coefficient is updated, the local matrices are scaled and then assembled into the final global system. Another issue is that the coarse systems generated by our coarsening method are of saddle point form, even though the original problem (graph Laplacian) is symmetric positive-definite. Consequently, in order to have faster simulation on coarse levels, we need an efficient solver for the coarse saddle point problems. To this end, we again make use of the availability of local matrices (during the coarsening steps) and transform the saddle point problems into symmetric positive-definite problems via hybridization, which is a classic technique for solving discrete problems arising from mixed finite element discretization of partial differential equations \cite{ArnoldBrezzi85}. Recently in \cite{lv16}, and \cite{dkltv19}, an algebraic multigrid with some proper diagonal rescaling is demonstrated be a competitive solver for the hybridized system. Because the solver is constructed algebraically, it can be applied to solve saddle point problems that are not coming from finite element discretizations, which is the case in the current paper.

The rest of the paper is organized as follows. In Section~\ref{sec:notation}, we state the problem of interest and notation. Then the multilevel spectral coarsening algorithm will be presented in detail in Section~\ref{sec:algorithm}. Lastly, we illustrate some applications of the coarsening algorithm in Section~\ref{sec:numerics} and conclude the paper in Section~\ref{sec:conclusion}.

\section{Some setup and notation}
\label{sec:notation}

Throughout the paper, we use the convention that level $\ell = 0$ corresponds to the finest resolution while $\ell = L$ corresponds to the coarsest resolution. The number of elements in a set $S$ will be denoted by $\#S$. For any matrix $A$, its null space and range space are denoted by $\Null{A}$ and $\Range{A}$ respectively.

We have a connected and undirected (fine) graph $\graph^0 = (\vertices^0, \edges^0)$, where $\vertices^0$ and $\edges^0$ are the set of vertices and edges respectively. Let $\Lapl^0$ be a weighted graph Laplacian associated with $G^0$. It is well-known that $\Lapl^0$ can be decomposed as $\Lapl^0 = D^0 (M^0)^{-1} (D^0)^T$, where $M^0$ is a diagonal matrix containing the inverse of edge weights and $(D^0)^T$ resembles some discrete gradient. More precisely, each row of $(D^0)^T$, which corresponds to some edge connecting vertex $v_i^0$ to $v_j^0$, has exactly two nonzeros, 1 at the corresponding column of $v_i^0$ and -1 at the other. By considering this decomposition, the graph Laplacian problem
\[
\Lapl^0 \bu^0 = \bff^0
\]
is equivalent to the mixed graph Laplacian problem
\[
  \begin{bmatrix} M^0& (D^0)^T \\ D^0& \end{bmatrix}  \begin{bmatrix} \bsigma^0 \\ \bu^0  \end{bmatrix}  =  \begin{bmatrix} 0 \\ -\bff^0 \end{bmatrix}
\]
where $\bu^0$ is in the vertex space $\vertexspace^0 := \mathbb{R}^{\#\vertices^0}$ and $\bsigma^0$ is in the edge space $\edgespace^0 := \mathbb{R}^{\#\edges^0}$. The coarsening scheme in this paper is based on the mixed setting.

Our main objective is to construct two sequences of matrices $\{ M^\ell \}_{\ell=0}^L$ and $\{ D^\ell \}_{\ell=0}^L$ for defining coarse mixed graph Laplacians. This is accomplished by constructing interpolation matrices $\{ P_\sigma^\ell \}_{\ell=0}^{L-1}$ and $\{ P_u^\ell \}_{\ell=0}^{L-1}$ and applying the recursive definition:
\begin{equation}
\begin{bmatrix}
M^{\ell+1} & (D^{\ell+1})^T\, \\ D^{\ell+1}
\end{bmatrix}
= \begin{bmatrix}
P_\sigma^\ell \\ & P_u^\ell
\end{bmatrix}^T \begin{bmatrix}
M^\ell & (D^\ell)^T\, \\ D^\ell
\end{bmatrix}
\begin{bmatrix}
P_\sigma^\ell \\ & P_u^\ell
\end{bmatrix}, \quad \text{for }\ell = 0, 1, \dots, L-1.
\label{eq:RAP}
\end{equation}

For multilevel Monte Carlo simulations and nonlinear multigrid solvers, we will also need projection operators that take quantities from a given level to a coarser level. To this end, we will construct $\{ Q_\sigma^\ell \}_{\ell=0}^{L-1}$ and $\{ Q_u^\ell \}_{\ell=0}^{L-1}$ such that $Q_u^\ell$ is a left inverse of $P_u^\ell$ and $Q_\sigma^\ell$ is a left inverse of $P_\sigma^\ell$.

\section{Multilevel spectral coarsening of graph Laplacians}
\label{sec:algorithm}

In \cite{blv17}, a two-level spectral coarsening algorithm for the graph Laplacian $D^0 (M^0)^{-1} (D^0)^T$ is introduced. The algorithm constructs interpolation matrices $P_u^0$ and $P_\sigma^0$ for vertex and edge spaces respectively. Then the coarse level operators $M^1$ and $D^1$ are defined through the ``RAP'' procedure \eqref{eq:RAP} when $\ell = 0$. Note that the coarse operator $D^1 (M^1)^{-1} (D^1)^T$ is not a graph Laplacian in the usual sense (for example, its row sums are not necessarily zero). Therefore, the coarsening procedure in \cite{blv17} cannot be applied directly to $D^1 (M^1)^{-1} (D^1)^T$. To extend the two-level algorithm to a recursive multilevel one, the key idea is to make a distinction between graph entities (vertices and edges) and their associated vector spaces. When defining local problems, we consider aggregation of degrees of freedom (dofs) instead of entities. Before we dive into the algorithm in detail, we shall give a high-level view of the coarsening procedure:
\begin{enumerate}[leftmargin=10mm]
  \item Coarsen graph entities: aggregate vertices and edges.
  \item Aggregate vertex dofs and edge dofs based on aggregation of vertices and edges.
  \item Solve spectral problems on aggregates of fine dofs to define coarse vertex dofs (for good approximation).
  \item Define coarse edge dofs associated with edges connecting aggregates (for coupling of coarse dofs).
  \item Define coarse edge dofs associated with aggregates of vertices (for stability of the coarse problem). \vspace{3mm}
\end{enumerate}

\begin{remark}
We stress here that, when $L = 1$, the multilevel algorithm in the current paper will result in exactly the two-level method proposed in \cite{blv17}.
\end{remark}

\subsection{Graph hierarchy}
\label{sec:hierarchy}

\begin{figure}[h]
\centering
\begin{tikzpicture}[scale=.35, auto]
\draw [fill=blue,blue] (1.8,5.8) circle [radius=5.8];
\draw [fill=blue,blue] (17,4.3) circle [radius=5];
\draw [line width=1.2cm,blue] (3,4.8)--(15,3.8);
\draw [fill=red,red] (1.7,1.5) circle [radius=1.2];
\draw [fill=red,red] (15,2.1) circle [radius=1.6];
\draw [fill=red,red] (5.3,4.2) circle [radius=1.2];
\draw [fill=red,red] (3,9.5) circle [radius=1.4];
\draw [fill=red,red] (16,5.8) circle [radius=1.4];
\draw [line width=0.2cm,red] (1.7,1.8)--(5.3,4.2);
\draw [line width=0.2cm,red] (15,3)--(6,4);
\draw [line width=0.38cm,red] (3.1,9)--(5.2,4.2);
\draw [line width=0.2cm,red] (16,5)--(5.5,4.5);
\draw [line width=0.18cm,red] (16,5)--(15,3);
\draw (5,4)--(5.3,4.5)--(6,4)--(5,4);
\foreach \pt in {(1,1), (2.5,1), (1.2, 2.2), (5,4)}
\draw \pt--(2,2);
\foreach \pt in {(16,2), (14.5,1), (15.5,3.3), (6,4),(16,3)}
\draw \pt--(15,3);
\foreach \pt in {(2,10), (2.7,8.5), (3.3,9.2), (5,4), (5.3,4.5)}
\draw \pt--(2.9,9);
\foreach \pt in {(3.1,10), (3.5,10.5), (5.3,4.5)}
\draw \pt--(3.3,9.2);
\foreach \pt in {(17,6), (15,5.5), (5.3,4.5), (15,3)}
\draw \pt--(16,5);
\foreach \pt in {(1,1),(2,2),(2.5,1),(1.2, 2.2),(5,4),(5.3,4.5),(6,4),(16,2),(14.5,1),(15.5, 3.3),(16,3),(15,3)}
\draw [fill] \pt circle [radius=0.13];
\foreach \pt in {(3.1,10),(3.5,10.5),(2.9,9),(2,10),(2.7,8.5),(3.3,9.2),(17,6),(15,5.5),(16,5)}
\draw [fill] \pt circle [radius=0.13];
\end{tikzpicture}
\caption{An example of a 3-level graph hierarchy. Black: level 0, red: level 1, blue: level 2.}
\label{fig:hierarchy}
\end{figure}
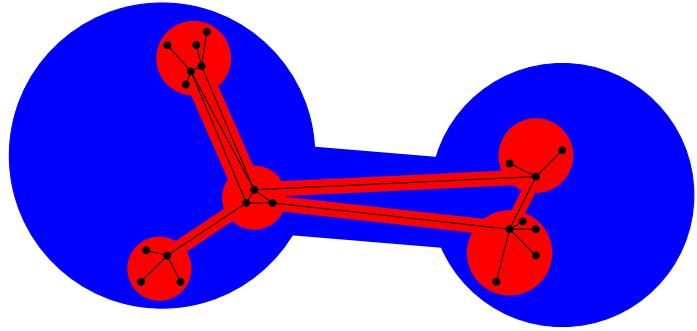

Starting with the original graph $\graph^0$, we will generate a hierarchy of graphs $\{G^\ell\}_{\ell=0}^L$ with $L+1$ levels of resolution. Given a graph $\graph^\ell = (V^{\ell}, E^{\ell})$ on level $\ell$, we construct a coarser graph $\graph^{\ell+1}$ as follows. First, we partition the set of vertices $\vertices^{\ell}$ into non-overlapping aggregates of vertices $\{\agg_j^\ell\}$ such that the subgraph within each aggregate is connected. We identify the set of aggregates on level $\ell$ as the set of vertices on level $\ell+1$, that is,
\[
V^{\ell+1} = \left\{ v^{\ell+1}_j \right\} := \left\{ \agg^\ell_j \right\}.
\]
Aggregates inherit their connections from the finer level. In particular we say that $ \agg^\ell_i $ is {\em connected} to $ \agg^\ell_j $ if there is some edge $ e^\ell = (v^\ell_i, v^\ell_j) $ of the finer graph $G^\ell$ such that $ v^\ell_i \in \agg^\ell_i $ and $ v^\ell_j \in \agg^\ell_j $, and if such a connection exists we define a {\em face} as the collection of all such edges
\[
\face^\ell := \left\{\; e^{\ell} = (v_i^{\ell}, v_j^{\ell}) \in \edges^{\ell} :\, v_i^{\ell} \in \agg^\ell_i \text{ and } v_j^{\ell} \in \agg^\ell_j \;\right\}.
\]
Note that the faces on level $\ell$ can naturally be identified as the set of edges on level $\ell+1$
\[
E^{\ell+1} :=  \left\{ e^{\ell+1} := \face^\ell  = ( \agg^\ell_i, \agg^\ell_j ) = ( v_i^{\ell+1}, v_j^{\ell+1} ) \right\}.
\]
This completes the construction of a coarser graph $\graph^{\ell+1} = (V^{\ell+1}, E^{\ell+1})$ from $\graph^\ell$. This process can easily be applied recursively, resulting in a hierarchy of graphs, as illustrated in Figure~\ref{fig:hierarchy}.

\subsection{Graph spaces}
\label{sec:space}

While the construction in Section~\ref{sec:hierarchy} preserves a graph structure on each level, the spectral coarsening method in \cite{blv17} does {\em not} preserve a graph Laplacian structure. The main problem is that the spectral method allows more than one dof per entity in the coarse vertex and edge spaces. Consequently, the matrix graph of the coarse Laplacian may have some pair of vertices that are connected through 3 or more edges. In other words, we cannot extract a normal undirected graph using the coarse Laplacian operator alone. To remedy this situation, we need to independently keep track of the aggregation of graph entities (as in Section~\ref{sec:hierarchy}), and their associated dofs. To this end we introduce local graph spaces, whose definitions are motivated by the coarse spaces of the two-level algorithm: 
\begin{itemize}
\item $\vertexspace^\ell(v^\ell_i)$: span of level-$\ell$ vertex-based dofs associated with vertex $v^\ell_i$.
\item $\edgespace^\ell(e^\ell_i)$: span of level-$\ell$ edge-based dofs associated with edge $e^\ell_i$.
\item $\edgespace^\ell(v^\ell_i)$: span of level-$\ell$ edge-based dofs associated with vertex $v^\ell_i$.
\end{itemize}

In our local spectral coarsening algorithm, we will need to extract submatrices corresponding to dofs associated with certain aggregates of vertices/edges, so it will be useful to define local graph spaces on some arbitrary aggregate of vertices $\bgg$ and face $\face$ (aggregate of edges): 
\begin{equation*}
\begin{split}
\vertexspace^\ell(\bgg) := \bigoplus_{v^\ell_i\in\bgg}\vertexspace^\ell(v^\ell_i), \qquad &
\edgespace^\ell(\bgg) := \left( \bigoplus_{e^\ell_i \in E^\ell(\bgg)}\edgespace^\ell(e^\ell_i) \right) \bigoplus \left( \bigoplus_{v^\ell_i \in \bgg}\edgespace^\ell(v^\ell_i) \right), \\
\edgespace^\ell(\face) := \bigoplus_{e^\ell_i\in\face}\edgespace^\ell(e^\ell_i), \qquad &
\bdr\edgespace^\ell(\bgg) := \bigoplus_{e^\ell_i\in\bdr E^\ell(\bgg)}\edgespace^\ell(e^\ell_i),
\end{split}
\end{equation*}
where
\[
E^\ell(\bgg) := \left\{ \, e = (v_i, v_j) \in E^\ell : v_i, v_j \in\bgg \, \right\} \;\; \text{ and }  \;\;\; \bdr E^\ell(\bgg) := \left\{ \, e = (v_i, v_j) \in E^\ell : v_i \in \bgg, v_j \not\in\bgg \, \right\}.
\]
The global spaces are simply when $\bgg$ is $V^\ell$, the set of all vertices in level $\ell$:
\[
\vertexspace^\ell := \vertexspace^\ell(V^\ell), \qquad 
\edgespace^\ell := \edgespace^\ell(V^\ell).
\]

On the finest level ($\ell = 0$), each of $\vertexspace^0(v^0_i)$ and $\edgespace^0(e^0_i)$ is just taken to be the space of real numbers $\mathbb{R}$, and $\edgespace^0(v^0_i)$ is empty. Hence, the global spaces match with our earlier definition $\vertexspace^0 = \bigoplus_{v^0_i\in V^0}\vertexspace^0(v^0_i) = \mathbb{R}^{\#\vertices^0}$, and $\edgespace^0 = \bigoplus_{e^0_i\in E^0}\edgespace^0(e^0_i) = \mathbb{R}^{\#\edges^0}$. On coarser levels ($\ell > 0$), $\vertexspace^\ell(v^\ell_i)$ contains eigenvectors of some local spectral problem associated with $v^\ell_i = \agg^{\ell-1}_i$. $\edgespace^\ell(e^\ell_i)$ is spanned by the trace extensions associated with $e^\ell_i = \face^{\ell-1}_i$, and $\edgespace^\ell(v^\ell_i)$ is spanned by the bubbles associated with $v^\ell_i = \agg^{\ell-1}_i$. The precise definition and detailed construction of all these components in the coarse levels will be given in Section~\ref{sec:vdof}--~\ref{sec:edof}.

Now, for an aggregate of vertices $\bgg$ and a face $\face$, we define its associated submatrices:
\begin{equation}
\begin{split}
M^\ell_\bgg & := M^\ell_{\edgespace^\ell(\bgg), \edgespace^\ell(\bgg)},  \quad
M^\ell_{\partial \bgg} := M^\ell_{\edgespace^\ell(\bgg), \partial \edgespace^\ell(\bgg)}, \quad
M^\ell_{\bgg, \face} := M^\ell_{\edgespace^\ell(\bgg), \edgespace^\ell(\face)}, \\
D^\ell_\bgg & := D^\ell_{\vertexspace^\ell(\bgg), \edgespace^\ell(\bgg)},  \quad
D^\ell_{\partial \bgg} := D^\ell_{\vertexspace^\ell(\bgg), \partial \edgespace^\ell(\bgg)}, \quad
D^\ell_{\bgg, \face} := D^\ell_{\vertexspace^\ell(\bgg), \edgespace^\ell(\face)}. \\
\end{split}
\label{eq:submatrices}
\end{equation}
Here, $M^\ell_{\edgespace^\ell(\bgg), \edgespace^\ell(\face)}$ is the submatrix of $M^\ell$ whose rows and column are respectively restricted to $\edgespace^\ell(\bgg)$ and $\edgespace^\ell(\face)$; other submatrices are defined similarly. An analogous notion for subvectors associated with aggregates and faces will also be used.

\begin{remark}
The concept of graph spaces is an analogue to finite element spaces in the context of discretizations of partial differential equations, in which each dof is associated with some geometrical entity (e.g. element, face) of a mesh.
The present graph-based construction can be seen as a version of the de Rham sequence corresponding only to $H(\text{div})$ and $L_2$ degrees of freedom on an agglomerated mesh. See \cite{lv08} for detailed descriptions of such mesh-based constructions.
\end{remark}

\subsection{Vertex degrees of freedom}
\label{sec:vdof}
Following \cite{blv17}, when generating coarse vertex degrees of freedom, we consider the local matrices associated with the neighborhood of each aggregate $N(\agg^\ell_i)$. That is, in this case $\bgg$ in \eqref{eq:submatrices} is taken to be $\bgg = N(\agg^\ell_i)$, where $N(\agg^\ell_i)$ is defined to be the set of all vertices connected by an edge to a vertex in $ \agg^\ell_i $. Define
\[
  \Lapl_{N(\agg^\ell_i)} := D^{\ell}_{N(\agg^\ell_i)} \big(M^{\ell}_{N(\agg^\ell_i)}\big)^{-1} \big(D^{\ell}_{N(\agg^\ell_i)}\big)^T
\]
and solve the eigenvalue problem
\begin{equation}
  \Lapl_{N(\agg^\ell_i)}  \bq = \lambda \bq.
  \label{eq:eigen}
\end{equation}
We select a few eigenvectors corresponding to small eigenvalues, restrict them to $\agg^\ell_i$, use SVD to remove possible linear dependence, and choose the remaining linearly independent vectors $\tilde{\bq}$ as columns of the matrix $P_{u, \agg_i^\ell}^{\ell}$. By restricting local basis functions to aggregates $\agg^\ell_i$, the support of basis functions associated with different aggregates will have no overlap. The final interpolation matrix $P_u^\ell$ for the vertex space is obtained by diagonally concatenating $P_{u, \agg_i^\ell}^{\ell}$ for all aggregates $\agg_i^\ell\in\vertices^{\ell+1}$:
\[
P_u^\ell := \begin{bmatrix}
\ddots & & \\ & P_{u, \agg_i^\ell}^{\ell} & \\ & & \ddots
\end{bmatrix}.
\]
The local coarse vertex space of each aggregate $\agg^\ell_i$ (vertex in level $\ell+1$) is the column space of $P_{u, \agg_i^\ell}^{\ell}$: 
\[
\vertexspace^{\ell+1}(v^{\ell+1}_i) = \vertexspace^{\ell+1}(\agg^\ell_i) := \Range{P_{u, \agg_i^\ell}^{\ell}}.
\]

Since $P_u^\ell$ is block diagonal and the columns of each block are orthonormal, we simply define the projection operator for the coarse vertex space as $\pi_u^\ell : = P_u^\ell Q_u^\ell$, where
\[
Q_u^\ell := \big( P_u^\ell \big)^T.
\]
Clearly, $Q_u^\ell$ is defined locally in each aggregate and
\begin{equation}
Q_u^\ell P_u^\ell = I_u^{\ell+1}.
\label{eq:P_u_orthogonal}
\end{equation}

\begin{remark}
\label{rmk:const_representation}
Let $\bone^0$ be the constant vector of 1's in $\vertexspace^0 = \mathbb{R}^{\#\vertices^0}$ (i.e., $\bone^0 := [1, 1, \cdots, 1]^T$). Then, its projection in coarse levels are defined recursively as
\begin{equation}
\bone^{\ell} :=  Q_u^{\ell-1} \,\bone^{\ell-1}, \quad \text{for }\ell = 1, \dots, L.
\label{eq:const_representation}
\end{equation}

\noindent From Lemma~\ref{lemma:D_null} and Lemma~\ref{lemma:const_preserve} in the appendix, we know our construction guarantees that $\bone_{\agg_i^\ell}^\ell$ (the restriction of $\bone^\ell$ to $\agg^\ell_i$) is always in $\Range{P_{u, \agg_i^\ell}^{\ell}}$. Therefore, $P_{u, \agg_i^\ell}^\ell$ is of the form
\begin{equation}
P_{u, \agg_i^\ell}^\ell = {\footnotesize
\begin{tikzpicture}[baseline=-\the\dimexpr\fontdimen22\textfont2\relax ]
\matrix (m)[matrix of math nodes,left delimiter={[},right delimiter={]}]
{
\bq_{\agg_i^\ell}^{PV} & \hspace{1mm} & P_{u, \agg_i^\ell}^{\ell,NPV}  \\
}; 
\end{tikzpicture} } , 
\label{eq:vertex_basis_decomposition}
\end{equation}
where columns of $P_{u, \agg_i^\ell}^{\ell,NPV}$ are orthogonal to $\bone_{\agg_i^\ell}^\ell$ and $\bq_{\agg_i^\ell}^{PV}$ is the unit vector parallel to $\bone_{\agg_i^\ell}^\ell$:
\begin{equation}
\bq_{\agg_i^\ell}^{PV} = \left(\|\bone_{\agg_i^\ell}^\ell\|^{-1}\right) \, \bone_{\agg_i^\ell}^\ell.
\label{eq:vertex_pv_vector}
\end{equation}
Note that we use the superscript ``$PV$" because $\bq_{\agg^\ell}^{PV}$ are the coarse basis in an earlier work of Pasciak and Vassilevski \cite{pasciak-vassilevski}. Together with ``PV" traces that will be defined in the next section, they are the building blocks of a pair of minimal coarse spaces such that the coarse problem is inf-sup stable, cf. condition \eqref{eq:pv_condition}.

\end{remark}

\begin{figure}[h!]
\centering
\subfigure[Basis 1 in $\agg_1^1$]{\includegraphics[scale=.165,clip,trim=.8cm 0 .8cm 0]{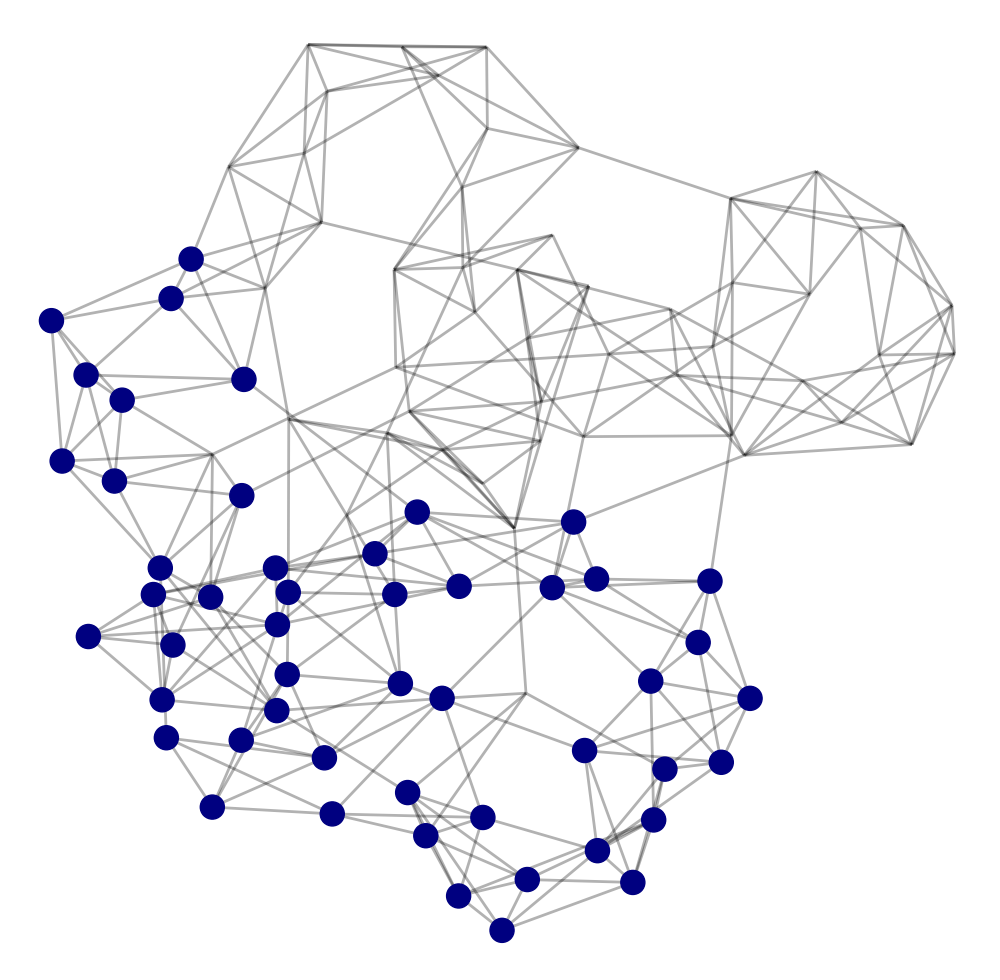}}
\subfigure[Basis 2 in $\agg_1^1$]{\includegraphics[scale=.165,clip,trim=.8cm 0 .8cm 0]{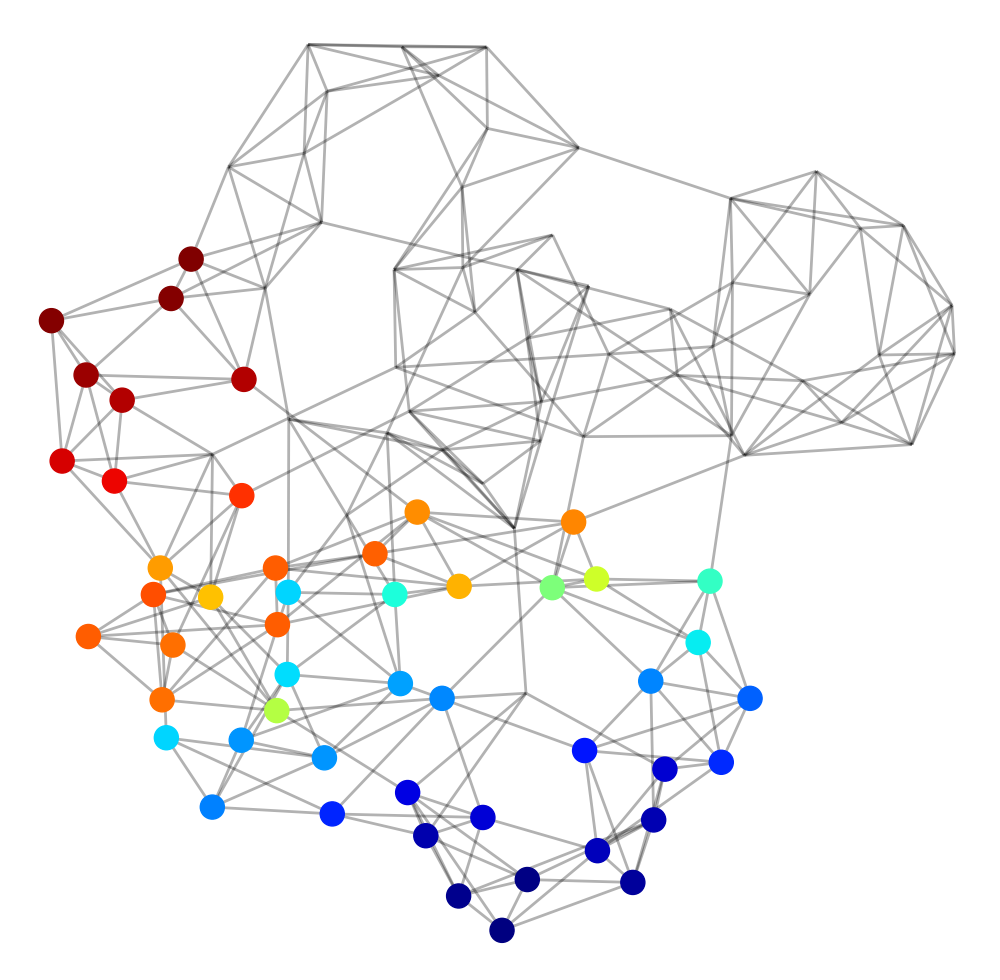}}
\subfigure[Basis 1 in $\agg_2^1$]{\includegraphics[scale=.165,clip,trim=.8cm 0 .8cm 0]{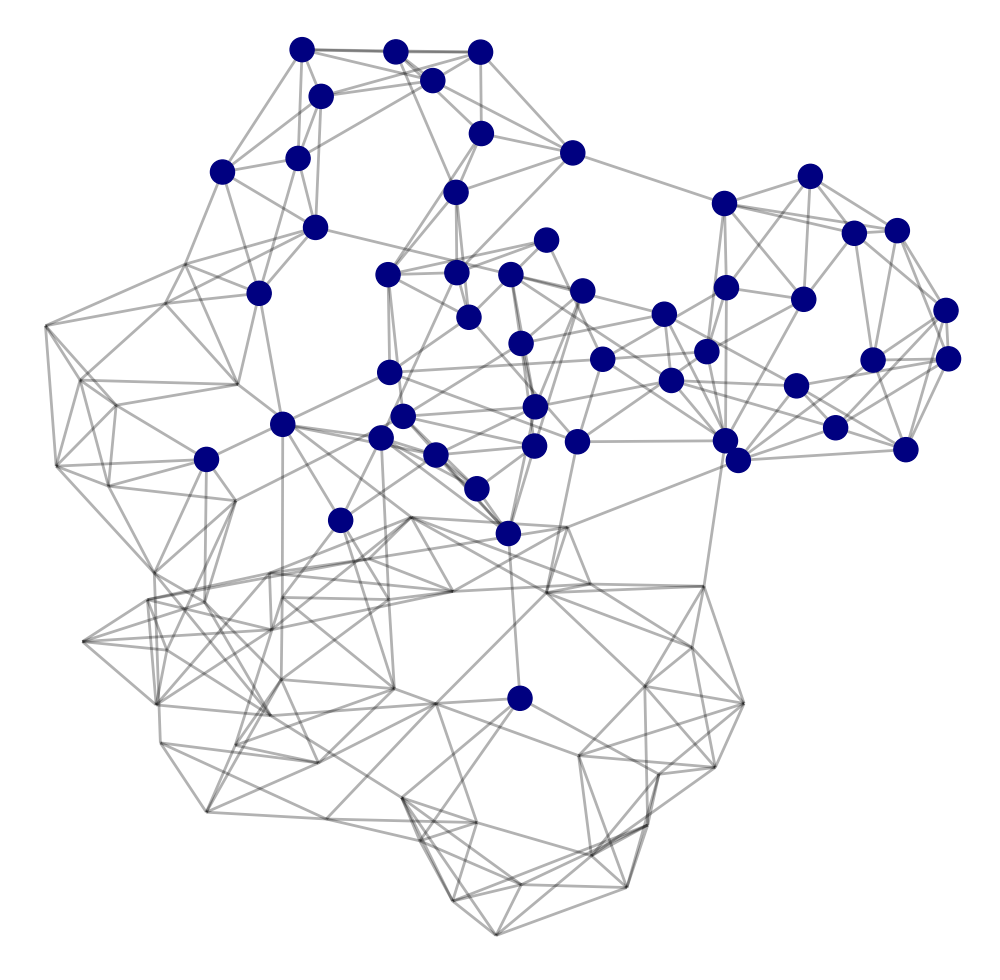}}
\subfigure[Basis 2 in $\agg_2^1$]{\includegraphics[scale=.165,clip,trim=.8cm 0 .8cm 0]{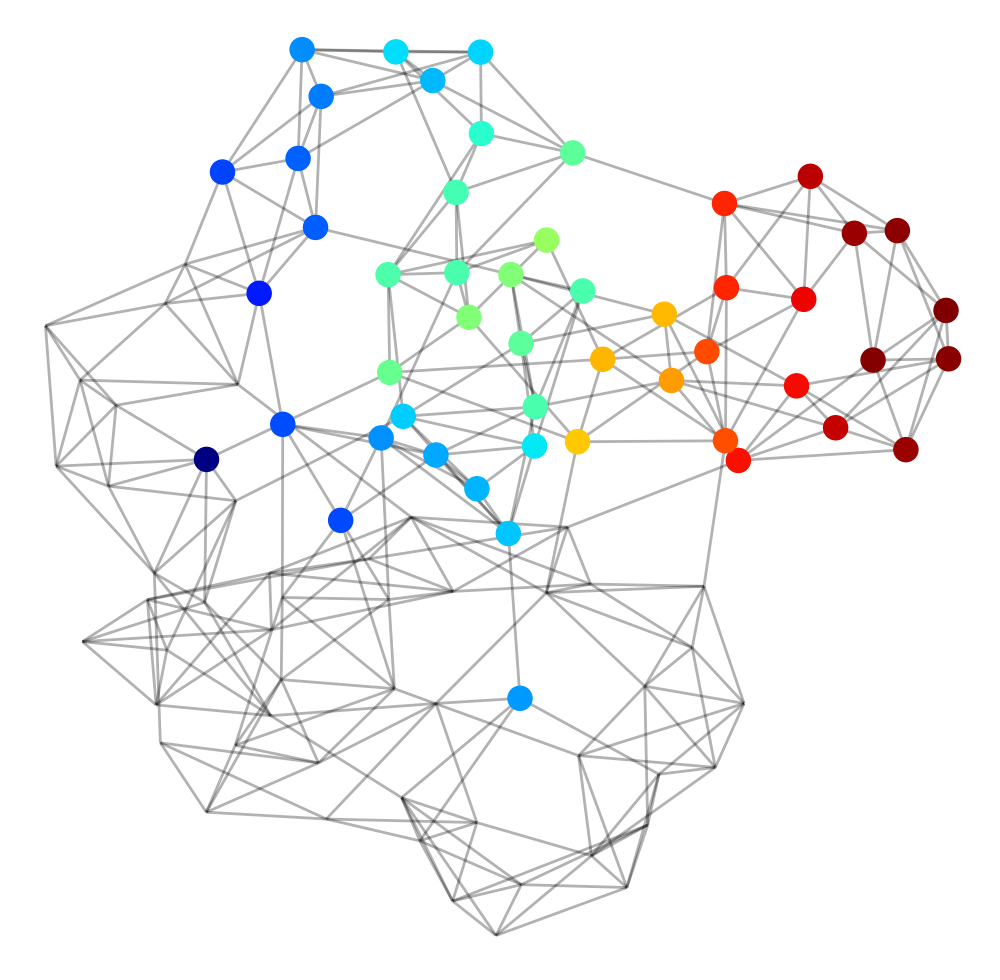}}
\caption{Examples of level 1 coarse vertex degrees of freedom of some random graph}
\label{fig:basis_level1}
\end{figure}

\begin{figure}[h!]
\centering
\subfigure[Basis 1 in $\agg_1^2$]{\includegraphics[scale=.165,clip,trim=.8cm 0 .8cm 0]{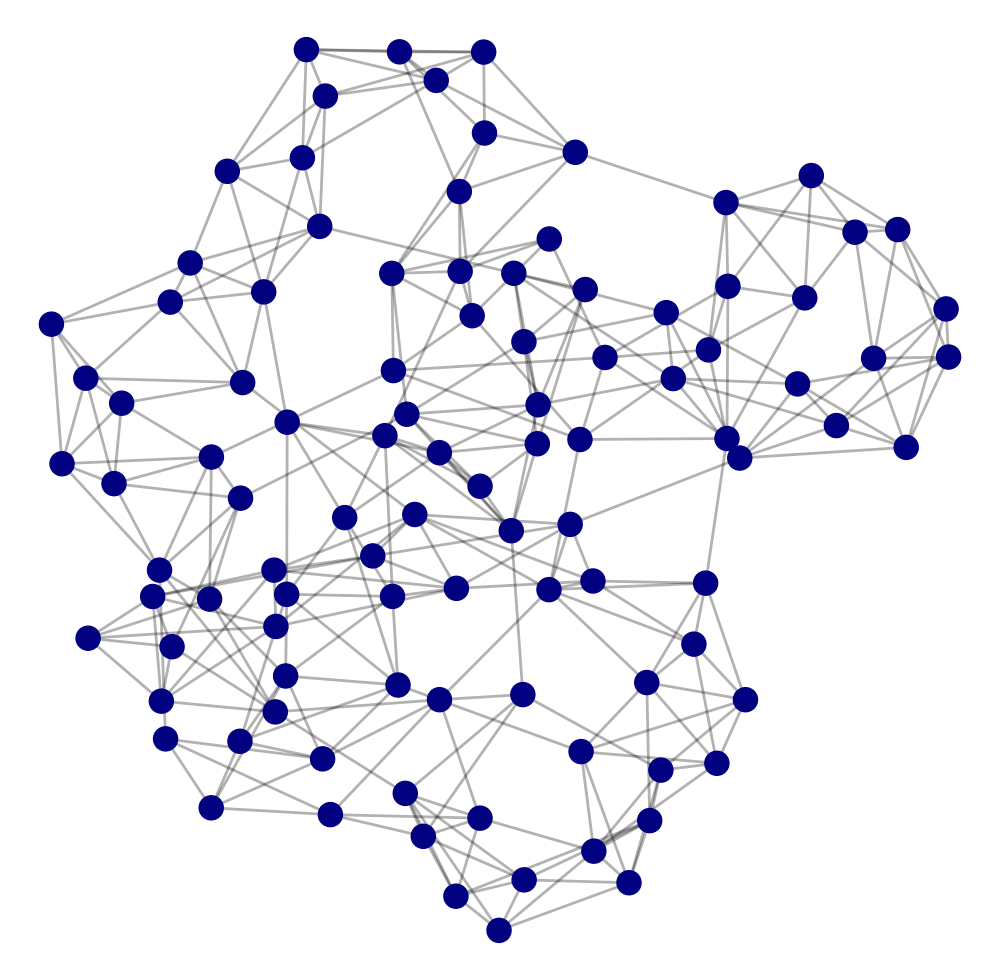}}
\subfigure[Basis 2 in $\agg_1^2$]{\includegraphics[scale=.165,clip,trim=.8cm 0 .8cm 0]{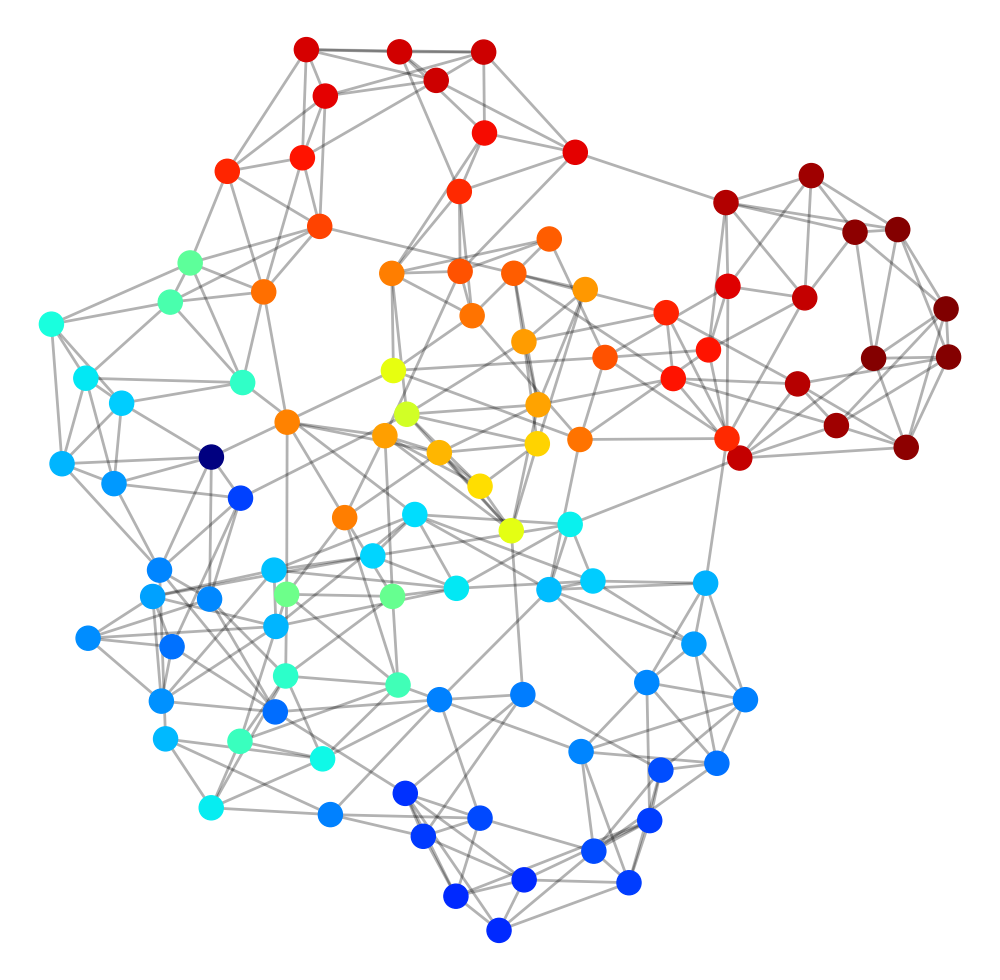}}
\subfigure[Basis 3 in $\agg_1^2$]{\includegraphics[scale=.165,clip,trim=.8cm 0 .8cm 0]{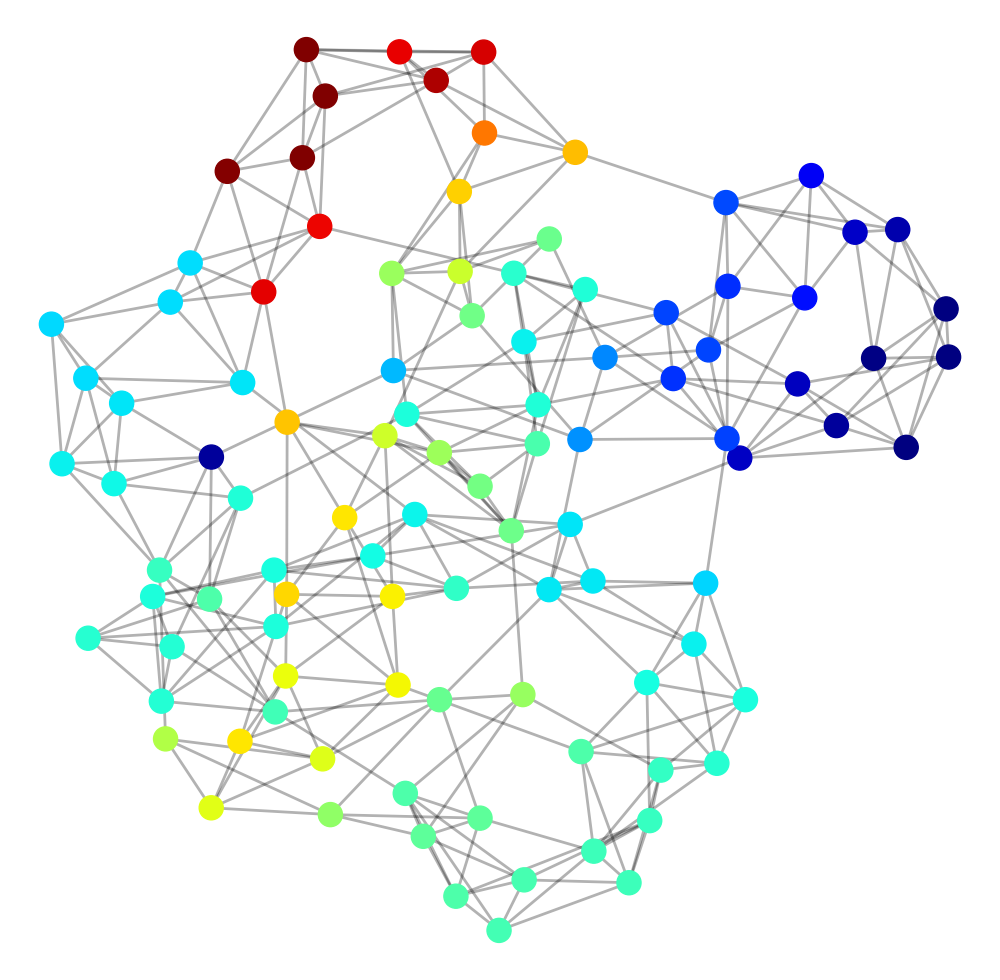}}
\subfigure[Basis 4 in $\agg_1^2$]{\includegraphics[scale=.165,clip,trim=.8cm 0 .8cm 0]{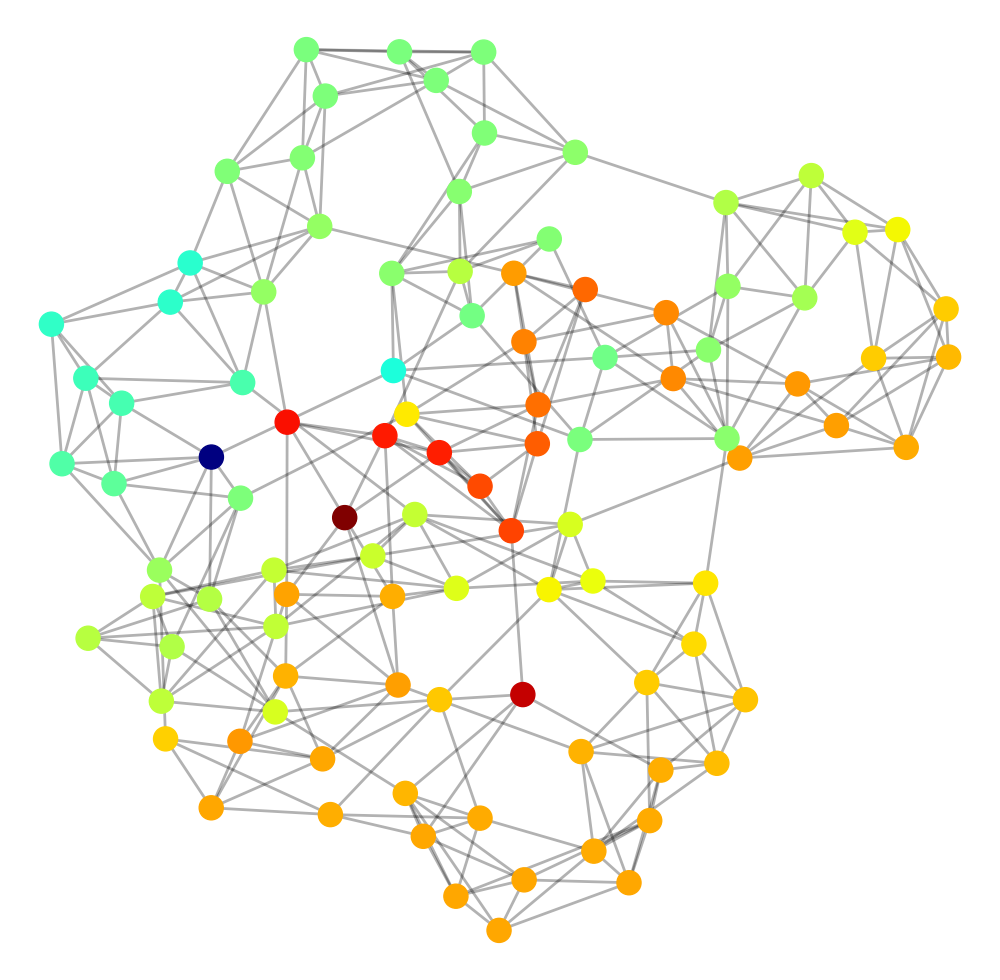}}
\caption{Examples of level 2 coarse vertex degrees of freedom of some random graph}
\label{fig:basis_level2}
\end{figure}

Figures~\ref{fig:basis_level1} and~\ref{fig:basis_level2} illustrate some example coarse vertex degrees of freedom in a 3-level hierarchy of some random graph. Notice that the first basis of any aggregate on any level is always constant.

\subsection{Edge degrees of freedom}
\label{sec:edof}

There are two types of coarse edge degrees of freedom, one is called trace extensions while the other one is called bubbles. Each trace extension is associated with a face, and is supported in the aggregates sharing the face. They are the only components of the coarse space that connect neighboring aggregates and hence couple the global problem. On the other hand, each bubble is associated with an aggregate. The support of a bubble is in one and only one aggregate, hence the name bubble.

\subsubsection{Trace extensions}
\label{sec:trace}

The idea here is to define some representative vectors (which are what we call traces) that live on faces, and then extend them harmonically into the neighboring aggregates. It suffices to focus on one face because trace extensions are associated with faces. Recall that by definition, each face is associated with two aggregates: $ \face^\ell = (\agg^\ell_i, \agg^\ell_j) $. For each associated aggregate and each eigenvector $\bq$ selected when solving \eqref{eq:eigen}, we define the trace
\[
  \left. \left( \big(M^{\ell}_{N(\agg^\ell)}\big)^{-1} \big(D^{\ell}_{N(\agg^\ell)}\big)^T \bq \right) \right|_{\edgespace^\ell(\face^\ell)}.
\]
Also, for inf-sup stability, we need a trace $\bsigma_{\face^\ell}^{PV}$ (which we call it the PV trace, cf. Remark~\ref{rmk:const_representation}) that satisfies
\begin{equation}
(\bq^{PV}_{\agg_i^\ell})^T D^\ell_{\agg_i^\ell, \face^\ell} \,\bsigma_{\face^\ell}^{PV} \not= 0.
\label{eq:pv_condition}
\end{equation}
Because there is no guarantee that $\big(M^{\ell}_{N(\agg^\ell)}\big)^{-1} \big(D^{\ell}_{N(\agg^\ell)}\big)^T \bq$ will satisfy \eqref{eq:pv_condition}, $\bsigma_{\face^\ell}^{PV}$ is constructed separately. One computationally inexpensive choice is to take
\[
\bsigma_{\face^\ell}^{PV} := (D^\ell_{\agg_i^\ell, \face^\ell})^T \bq^{PV}_{\agg_i^\ell}.
\]
Another choice of PV trace takes edge weights into account and requires solving a local problem in $\agg^\ell_i \cup \agg^\ell_j$. The detailed construction of the second choice will be given in Appendix~\ref{sec:pv_trace}. We observe empirically that the second option typically produces higher quality coarse spaces when the edge weights have high contrast (which is usually the case for graphs coming from reservoir simulations), so the numerical results in this paper are generated using the second option.

After collecting all the traces from both sides of the face as well as the PV trace, we perform SVD to remove possible linear dependence due to restriction and union. Let the set of traces associated with $\face^\ell$ after SVD be $Tr(\face^\ell)$. Then for each trace $\bsigma_{\face^\ell}\in Tr(\face^\ell)$, we solve
\begin{equation}
 \begin{bmatrix}
    M^{\ell}_{\agg} &  \left(D^{\ell}_{\agg}\right)^T \\
    D^{\ell}_{\agg} &
  \end{bmatrix} 
  \begin{bmatrix}
    \bsigma_\agg \\ \bu
  \end{bmatrix} 
  =
  \begin{bmatrix}
    -M^\ell_{\agg,\face^\ell} \,\bsigma_{\face^\ell} \\ c_\agg\bq^{PV}_\agg - D^\ell_{\agg,\face^\ell} \,\bsigma_{\face^\ell} 
  \end{bmatrix}.
  \label{eq:extension}
\end{equation}
Here, $c_\agg$ is the constant such that $\big(\bq^{PV}_\agg\big)^T\big( c_\agg\bq^{PV}_\agg - D^\ell_{\agg,\face^\ell} \,\bsigma_{\face^\ell} \big) = 0$. By \eqref{eq:vertex_pv_vector}, $\| \bq^{PV}_\agg \| = 1$, so
\begin{equation}
c_\agg = \big(\bq^{PV}_\agg\big)^T D^\ell_{\agg,\face^\ell} \,\bsigma_{\face^\ell} .
\label{eq:extension_rhs_constant}
\end{equation}
Each trace extension is then formed by glueing together $\bsigma_{\face^\ell}$, $\bsigma_{\agg_i^\ell}$, and $\bsigma_{ \agg_j^\ell}$.
These extended traces are what ultimately end up as columns of $ P_\sigma^{\ell} $.

\begin{remark}
Note that adding $c_\agg\bq^{PV}_\agg$ to the right hand side of \eqref{eq:extension} guarantees the solvability of \eqref{eq:extension} because the nullspace of $(D^{\ell}_{\agg})^T$ is spanned by $\bq^{PV}_\agg$ (see Lemma~\ref{lemma:D_null}).
\end{remark}

\subsubsection{Bubbles}
\label{sec:bubble}

The definition of bubbles are exactly analogous to the two-level case --- solve \eqref{eq:extension} with the right-hand-side vector replaced by $ \begin{bmatrix} 0 \\ \bq \end{bmatrix} $ where $\bq$ are just the columns of $P_{u, \agg^\ell}^{\ell,NPV}$. Let $P_{\sigma, \agg^\ell}^{\ell,B}$ be the matrix whose columns are all bubbles in $\agg^\ell$. Thus, by construction 
\begin{equation}
D_\agg^\ell \; P_{\sigma, \agg^\ell}^{\ell,B} = P_{u, \agg^\ell}^{\ell,NPV} .
\label{eq:bubble_divergence}
\end{equation}
We need these bubbles in the coarse edge space to make sure the coarse saddle point problem is inf-sup stable.

\subsubsection{Structure of $P_\sigma^\ell$}
\label{sec:P_sigma}

To summarize Section~\ref{sec:trace} and~\ref{sec:bubble}, we give a high level view on the structure of the edge space interpolation matrix $P_\sigma^\ell$ here. Also, for the purpose of a clearer matchup when defining the projection operator later, different components are highlighted with different colors.

First, on each face $\face^\ell$, let $P_{\sigma, \face^\ell}^{\ell,T}$ be the matrix whose columns are all traces in $Tr(\face^\ell)$:
\begin{equation}
P_{\sigma, \face^\ell}^{\ell,T} = {\footnotesize
\begin{tikzpicture}[baseline=-\the\dimexpr\fontdimen22\textfont2\relax ]
\matrix (m)[matrix of math nodes,left delimiter={[},right delimiter={]}]
{
\bsigma_{\face^\ell}^{PV} & \hspace{1mm} & P_{\sigma, \face^\ell}^{\ell,NPV}  \\
};
\begin{pgfonlayer}{myback}
\fhighlight[blue!30]{m-1-1}{m-1-1}
\fhighlight[red!30]{m-1-3}{m-1-3}
\end{pgfonlayer}
\end{tikzpicture} }
\label{eq:trace_decomposition}
\end{equation}
where columns of $P_{\sigma, \face^\ell}^{\ell,NPV}$ are ``non-PV" traces on $\face^\ell$. Then the extension part of all traces are collected in $P_{\sigma}^{\ell,E}$. More precisely, for each column of $P_{\sigma, \face^\ell}^{\ell,T}$, which is a trace $\bsigma_{\face^\ell}$, the corresponding column in $P_{\sigma}^{\ell,E}$ is the extensions of $\bsigma_{\face^\ell}$ into the neighboring aggregates sharing $\face^\ell$ (cf. \eqref{eq:extension}). The final interpolation matrix $P_{\sigma}^\ell$ has the following structure:
\begin{equation}
P_\sigma^\ell = {\footnotesize
\begin{tikzpicture}[baseline=-\the\dimexpr\fontdimen22\textfont2\relax ]
\matrix (m)[matrix of math nodes,left delimiter={[},right delimiter={]}]
{
\ddots & {\color{red!30}\ddots} & & & & {\color{blue!0}\ddots} & & & {\color{blue!0}\ddots}  \\
{\color{blue!30}\ddots} & \ddots & & &  \\
& & \bsigma_{\face^\ell{\color{blue!30},i}}^{\ell,PV} & P_{\sigma, \face^\ell}^{\ell,NPV} & \\
& & & & \ddots & {\color{red!30}\ddots} & & & \\
\hspace{1mm} & & & & {\color{blue!30}\ddots} & \ddots & & & \hspace{1mm} \\
& & & & & & \ddots & & {\color{blue!0}\ddots} \\
& & \hspace{1mm}  & & & & & P_{\sigma, \agg^\ell}^{\ell, B} \\
\hspace{1mm} & & & & & {\color{blue!0}\ddots} & \hspace{1mm} & & \ddots \\
};
\draw[black!30,thick] (m-5-1.south west) -- (m-5-9.south east);
\draw[black!30,thick] (m-1-6.north east) -- (m-8-6.south east);
\begin{pgfonlayer}{myback}
\fhighlight[blue!30]{m-1-1}{m-2-1}
\fhighlight[red!30]{m-1-2}{m-2-2}
\fhighlight[blue!30]{m-3-3}{m-3-3}
\fhighlight[red!30]{m-3-4}{m-3-4}
\fhighlight[blue!30]{m-4-5}{m-5-5}
\fhighlight[red!30]{m-4-6}{m-5-6}
\fhighlight[green!30]{m-6-7}{m-6-7}
\fhighlight[green!30]{m-7-8}{m-7-8}
\fhighlight[green!30]{m-8-9}{m-8-9}
\node (a) at (m-1-9.north) [right=20pt]{};
\node (b) at (m-5-9.south) [right=20pt]{};
\node (c) at (m-6-9.north) [right=20pt]{};
\node (d) at (m-8-9.south) [right=20pt]{};
\draw [decorate, decoration={brace, amplitude=10pt}] (a) -- (b) 
node[midway, right=10pt] {$\displaystyle\bigoplus_{\face^\ell\in E^{\ell+1}}\bSigma^{\ell}(\face^\ell)$};
\draw [decorate, decoration={brace, amplitude=10pt}] (c) -- (d) 
node[midway, right=10pt] {$\displaystyle\bigoplus_{\agg^\ell\in V^{\ell+1}}\bSigma^{\ell}(\agg^\ell)$};
\node (e) at (m-8-1.south west) [below=10pt]{};
\node (f) at (m-8-6.south east) [below=10pt]{};
\node (g) at (m-8-7.south west) [below=10pt]{};
\node (h) at (m-8-9.south east) [below=10pt]{};
\draw [decorate, decoration={brace, amplitude=10pt, mirror}] (e) -- (f) 
node[midway, below=10pt] {$\displaystyle\bigoplus_{\face^\ell\in E^{\ell+1}}\bSigma^{\ell+1}(\face^\ell)$};
\draw [decorate, decoration={brace, amplitude=10pt, mirror}] (g) -- (h) 
node[midway, below=10pt] {$\displaystyle\bigoplus_{\agg^\ell\in V^{\ell+1}}\bSigma^{\ell+1}(\agg^\ell)$};
\node (i) at (m-7-3) [right=0pt]{\normalsize$P_{\sigma}^{\ell,E}$};
\end{pgfonlayer}
\end{tikzpicture} }
\label{eq:Pedge_structure}
\end{equation}

\subsection{Coarse projections and their properties}

We have defined the interpolation matrices $P_u^\ell$ and $P_\sigma^\ell$, and the coarse vertex and edge spaces in previous sections. To use the hierarchy in multilevel Monte Carlo simulations and nonlinear multigrid solvers, we will also need projection operators that take quantities from a given level to a coarser level. We have already discussed the projection operator $Q_u^\ell = (P_u^\ell)^T$ for the vertex space in Section~\ref{sec:vdof}. Due to our aggregation-based coarsening, it turns out that projection operators for edge spaces can also be defined locally in each aggregate and face.\\

\begin{proposition}
There exists a locally constructed projection $\pi_{\sigma}^\ell := P_{\sigma}^\ell \, Q_{\sigma}^\ell $ to the coarse edge space , where $Q_{\sigma}^\ell: \bSigma^{\ell} \to \bSigma^{\ell+1}$ satisfies 
\[
Q_{\sigma}^\ell P_{\sigma}^\ell = I_\sigma^{\ell+1}.
\]
\label{prop:projection}
\end{proposition}

\begin{proposition}
\label{prop:commutativity}
On each level $\ell$, we have
\begin{equation}
Q_u^\ell D^\ell = D^{\ell+1} Q_\sigma^\ell.
\label{eq:commutativity}
\end{equation}

In other words, the following commutative diagram holds:
\[ 
\begin{tikzpicture}
  \matrix (m) [matrix of math nodes,row sep=3em,column sep=4em,minimum width=2em]
  {
     \edgespace^\ell &  \vertexspace^\ell \\
     \edgespace^{\ell+1} & \vertexspace^{\ell+1} \\};
  \path[-stealth]
    (m-1-1) edge node [above] {$ D^\ell $} (m-1-2) 
    (m-1-1) edge node [left] {$Q_\sigma^\ell \;$} (m-2-1)
    (m-2-1) edge node [above] {$D^{\ell+1}$} (m-2-2) 
    (m-1-2) edge node [right] {$\; Q_u^\ell$} (m-2-2);
\end{tikzpicture}.\vspace{-5mm}
\]
\end{proposition}

{\it Proof.} 
The proposition is a direct consequence of Lemma~\ref{lemma:commutativity} and~\ref{lemma:D_null}. \hfill$\square$\\
\\
We defer the detailed construction of $Q_{\sigma}^\ell$, as well as the proof of Proposition~\ref{prop:projection} to Appendix~\ref{sec:Q_sigma}. With Proposition~\ref{prop:commutativity}, the well-posedness of the coarse saddle point problems can be shown as in Section 7 of \cite{blv17}.

\section{Applications of the multilevel hierarchy}
\label{sec:numerics}

In this section, we demonstrate some applications of our multilevel spectral coarsening algorithm, including multilevel upscaling, multilevel Monte Carlo simulations, and nonlinear multigrid solvers. The numerical results here are obtained using  the open source software {\it smoothG} \cite{smoothg}, which contains a C++ implementation with MPI parallelization of the multilevel spectral coarsening algorithm. It is available on the github repository: \url{github.com/LLNL/smoothG}. The graph visualizations in Section~\ref{sec:upscaling} are generated by NetworkX \cite{networkx}, while the visualizations for PDE related applications in Section~\ref{sec:MLMC} and~\ref{sec:FAS} are generated using MFEM \cite{mfem} and GLVis \cite{glvis}.

\subsection{Multilevel upscaling versus two-level}\label{sec:upscaling}

\begin{figure}[h!]
\centering
\includegraphics[scale=.5]{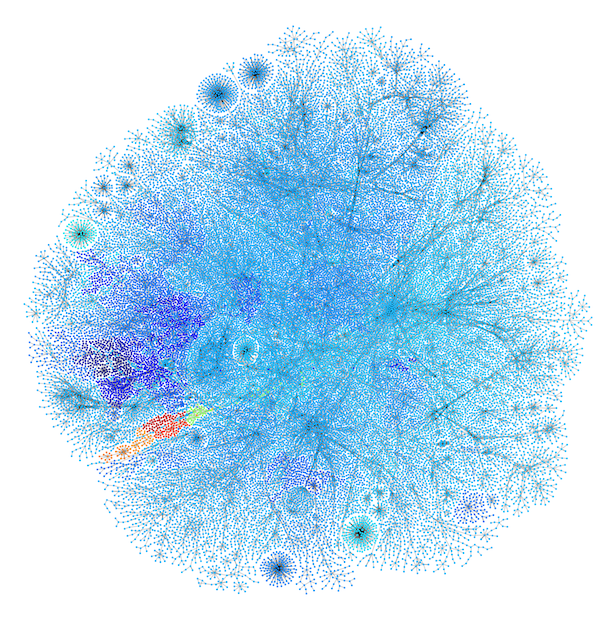}
\caption{Fiedler vector of an internet graph (dim($\bU^0$) = 35638).}
\label{fig:fiedler}
\end{figure}

\begin{figure}[h!]
\centering
\includegraphics[scale=.07,clip,trim=4cm 5.5cm 3cm 3.5cm]{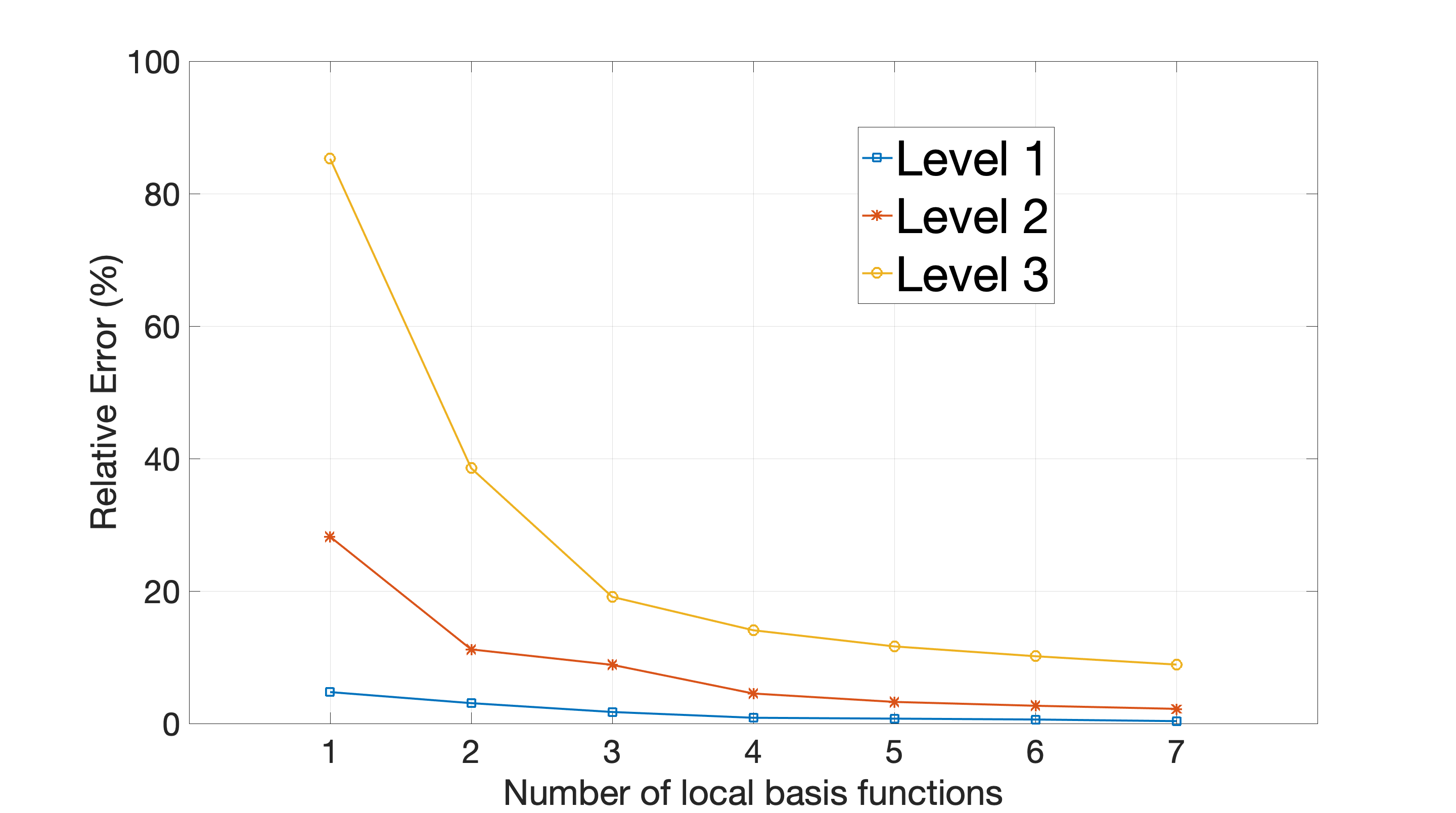}
\includegraphics[scale=.07,clip,trim=4cm 5.5cm 3cm 3.5cm]{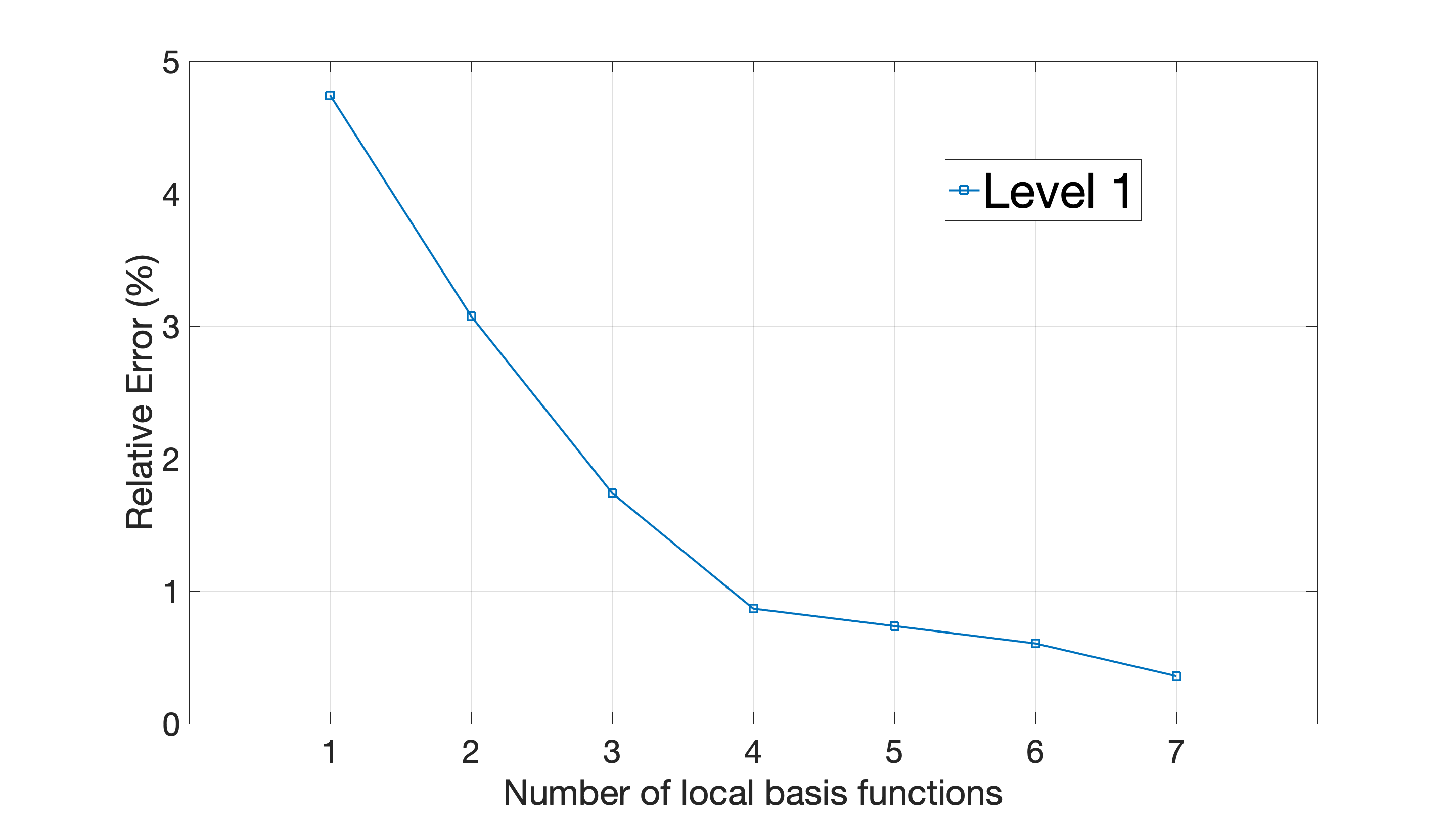}\\
\hspace{1.5mm}\includegraphics[scale=.07,clip,trim=4cm 0 3cm 1cm]{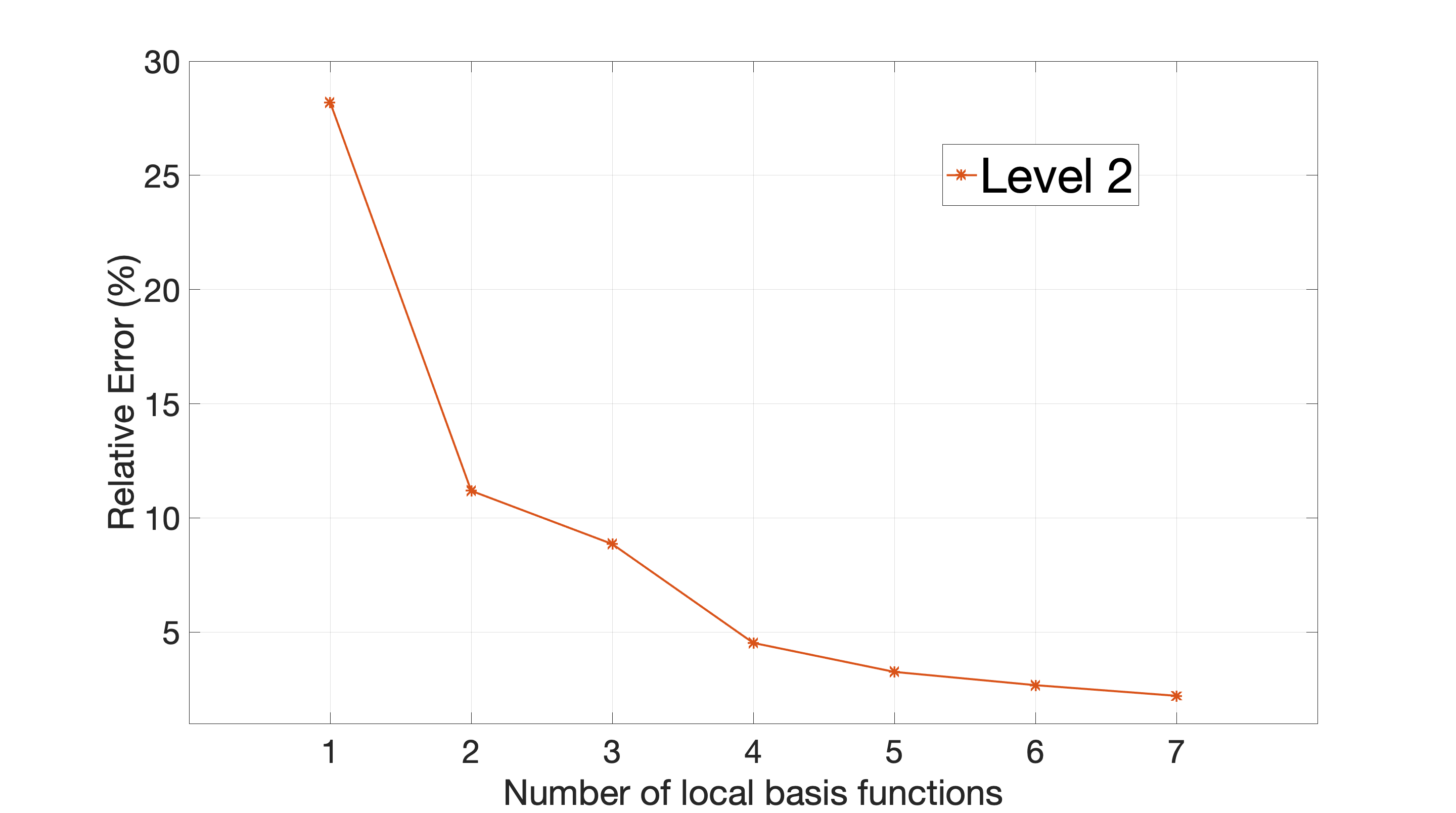}
\includegraphics[scale=.07,clip,trim=4cm 0 3cm 1cm]{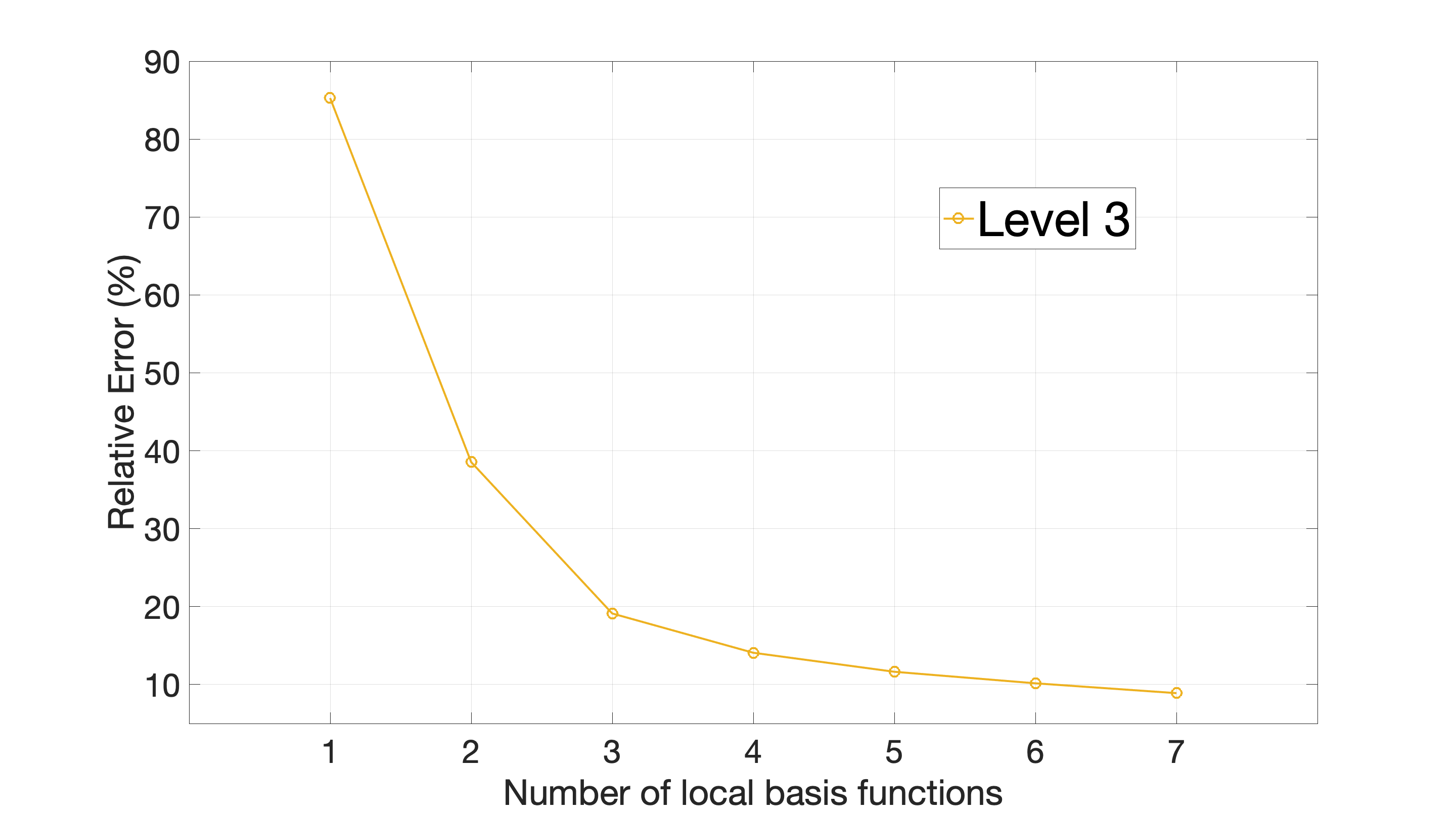}
\caption{Rapid decay of upscaling error in all levels as number of local basis functions increases.}
\label{fig:ml_upscale_error}
\end{figure}

We consider an internet graph from \url{http://opte.org}. The graph has 35638 vertices and 42827 edges, and its Fielder vector is depicted in Figure~\ref{fig:fiedler}. In this example, the right hand side $\bff$ is set to be the Fiedler vector (scaled by the smallest positive eigenvalue of the associated graph Laplacian). We apply our multilevel upscaling algorithm to obtain 4-level hierarchies of approximation spaces by varying the number of local eigenvectors ($m_\agg$) selected as basis functions. The coarsening factor is 16 for all hierarchies. The upscaling errors are summarized in Table~\ref{tab:ml_upscale_error} and Figure~\ref{fig:ml_upscale_error}. In particular, Figure~\ref{fig:ml_upscale_error} clearly shows that the upscaling errors decay rapidly in all levels when more local eigenvectors are added to the coarse spaces. 

\begin{table}[h!]
\centering
\caption{Relative upscaling errors in all levels (\%)}
\label{tab:ml_upscale_error}\scalebox{1.}{
\begin{tabular}{cccccccc}
\toprule
      &  \multicolumn{7}{c}{$m_\agg$} \\
   \cmidrule(l){2-8}
& 1 & 2 & 3 & 4 & 5 & 6 & 7 \\
\midrule
Level 1 & 4.744 & 3.075 & 1.739 & 0.8672 & 0.7362 & 0.6041 & 0.3576 \\
Level 2 & 28.20 & 11.19 & 8.859 & 4.529 & 3.260 & 2.679 & 2.213 \\
Level 3 & 85.27 & 38.55 & 19.12 & 14.08 & 11.64 & 10.16 & 8.897 \\
\bottomrule
\end{tabular}}
\end{table}

\begin{figure}[h!]
\centering
\subfigure[$m_\agg=1$, level 1: \#dofs = 2218]{\includegraphics[scale=.5]{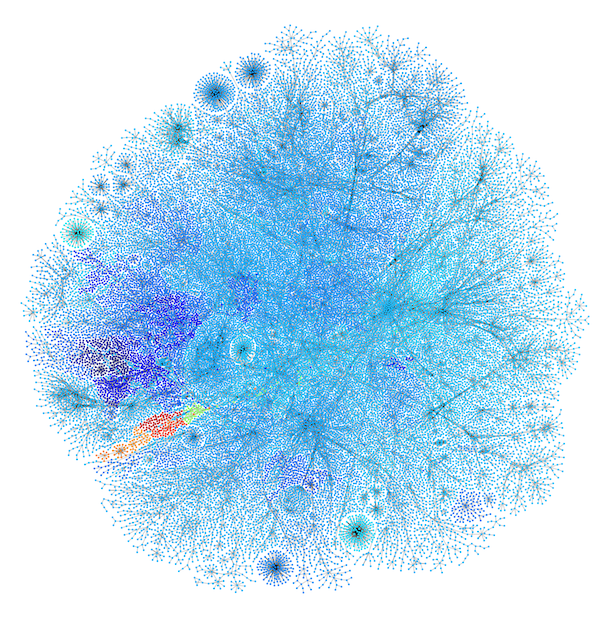}}
\subfigure[$m_\agg=1$, level 2: \#dofs = 139]{\includegraphics[scale=.5]{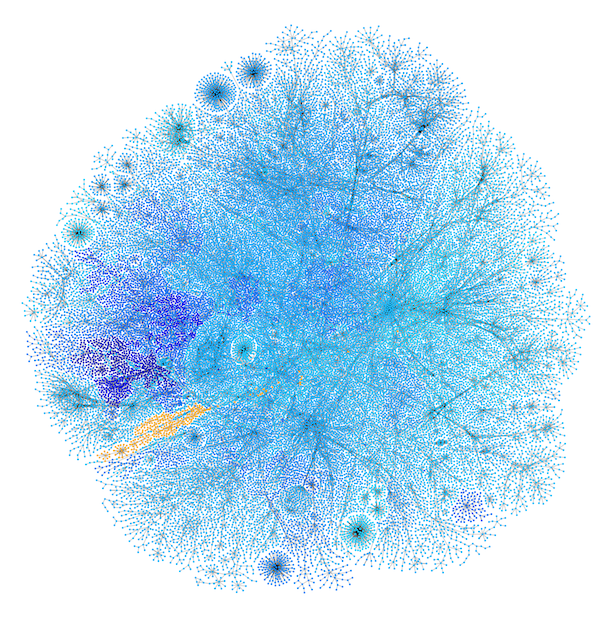}}
\subfigure[$m_\agg=1$, level 3: \#dofs = 9]{\includegraphics[scale=.5]{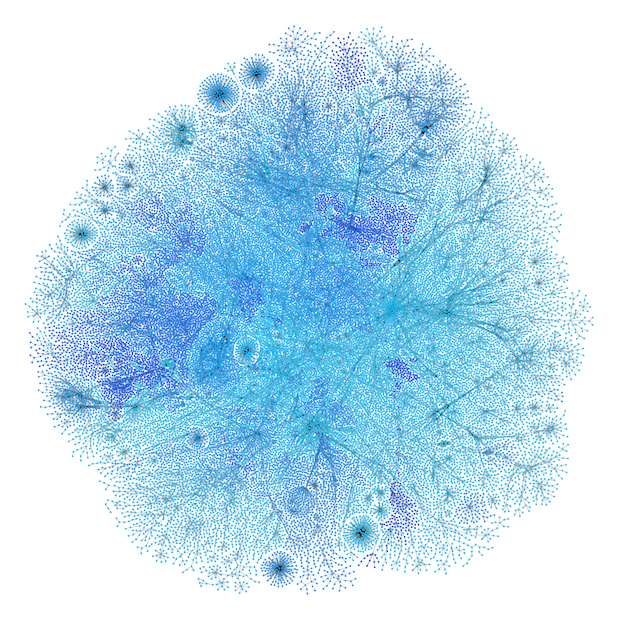}}
\subfigure[$m_\agg=4$, level 1: \#dofs = 8625]{\includegraphics[scale=.5]{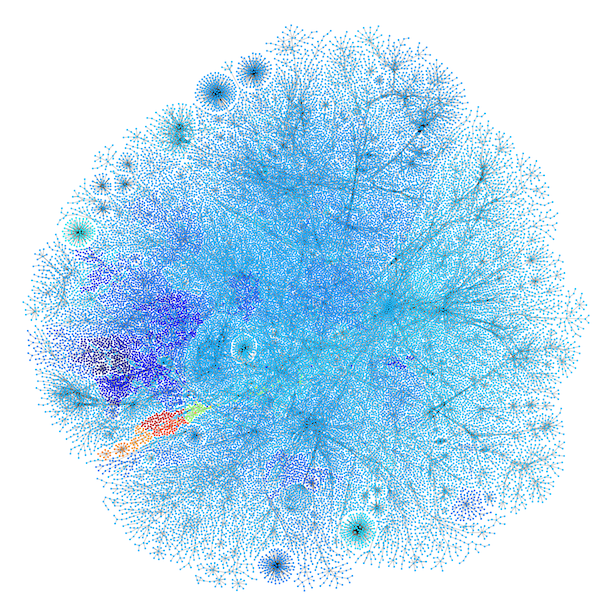}}
\subfigure[$m_\agg=4$, level 2: \#dofs = 552]{\includegraphics[scale=.5]{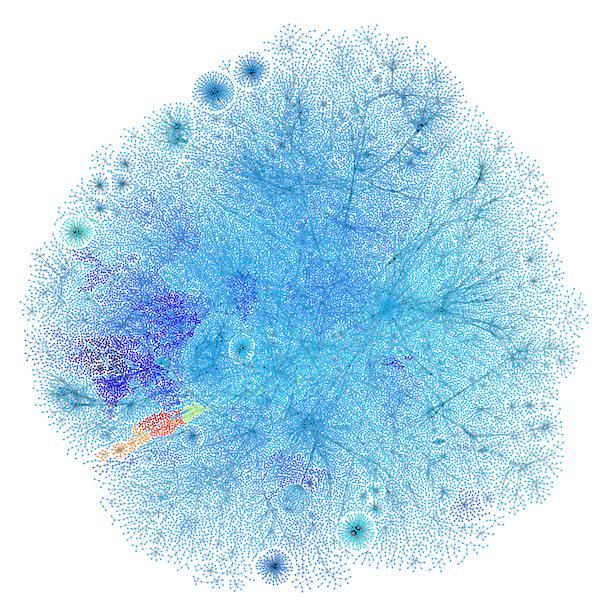}}
\subfigure[$m_\agg=4$, level 3: \#dofs = 36]{\includegraphics[scale=.5]{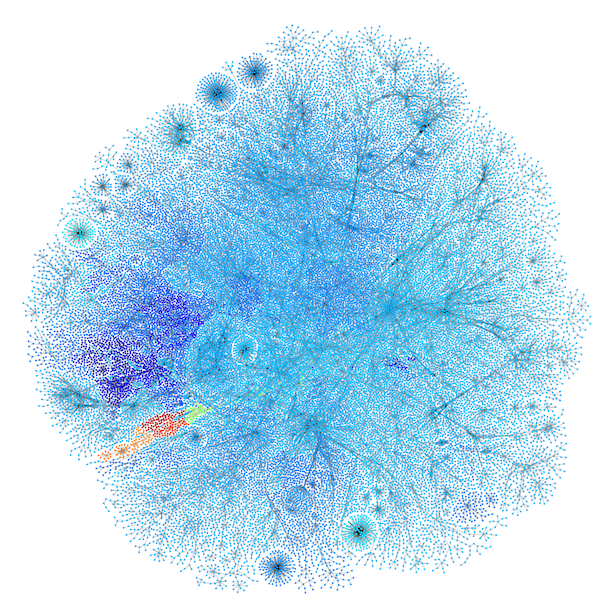}}
\caption{Multilevel coarse approximations of Fiedler vector; top: $m_\agg=1$, bottom: $m_\agg=4$.}
\label{fig:ml_fiedler}
\end{figure}

In Figure~\ref{fig:ml_fiedler}, the approximated Fiedler vector is shown on each coarse level. When $m_\agg=1$, we see that the coarse approximations start to deviate from the fine level solution (Figure~\ref{fig:fiedler}); while the coarse solutions corresponding to $m_\agg=4$ approximate the fine level solution very well across all levels.

\subsection{Multilevel Monte Carlo for finite volume schemes}\label{sec:MLMC}

In \cite{ovv17}, a multilevel Monte Carlo method for forward uncertainty quantification was shown to have efficiency gains over standard Monte Carlo in a mixed finite element setting.
Here we present similar results for a finite volume method, which can be understood as a graph problem by considering each finite volume as a graph vertex and each face between volumes as a graph edge, see \cite{blv17} for details.

We are interested in the (linear) Darcy equations for a pressure $p$
\begin{equation}
  -\nabla \cdot \left( \kappa(\omega) \nabla p \right) = f \text{ in } \Omega \label{eq:linear-darcy}
\end{equation}
with appropriate boundary conditions and forcing function $f$ to be described in the test cases below.
Here $ \omega $ is an uncertain parameter describing uncertainty about the coefficient $ \kappa $.
The solution $p$ also depends implicitly on $\omega$, and that uncertainty is what we seek to quantify here; in particular, we want to propagate the uncertainty in $ \kappa $ via the forward problem \eqref{eq:linear-darcy}
and compute statistics (e.g. expected value) of a functional of the solution $p$. 

We assume $ \kappa(\omega) $ is a lognormal random field, in particular $ \kappa = \exp(\theta) $ where $\theta$ is a stationary isotropic Gaussian field with exponential covariance given by
\[
  \text{cov}(x, y) = \sigma^2 \exp(-\corr \| y - x \|)
\]
where $\sigma^2$ is the marginal variance and $\corr$ is the correlation length.

To be very brief, the basic idea of multilevel Monte Carlo is to estimate the expected value
$\mathbb{E}[Q_0]$ of $Q_0$, a quantity of interest defined on the finest level, by computing 
standard Monte Carlo estimators (using a sample and average technique) of the telescoping sum 
\begin{equation}
  \mathbb{E}[Q_0] =  \mathbb{E}[Q_L] + \sum_{\ell = 0}^{L-1} \mathbb{E}[Q_\ell - Q_{\ell + 1}] .
\end{equation}

The correction terms $ Y_\ell = Q_\ell - Q_{\ell + 1} $ have smaller variance than the original quantities of interest, and so a Monte Carlo estimator for $ Y_\ell $ converges relatively quickly.
At the same time, though the coarsest estimator for $ \mathbb{E}[Q_L] $ may need a large number of samples, those samples are cheap because they are performed at a coarse level.
See, for example, \cite{cliffe-giles-mlmc} for more details on the efficiency and accuracy of multilevel Monte Carlo.

Because the emphasis here is on the finite volume discretization, the setting for the experiments is somewhat simpler than in \cite{ovv17}.
In particular, although we use the stochastic PDE sampler from \cite{lrl11, ovv17} to generate realizations of the lognormal random field, we do not embed the domain of interest in a larger domain, and so our sampler is subject to some modest boundary effects.
There is no real theoretical or computational barrier to using such mesh embedding, however.

\subsubsection{Unit square results}

\input{dirichlet_bc_square.tex}

For demonstration, we employ multilevel Monte Carlo to estimate the expectation of a quantity of interest on the unit square $ \Omega = [0, 1]^2 $ with boundary conditions $ \nabla p \cdot n = 0 $ on $ x = 0, 1 $, $ p = -1 $ on $ y = 1 $, and $ p = 0 $ on $ y = 0 $.
The quantity of interest is the average flux on the top boundary, that is
\[
  Q = \int_{x=0}^{x=1} \kappa \nabla p(x, 1) \cdot n \, dx .
\]

Some results are shown in Figures \ref{fig:dirichlet_bc_square_coeff}, \ref{fig:dirichlet_bc_square}, where we see that the correction terms have smaller variance than the fine-level estimator, as expected, and that the MLMC algorithm scales well with respect to problem size.
Computational cost in comparison to a standard single-level Monte Carlo algorithm is shown in Table \ref{tab:dirichlet_square_speedup}, where significant speedup is demonstrated.

\input{dirichlet_square_speedup.tex}




\subsubsection{Unit cube results}

\input{dirichlet_bc_cube.tex}

For another demonstration, we apply the multilevel Monte Carlo procedure on the unit cube $ \Omega = [0, 1]^3 $ with boundary conditions $ \nabla p \cdot n = 0 $ on $ x, y = 0, 1 $, $ p = -1 $ on $ z = 1 $, and $ p = 0 $ on $ z = 0 $. Like the unit square problem, the quantity of interest is the average flux across the top boundary ($z = 1$). 

Our results are shown in Figures \ref{fig:dirichlet_bc_cube_coeff}, \ref{fig:dirichlet_bc_cube}.
Computational cost in comparison to a one-level Monte Carlo algorithm is reported in Table \ref{tab:dirichlet_cube_speedup}, where significant speedup is demonstrated in reaching a target mean square error using 
an adaptive algorithm where the number of samples computed on each level is chosen 
to minimize the variance for a fixed computational cost (see, e.g., \cite{cliffe-giles-mlmc} for details). 
In particular, the variance results in \ref{fig:dirichlet_bc_cube} show that the most expensive levels have the smallest variance, which is the key to the speedup in multilevel Monte Carlo as fewer samples are necessary on the most expensive levels. 

\input{dirichlet_cube_speedup.tex}

\subsubsection{Results inspired by the SPE10 dataset}

\input{spe10.tex}

The SPE10 dataset \cite{spe10} is very often used as a challenging example of a complicated discontinuous coefficient.
As given, this dataset is deterministic, and so we adapt one of its prominent features for our sampling and uncertainty quantification tests, namely its layered structure.
We produce samples with a similar structure on the SPE10 domain by imposing an anisotropic covariance structure, where the correlation length in horizontal directions is much longer than the correlation length in the vertical direction \cite{simpson-lindgren-continuous}.

The results for this example are performed on a parallel computer with eight processors.
We set $ \kappa = 0.1 $ in the vertical direction and $ \kappa = 0.001 $ in the horizontal direction.
The boundary conditions are similar to those for the unit cube case, with $ p = 0 $ imposed on one side of the domain and $ p = -1 $ on the opposite side, while the flux across this latter boundary is the quantity of interest.
We demonstrate the multilevel Monte Carlo method with some results shown in Figures \ref{fig:spe10_coeff}, \ref{fig:spe10} which show that the multilevel Monte Carlo procedure is effective in reducing variance also for this problem with an anisotropic covariance.

\begin{table}
\caption{Sample counts and computational cost for the SPE10 problem, including one-level and multilevel Monte Carlo for various tolerances on the mean-square error (MSE).
For multilevel Monte Carlo, the sample counts are listed from finest to coarsest level.}
\label{tab:spe10_speedup}
\begin{center}\begin{tabular}{llllll}
\toprule
      &  \multicolumn{2}{c}{one-level}&  \multicolumn{3}{c}{multilevel} \\
   \cmidrule(l){2-3}  \cmidrule(l){4-6}
MSE&  samples&  time(s)&  samples&  time(s)&  speedup \\
\midrule
1e+00&  284&  1021.1&  [20, 20, 49, 274]&  73.2&  13.9 \\
5e-01&  537&  1911.1&  [20, 31, 128, 663]&  72.7&  26.3 \\
2e-01&  1323&  4745.0&  [20, 107, 454, 2657]&  74.9&  63.4 \\
1e-01&  2548&  9175.3&  [24, 650, 2448, 13303]&  107.3&  85.5 \\
1e-02&  ---&  ---&  [136, 5827, 20561, 108840]&  671.1&  --- \\
\bottomrule
\end{tabular}\end{center}
\end{table}

In Table \ref{tab:spe10_speedup}, we show the computational efficiency gains of the multilevel Monte Carlo algorithm for this problem.
In this problem setting, the computational time is completely dominated by the number of samples taken on the finest level.
As a result, MLMC results in very large performance gains, offering speedups of greater than 50 compared to single-level Monte Carlo.
The MLMC algorithm is able to reach an empirical mean-square error level of 0.01 in less time than the single-level algorithm takes to get to an (absolute) error of 1.0 in the quantity of interest.

\subsection{Nonlinear multigrid (full approximation scheme)}\label{sec:FAS}

Consider the following nonlinear elliptic problem:
\begin{equation}
\begin{split}
-\nabla\cdot\big(\kappa(p)\nabla p\big) = f  & \text{ in } \Omega, \\
p = 0 & \text{ on } \bdr \Omega.
\end{split}
\label{eq:nonlinear}
\end{equation}
The permeability function $\kappa(p)$ takes the form:
\[
\kappa(p) = \kappa_0 e^{\alpha p}
\]
where $\kappa_0$ is some reference permeability field when $p = 0$. As an example, we take the domain $\Omega$ to be the Egg model from \cite{EggModel}, with $\alpha = 5$, $f = 1$, and $\kappa_0$ is realization 27 of the Egg model. The problem is discretized by the finite volume method with two-point flux approximation (TPFA). The discrete problem has the form:
\begin{equation}
D \big(M(\bp)\big)^{-1} D^T \bp = \bff.
\label{eq:discrete_nonlinear}
\end{equation}
Since our coarsening scheme is in the mixed setting, we reformulate \eqref{eq:discrete_nonlinear} as:
\begin{equation}
A(\bx) := \begin{bmatrix} M(\bp) & D^T \\ D &  \end{bmatrix}
  \begin{bmatrix} \bv \\ \bp \end{bmatrix} =  \begin{bmatrix} \bzero \\ -\bff \end{bmatrix} =: \bb
\label{eq:discrete_nonlinear_mixed}
\end{equation}
where $\bx := (\bv, \bp)$.

To solve \eqref{eq:discrete_nonlinear_mixed}, we use full approximation scheme (FAS) \cite{h03}, which applies the idea of multigrid directly to the nonlinear problem. In order for FAS to have rapid convergence, it is important to have coarse spaces with good approximation properties. To this end, we utilize our multilevel coarsening algorithm presented in Section~\ref{sec:algorithm} to construct the multigrid hierarchy. In the setup phase, the hierarchy is constructed by considering the reference permeability $\kappa(0) = \kappa_0$. In the solve phase, one cycle of FAS is summarized in Algorithm~\ref{alg:FAS}. We compare numerically the performance of FAS to Picard iteration for solving \eqref{eq:discrete_nonlinear}. For both methods, nonlinear iteration stops when (relative/absolute) $\ell_2$ norm of the residual is less than 1e-8/1e-10. \\

\begin{algorithm}
\begin{algorithmic}[1]
\Function{\tt NonlinearMG}{$\bx^\ell, \bb^\ell$}
    \If{$\ell$ == coarsest level}
        \State Solve $A^\ell(\bx^\ell) = \bb^\ell$ for $\bx^\ell$
    \Else
        \State $\bx^\ell \leftarrow$ \Call{\tt NonlinearRelaxation}{$\bx^\ell, \bb^\ell$}
        \State $\bx^{\ell+1} \leftarrow Q^\ell \bx^\ell$
        \State $\bb^{\ell+1} \leftarrow (P^\ell)^T \big(\bb^\ell - A^\ell(\bx^\ell)\big) + A^{\ell+1}(\bx^{\ell+1})$
        \State $\by^{\ell+1} \leftarrow$ \Call{\tt NonlinearMG}{$\bx^{\ell+1}, \bb^{\ell+1}$}
        \State $\bx^\ell \leftarrow \bx^\ell + P^\ell(\by^{\ell+1} - \bx^{\ell+1}) $
        \State $\bx^\ell \leftarrow$ \Call{\tt NonlinearRelaxation}{$\bx^\ell, \bb^\ell$}
    \EndIf
        \State \Return{$\bx^\ell$}
\EndFunction
\end{algorithmic}
\caption{Full approximation scheme (V-cycle)}
\label{alg:FAS}
\end{algorithm}


\noindent{\bf Multigrid and coarsening setting:}\vspace{2mm}\\
\begin{tabular}{ll}
Coarsening factor: 28  &   Number of vertex basis per aggregate: $m_\agg$ \\
Nonlinear relaxation: 20 linearized iterations \hspace{5mm}  &  Number of edge basis per face: 1   \\
\end{tabular}
\vspace{5mm}

\begin{table}[h!]
\centering
\caption{Comparison of solve time (\#iter) of FAS and Picard method for solving \eqref{eq:nonlinear}}
\label{tab:nonlinear} \vspace{0mm}\scalebox{1.}{
\begin{tabular}{cccccccc}
\toprule
& \multicolumn{2}{c}{FAS info} &  \multicolumn{4}{c}{FAS (multilevel Picard)} & one-level   \\
 \cmidrule(l){2-3}  \cmidrule(l){4-7}
 \#dofs & \#levels & OC & $m_\agg$=1 & $m_\agg$=5 & $m_\agg$=10 & $m_\agg$=15 & Picard \\
 \midrule
 84,580 & 3 & 1.22 & 1.24s (7) & 0.95s (5) & 0.71s (4) & 0.52s (3)  & 0.96s (11) \\
 636,248 & 4 & 1.26 & 14.79s (9) & 5.07s (4) & 4.24s (3) & 4.52s (3) & 9.58s (13) \\
4,919,776 & 4 & 1.3 & 114.8s (9) & 53.4s (4) & 41.0s (3) & 57.5s (4)& 97.3s (13)  \\
\bottomrule
\end{tabular}}
\end{table}

\def\scaleone{.5}
\begin{figure}[h!]
\centering
\subfigure{
\includegraphics[scale=\scaleone,clip,trim=5.5cm 8cm 9.7cm 10cm]{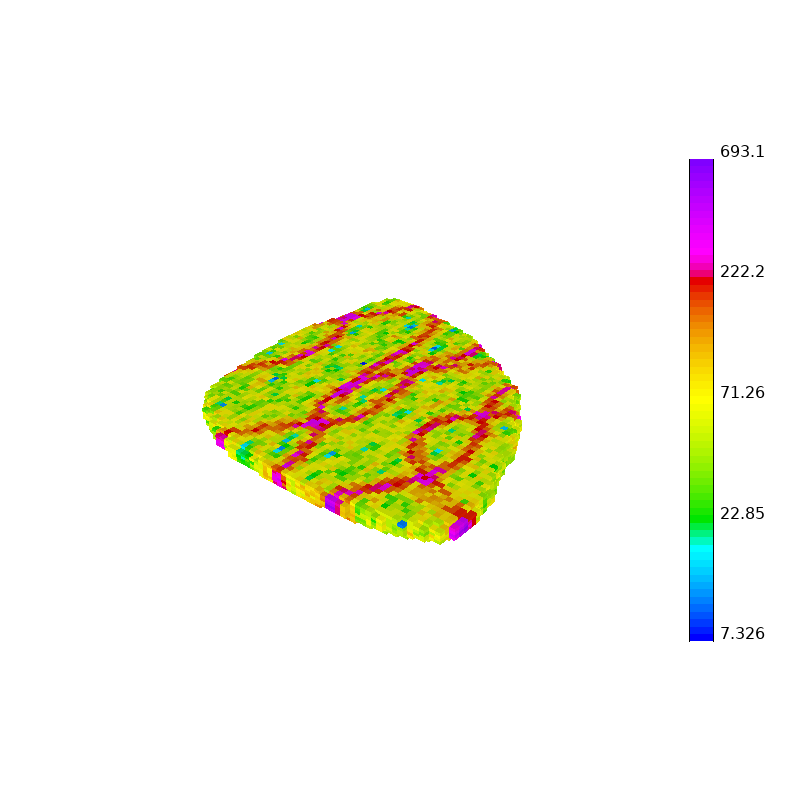}
\includegraphics[scale=.25,clip,trim=23cm 4.6cm 1cm 4.6cm]{egg_perm27.png}}
\subfigure{
\includegraphics[scale=\scaleone,clip,trim=5.5cm 8cm 9.7cm 10cm]{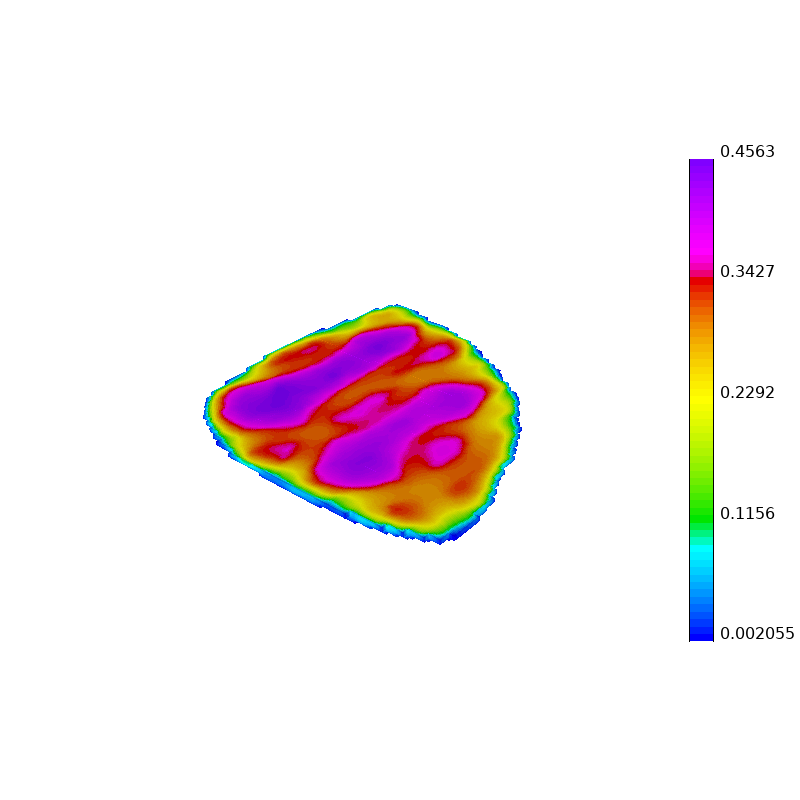}
\includegraphics[scale=.25,clip,trim=23cm 4.6cm 0 4.6cm]{egg_sol.png}}
\caption{Reference permeability $k_0$ (left) and pressure solution (right)}
\label{fig:egg_perm}
\end{figure}

\def\scaletwo{.2}
\begin{figure}[h!]
\centering
\subfigure{
\includegraphics[scale=\scaletwo,clip,trim=10cm 0cm 15cm 8cm]{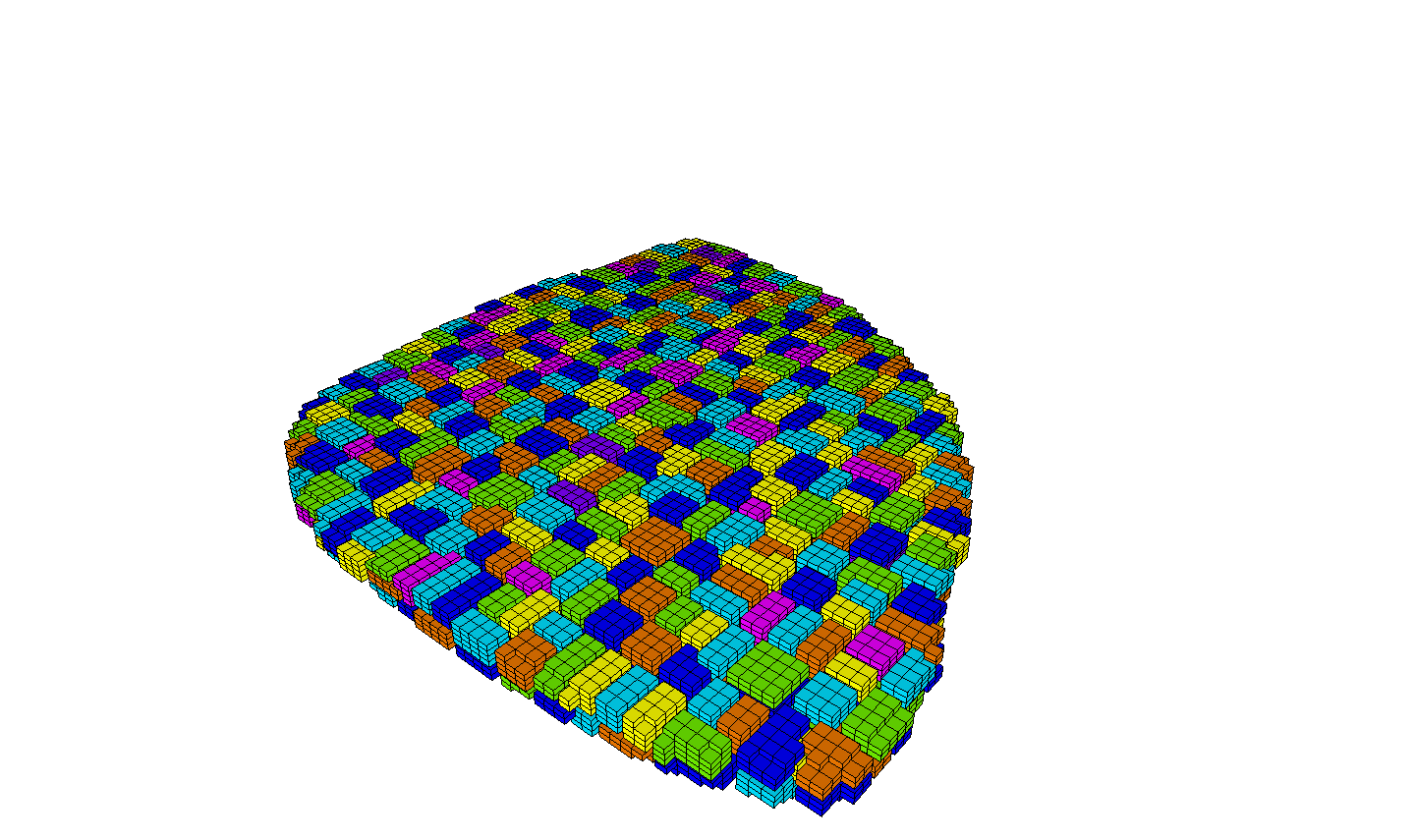}}
\subfigure{
\includegraphics[scale=\scaletwo,clip,trim=10cm 0cm 15cm 8cm]{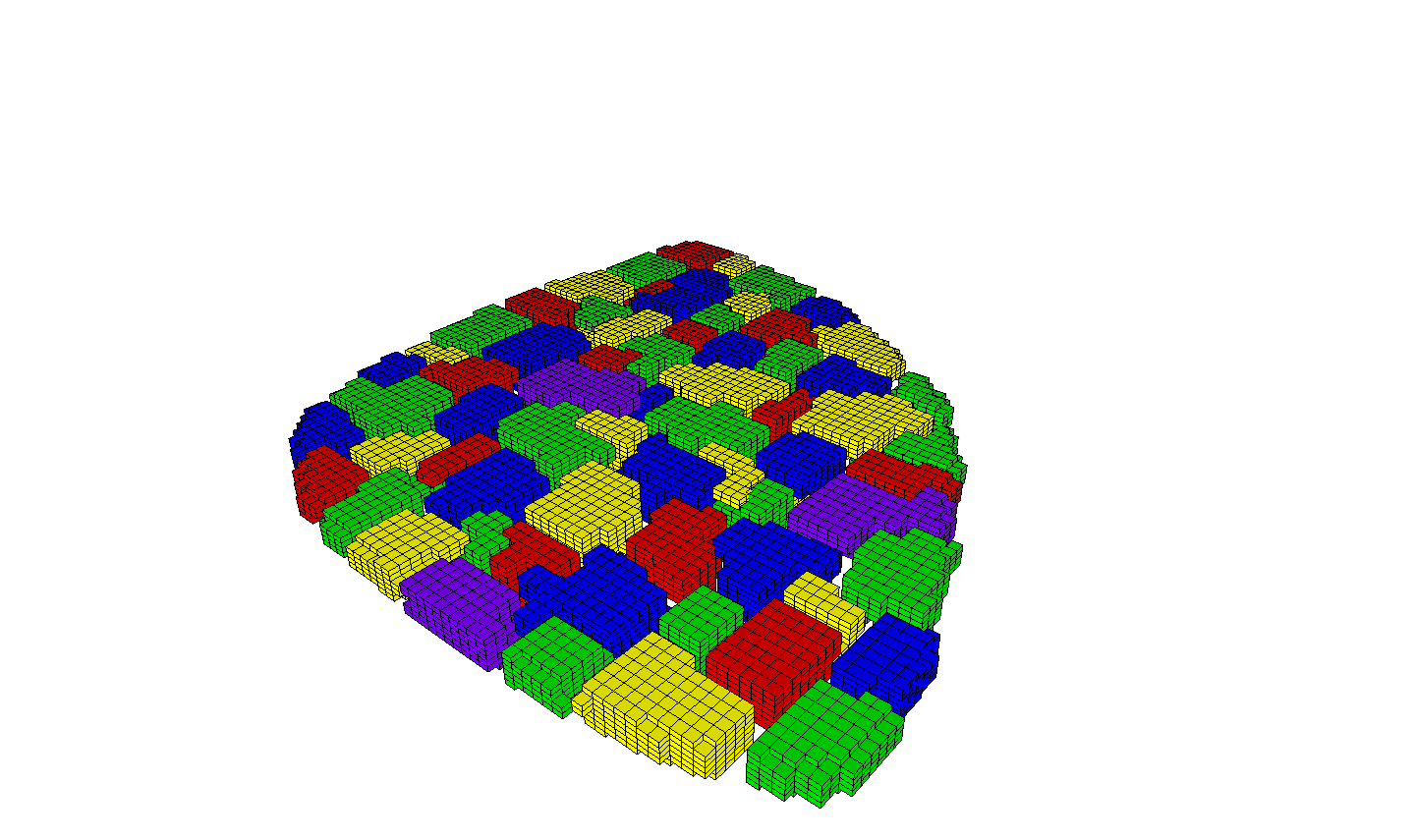}}
\subfigure{
\includegraphics[scale=\scaletwo,clip,trim=10cm 0cm 15cm 8cm]{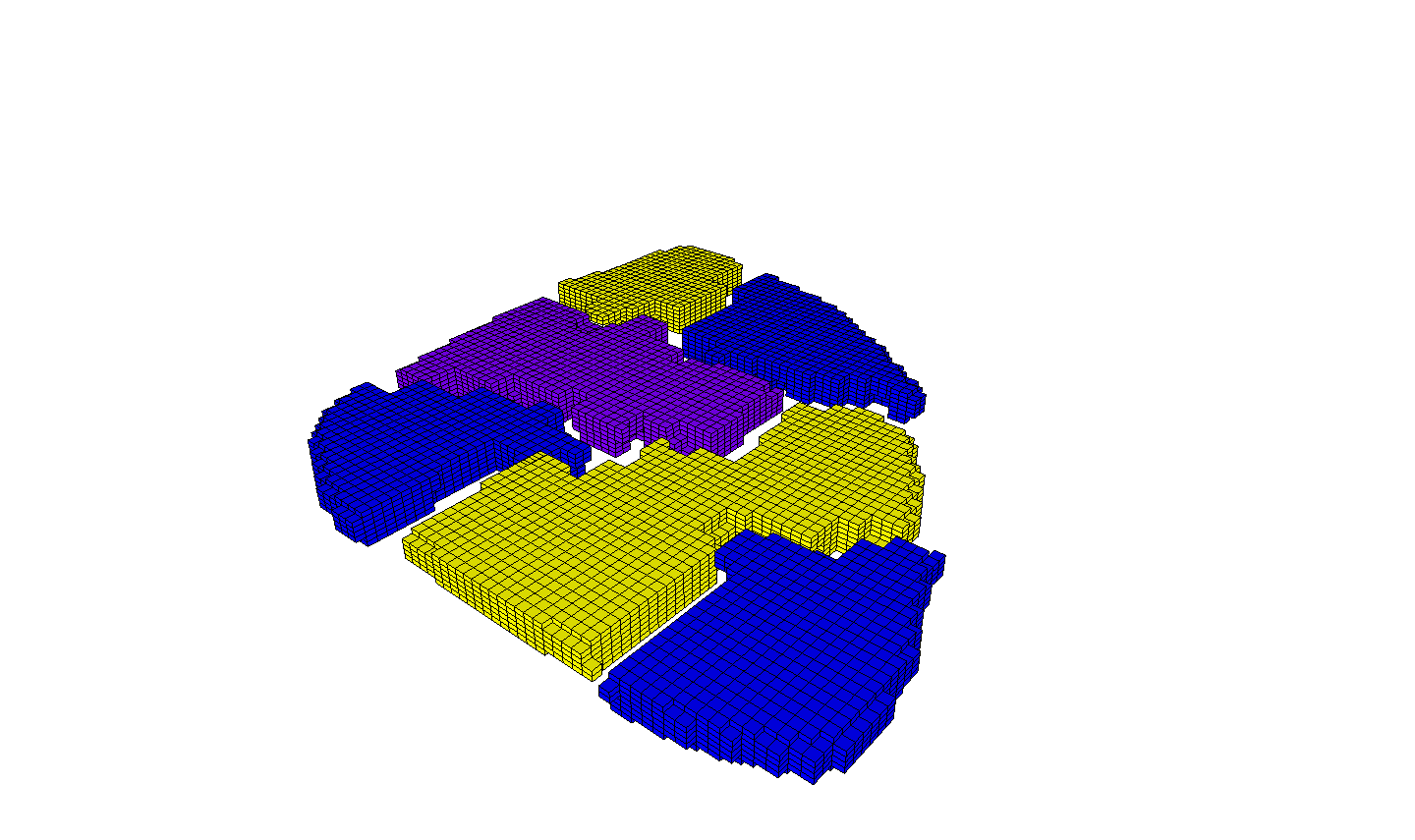}}
\caption{Hierarchical partitioning}
\label{fig:egg_partition}
\end{figure}

\begin{figure}[h!]
\centering
\subfigure[Nonlinear iterations w.r.t. \# local basis]{\includegraphics[scale=.088,clip,trim=5cm 0cm 9cm 0cm]{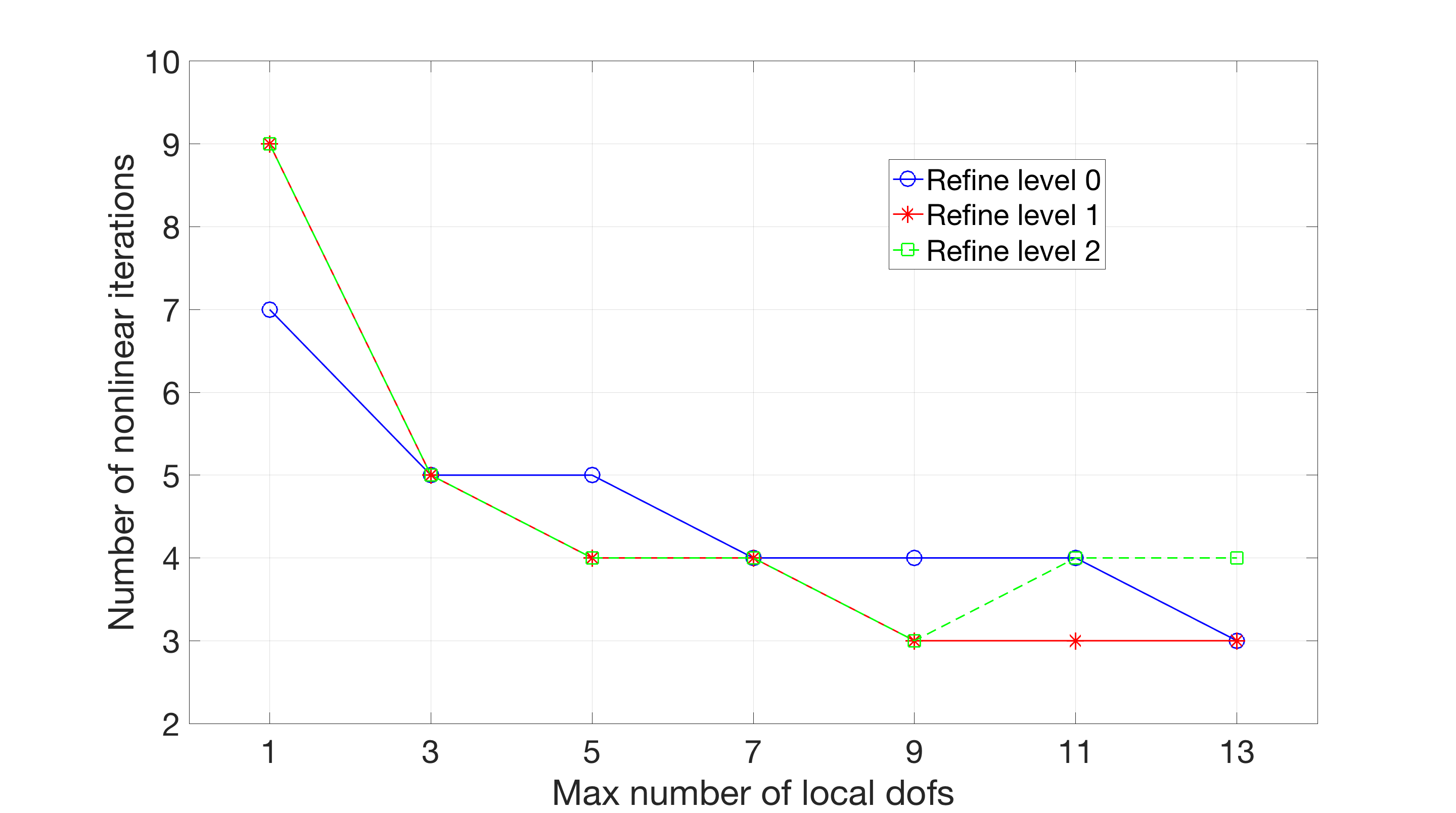}}
\hspace{5mm}
\subfigure[Solve time comparison with Picard method]{\includegraphics[scale=.088,clip,trim=5cm 0cm 9cm 0cm]{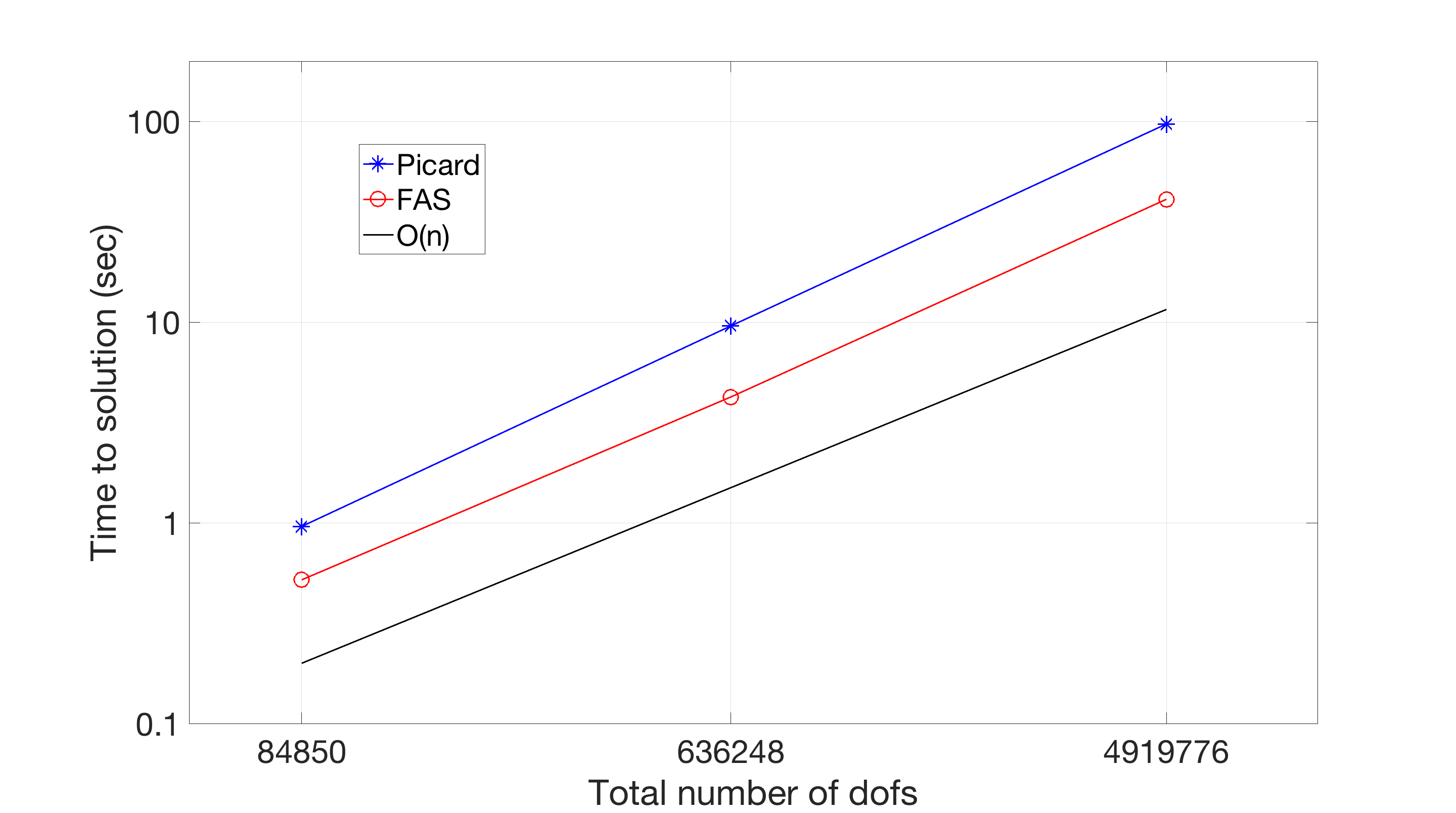}}
\caption{Performance of FAS}
\label{fig:egg_conv}
\end{figure}

From Table~\ref{tab:nonlinear} and Figure~\ref{fig:egg_conv}, it is evident that using more local spectral basis functions in the hierarchy does reduce the number of nonlinear iteration (and solving time) of FAS, which clearly shows the benefit of using higher order coarse approximation. Both Picard iteration and FAS exhibit optimal scaling ($\mathcal{O}(n)$), while FAS is about 2 times faster than Picard iteration if more local spectral basis functions are used ($m_\agg \ge 9$).
Notice that we solve the coarse problems by hybridization \cite{ArnoldBrezzi85, lv16, dkltv19}, which reduces the original saddle point problem to a symmetric positive definite problem on edges, so the sparsity patterns of the global systems in coarse levels are determined by the edge-based dofs. In our coarsening setting, we fix the number of edge-based dofs to be 1 per face, so when we increase $m_\agg$, the numbers of non-zeros of the global systems and therefore operator complexity (OC) are not affected. Consequently, in the solving stage, the extra cost of using a hierarchy with higher order spectral coarsening is only in the transformation of right hand sides and solutions between the original saddle point problem and the hybridized system, which is local and highly scalable. A more comprehensive study on building nonlinear multigrid solvers using local spectral coarsening is carried out in \cite{fas-spectral}.

\section{Conclusions}
\label{sec:conclusion}

This paper has extended the two-level graph coarsening algorithm of \cite{blv17} to a recursive multilevel framework, using the key step of coarsening degrees of freedom rather than graph edges and vertices.
We demonstrated that the coarsening procedure yields coarse representatives of a problem that are much cheaper to solve than the fine problem, but still retains some useful approximation properties that can be used to speed up nonlinear solves and forward uncertainty quantification problems in a finite volume context.
Possible future work could include developing faster solvers for the resulting saddle point problems or ways to speed up the (fairly expensive) setup phase by approximating local eigenvalue problems.

\bibliographystyle{abbrv}
\bibliography{multilevel}

\begin{appendices}

\section{A locally constructed projection for the edge space}\label{sec:Q_sigma}

We will construct $Q_{\sigma}^\ell: \bSigma^{\ell} \to \bSigma^{\ell+1}$ that is comprised of local operators on each face $\face^\ell$ and each aggregate $\agg^\ell$. Some of the notations here have been defined in Section~\ref{sec:edof}.

We first restrict ourselves to a face $\face^\ell = (\agg_i^\ell, \agg_j^\ell)$ on level $\ell$. Define a projection onto the PV trace:
\begin{equation}
\pi_{\sigma, \face^\ell}^{\ell,PV} := \bsigma_{\face^\ell}^{PV} Q_{\sigma, \face^\ell}^{\ell,PV}, \qquad \text{ where }  \;\; Q_{\sigma, \face^\ell}^{\ell,PV} :=  \left((\bq^{PV}_{\agg^\ell_i})^TD^\ell_{\agg^\ell_i, \face^\ell} \bsigma_{\face^\ell}^{PV}\right)^{-1} (\bq^{PV}_{\agg^\ell_i})^T D^\ell_{\agg^\ell_i, \face^\ell}.
\label{eq:def_proj_pv}
\end{equation}
Note that condition \eqref{eq:pv_condition} guarantees that $Q_{\sigma, \face^\ell}^{\ell,PV} $ is well-defined.

\begin{lemma}
If $(D_{\bgg^\ell}^\ell)^T \bone_{\bgg^\ell}^\ell = \bzero$ for any aggregate $\bgg^\ell$ of vertices on level $\ell$, then
\[
\left((\bq^{PV}_{\agg^\ell_i})^TD^\ell_{\agg^\ell_i, \face^\ell} \bsigma_{\face^\ell}^{PV}\right)^{-1} (\bq^{PV}_{\agg^\ell_i})^T D^\ell_{\agg^\ell_i, \face^\ell} = \left((\bq^{PV}_{\agg^\ell_j})^TD^\ell_{\agg^\ell_j, \face^\ell} \bsigma_{\face^\ell}^{PV}\right)^{-1} (\bq^{PV}_{\agg^\ell_j})^T D^\ell_{\agg^\ell_j, \face^\ell}.
\]
\label{lemma:proj_pv_unique}
\end{lemma}

{\it Proof.}
Let $\agg_{\face^\ell} := \agg_i^\ell \cup \agg_j^\ell$. Then
\begin{equation*}
\begin{split}
(\bzero_{\agg_{\face^\ell}})^T = (\bone^\ell_{\agg_{\face^\ell}})^T D^\ell_{\agg_{\face^\ell}} & = \begin{bmatrix}
(\bone^\ell_{\agg_i^\ell})^T D^\ell_{\agg_i^\ell} & (\bone^\ell_{\agg_i^\ell})^T D^\ell_{\agg_i^\ell, \face^\ell} + (\bone^\ell_{\agg_j^\ell})^T D^\ell_{\agg_j^\ell, \face^\ell} & (\bone^\ell_{\agg_j^\ell})^T D^\ell_{\agg_j^\ell}
\end{bmatrix}\\
& = \begin{bmatrix}
(\bzero_{\agg_i^\ell})^T & (\bone^\ell_{\agg_i^\ell})^T D^\ell_{\agg_i^\ell, \face^\ell} + (\bone^\ell_{\agg_j^\ell})^T D^\ell_{\agg_j^\ell, \face^\ell} & (\bzero_{\agg_j^\ell})^T
\end{bmatrix}.
\end{split}
\end{equation*}  
Thus, we have 
\begin{equation}
(\bone^\ell_{\agg_i^\ell})^T D^\ell_{\agg_i^\ell, \face^\ell} = -(\bone^\ell_{\agg_j^\ell})^T D^\ell_{\agg_j^\ell, \face^\ell}. 
\label{eq:proj_pv_unique}
\end{equation}
Recall from \eqref{eq:vertex_pv_vector} that $\bq^{PV}_{\agg^\ell_i} = \| \bone^\ell_{\agg_i^\ell} \|^{-1} \bone^\ell_{\agg_i^\ell}$, the lemma now follows easily from this and \eqref{eq:proj_pv_unique}.  \hfill$\square$ \\

\begin{remark}\label{rmk:proj_pv_unique}
We will show later in Lemma~\ref{lemma:D_null} that the hypothesis in Lemma~\ref{lemma:proj_pv_unique} is actually valid on all levels. Therefore, although we formally use one of the aggregates $\agg_i^\ell$ to define $Q_{\sigma,\face^\ell}^{\ell, PV}$ in \eqref{eq:def_proj_pv}, Lemma~\ref{lemma:proj_pv_unique} tells us that the definition of $Q_{\sigma,\face^\ell}^{\ell, PV}$ is in fact independent of the choice of $\agg_i^\ell$ or $\agg_j^\ell$.\\
\end{remark}
\\
Next, we will construct a projection for ``non-PV" traces. Without loss of generality, it is assumed that,
\begin{equation}
(\bq^{PV}_{\agg^\ell_i})^TD^\ell_{\agg^\ell_i, \face^\ell} \bsigma_{\face^\ell} = 0, \qquad \forall\, \bsigma_{\face^\ell} \not= \bsigma_{\face^\ell}^{PV}.
\label{eq:trace_orthogonality}
\end{equation}
Indeed, \eqref{eq:trace_orthogonality} can always be achieved by replacing 
$\bsigma_{\face^\ell}$ with $\bsigma_{\face^\ell} - \left( Q_{\sigma, \face^\ell}^{\ell,PV} \bsigma_{\face^\ell} \right) \bsigma_{\face^\ell}^{PV}$ for $\bsigma_{\face^\ell} \not= \bsigma_{\face^\ell}^{PV}$. Recall that $P_{\sigma, \face^\ell}^{\ell,NPV}$ defined in \eqref{eq:trace_decomposition} is the matrix whose columns are the non-PV traces. Now define a projection $\pi_{\sigma, \face^\ell}^{\ell,NPV} := P_{\sigma, \face^\ell}^{\ell,NPV} Q_{\sigma, \face^\ell}^{\ell,NPV}$ onto the span of non-PV traces, where
\begin{equation}
Q_{\sigma, \face^\ell}^{\ell,NPV} :=  \left((P_{\sigma, \face^\ell}^{\ell,NPV})^TP_{\sigma, \face^\ell}^{\ell,NPV}\right)^{-1}(P_{\sigma, \face^\ell}^{\ell,NPV})^T (I-\bsigma_{\face^\ell}^{PV} Q_{\sigma, \face^\ell}^{PV} ).
\label{eq:def_proj_npv}
\end{equation}
Combining $\pi_{\sigma, \face^\ell}^{\ell,PV}$ and $\pi_{\sigma, \face^\ell}^{\ell,NPV}$, we obtain a projection for all trace extensions associated with $\face^\ell$:
\[
\pi_{\sigma, \face^\ell}^{\ell, T} := \pi_{\sigma, \face^\ell}^{\ell,PV} + \pi_{\sigma, \face^\ell}^{\ell,NPV} = P_{\sigma, \face^\ell}^{\ell, T} \;Q_{\sigma, \face^\ell}^{\ell, T},
\]
where $P_{\sigma, \face^\ell}^{\ell,T}$ is defined in \eqref{eq:trace_decomposition} and $Q_{\sigma, \face^\ell}^{\ell, T} := {\footnotesize
\begin{tikzpicture}[baseline=-\the\dimexpr\fontdimen22\textfont2\relax ]
\matrix (m)[matrix of math nodes,left delimiter={[},right delimiter={]}]
{
\,\,Q_{\sigma, \face^\ell}^{\ell,PV}\,\, \\
 \hspace{1mm} \\
Q_{\sigma, \face^\ell}^{\ell,NPV}  \\
};
\begin{pgfonlayer}{myback}
\fhighlight[blue!30]{m-1-1}{m-1-1}
\fhighlight[red!30]{m-3-1}{m-3-1}
\end{pgfonlayer} 
\end{tikzpicture} }.
$

\vspace{3mm}

\begin{lemma}\label{lemma:trace_projection}
$Q_{\sigma, \face^\ell}^{\ell, T}$ is a left inverse of $P_{\sigma, \face^\ell}^{\ell,T}$. That is,
$Q_{\sigma, \face^\ell}^{\ell, T} \, P_{\sigma, \face^\ell}^{\ell,T} = I_{\sigma, \face^\ell}^{\ell+1}$,
where $I_{\sigma, \face^\ell}^{\ell+1}$ is the identity matrix of size number of traces in $Tr(\face^\ell)$.\\
\end{lemma}

{\it Proof.}
Recall from \eqref{eq:trace_decomposition} that columns of $P_{\sigma, \face^\ell}^{\ell,T}$ can be split into PV and non-PV traces
\[
P_{\sigma, \face^\ell}^{\ell,T} = {\footnotesize
\begin{tikzpicture}[baseline=-\the\dimexpr\fontdimen22\textfont2\relax ]
\matrix (m)[matrix of math nodes,left delimiter={[},right delimiter={]}]
{
\bsigma_{\face^\ell}^{PV} & \hspace{1mm} & P_{\sigma, \face^\ell}^{\ell,NPV}  \\
};
\begin{pgfonlayer}{myback}
\fhighlight[blue!30]{m-1-1}{m-1-1}
\fhighlight[red!30]{m-1-3}{m-1-3}
\end{pgfonlayer}
\end{tikzpicture} }.
\] 
It follows clearly from definition \eqref{eq:def_proj_pv} and condition \eqref{eq:trace_orthogonality} that
\begin{equation}
Q_{\sigma, \face^\ell}^{\ell,PV} \, P_{\sigma, \face^\ell}^{\ell,T} 
= {\footnotesize
\begin{tikzpicture}[baseline=-\the\dimexpr\fontdimen22\textfont2\relax ]
\matrix (m)[matrix of math nodes,left delimiter={[},right delimiter={]}]
{
Q_{\sigma, \face^\ell}^{\ell,PV} \, \bsigma_{\face^\ell}^{PV} &  \hspace{1mm} & Q_{\sigma, \face^\ell}^{\ell,PV} \,P_{\sigma, \face^\ell}^{\ell,NPV}  \\
};
\begin{pgfonlayer}{myback}
\fhighlight[blue!30]{m-1-1}{m-1-1}
\fhighlight[red!30]{m-1-3}{m-1-3}
\end{pgfonlayer}
\end{tikzpicture} }
= {\small
\begin{tikzpicture}[baseline=-\the\dimexpr\fontdimen22\textfont2\relax ]
\matrix (m)[matrix of math nodes,left delimiter={[},right delimiter={]}]
{
1 & \hspace{1mm} &  0 & \cdots & 0  \\
};
\begin{pgfonlayer}{myback}
\fhighlight[blue!30]{m-1-1}{m-1-1}
\fhighlight[red!30]{m-1-3}{m-1-5}
\end{pgfonlayer}
\end{tikzpicture} },
\label{eq:proj_trace_proof1}
\end{equation}
and therefore
\begin{equation}
(I-\bsigma_{\face^\ell}^{PV} Q_{\sigma, \face^\ell}^{PV} ) \, P_{\sigma, \face^\ell}^{\ell,T} 
= P_{\sigma, \face^\ell}^{\ell,T} - \bsigma_{\face^\ell}^{PV} {\small
\begin{tikzpicture}[baseline=-\the\dimexpr\fontdimen22\textfont2\relax ]
\matrix (m)[matrix of math nodes,left delimiter={[},right delimiter={]}]
{
1 & \hspace{1mm} &  0 & \cdots & 0  \\
};
\begin{pgfonlayer}{myback}
\fhighlight[blue!30]{m-1-1}{m-1-1}
\fhighlight[red!30]{m-1-3}{m-1-5}
\end{pgfonlayer}
\end{tikzpicture} }
= {\footnotesize
\begin{tikzpicture}[baseline=-\the\dimexpr\fontdimen22\textfont2\relax ]
\matrix (m)[matrix of math nodes,left delimiter={[},right delimiter={]}]
{
\bzero^\ell_{\sigma, \face^\ell} &  \hspace{1mm} & P_{\sigma, \face^\ell}^{\ell,NPV}   \\
};
\begin{pgfonlayer}{myback}
\fhighlight[blue!30]{m-1-1}{m-1-1}
\fhighlight[red!30]{m-1-3}{m-1-3}
\end{pgfonlayer}
\end{tikzpicture} }.
\label{eq:proj_trace_proof2}
\end{equation}
By \eqref{eq:def_proj_npv} and \eqref{eq:proj_trace_proof2},
\begin{equation}
\begin{split}
Q_{\sigma, \face^\ell}^{\ell,NPV} \, P_{\sigma, \face^\ell}^{\ell,T} 
 & = {\footnotesize
\begin{tikzpicture}[baseline=-\the\dimexpr\fontdimen22\textfont2\relax ]
\matrix (m)[matrix of math nodes,left delimiter={[},right delimiter={]}]
{
\bzero^{\ell+1, NPV}_{\sigma, \face^\ell}  \pgfmatrixnextcell  \hspace{1mm} \pgfmatrixnextcell \left( (P_{\sigma, \face^\ell}^{\ell,NPV})^TP_{\sigma, \face^\ell}^{\ell,NPV} \right)^{-1}(P_{\sigma, \face^\ell}^{\ell,NPV})^TP_{\sigma, \face^\ell}^{\ell,NPV}  \\
};
\begin{pgfonlayer}{myback}
\fhighlight[blue!30]{m-1-1}{m-1-1}
\fhighlight[red!30]{m-1-3}{m-1-3}
\end{pgfonlayer}
\end{tikzpicture} }\\
 & = {\footnotesize
\begin{tikzpicture}[baseline=-\the\dimexpr\fontdimen22\textfont2\relax ]
\matrix (m)[matrix of math nodes,left delimiter={[},right delimiter={]}]
{
\bzero^{\ell+1, NPV}_{\sigma, \face^\ell} \pgfmatrixnextcell  \hspace{1mm} \pgfmatrixnextcell I^{\ell+1, NPV}_{\sigma, \face^\ell}   \\
};
\begin{pgfonlayer}{myback}
\fhighlight[blue!30]{m-1-1}{m-1-1}
\fhighlight[red!30]{m-1-3}{m-1-3}
\end{pgfonlayer}
\end{tikzpicture} },
\end{split}
\label{eq:proj_trace_proof3}
\end{equation}
where $I_{\sigma, \face^\ell}^{\ell+1, NPV}$ is the identity matrix.
Combining \eqref{eq:proj_trace_proof1} and \eqref{eq:proj_trace_proof3},
\begin{equation*}
Q_{\sigma, \face^\ell}^{\ell, T} \, P_{\sigma, \face^\ell}^{\ell,T} 
= {\footnotesize
\begin{tikzpicture}[baseline=-\the\dimexpr\fontdimen22\textfont2\relax ]
\matrix (m)[matrix of math nodes,left delimiter={[},right delimiter={]}]
{
Q_{\sigma, \face^\ell}^{\ell,PV} \, P_{\sigma, \face^\ell}^{\ell, T}  \\
 \hspace{1mm} \\
Q_{\sigma, \face^\ell}^{\ell,NPV} \, P_{\sigma, \face^\ell}^{\ell, T}  \\
};
\begin{pgfonlayer}{myback}
\fhighlight[blue!30]{m-1-1}{m-1-1}
\fhighlight[red!30]{m-3-1}{m-3-1}
\end{pgfonlayer}
\end{tikzpicture} }
= {\footnotesize
\begin{tikzpicture}[baseline=-\the\dimexpr\fontdimen22\textfont2\relax ]
\matrix (m)[matrix of math nodes,left delimiter={[},right delimiter={]}]
{
1  \\
 \hspace{1mm} \\
& \hspace{1mm} & I^{\ell+1, NPV}_{\sigma, \face^\ell}  \\
};
\begin{pgfonlayer}{myback}
\fhighlight[blue!30]{m-1-1}{m-1-1}
\fhighlight[red!30]{m-3-3}{m-3-3}
\end{pgfonlayer}
\end{tikzpicture} }.
\end{equation*}
     
\vspace{-7mm}   
\hfill$\square$

\vspace{6mm}

Next, for each aggregate $\agg^\ell$, define a projection $\pi_{\sigma, \agg^\ell}^{\ell, B} := P_{\sigma, \agg^\ell}^{\ell, B}\,Q_{\sigma, \agg^\ell}^{\ell, B} $ to the bubble space in $\agg^\ell$, where
\begin{equation}
Q_{\sigma, \agg^\ell}^{\ell, B} := (P_{u, \agg^\ell}^{\ell, NPV})^T D^\ell_{\agg^\ell}.
\label{eq:bubble_project}
\end{equation}
We recall that $P_{u, \agg^\ell}^{\ell, NPV}$ defined in \eqref{eq:vertex_basis_decomposition} contains the coarse vertex space basis vectors in $\agg^\ell$ that are orthogonal to $\bone_{\ell^\ell}^\ell$. By \eqref{eq:bubble_divergence} and the fact that columns of $P_{u, \agg^\ell}^{\ell,NPV}$ are orthonormal, we have
\begin{equation}
Q_{\sigma, \agg^\ell}^{\ell, B} \, P_{\sigma, \agg^\ell}^{\ell, B} =  (P_{u, \agg^\ell}^{\ell, NPV})^T \, (D^\ell_{\agg^\ell}P_{\sigma, \agg^\ell}^{\ell, B}) = (P_{u, \agg^\ell}^{\ell, NPV})^T \, P_{u, \agg^\ell}^{\ell, NPV} = I_{\sigma, \agg^\ell}^{\ell+1, B}
\label{eq:bubble_projection}
\end{equation}
where $I_{\sigma, \agg^\ell}^{\ell+1, B}$ is the identity. Finally, with $D_{E^{\ell+1}}^\ell$ being the column restriction of $D^\ell$ to $\displaystyle\bigoplus_{\face^\ell\in E^{\ell+1}}\bSigma^{\ell}(\face^\ell)$, we define the projection to the coarse edge space $\pi_\sigma^\ell := P_\sigma^\ell\,Q_\sigma^\ell$, where
\begin{equation}
Q_\sigma^\ell := {\footnotesize
\begin{tikzpicture}[baseline=-\the\dimexpr\fontdimen22\textfont2\relax ]
\matrix (m)[matrix of math nodes,left delimiter={[},right delimiter={]}]
{
\;\;\ddots\;\; & {\color{blue!30}\ddots} & & & {\color{blue!0}\ddots} & & & {\color{blue!0}\ddots}  \\
{\color{red!30}\;\;\ddots\;\;} & \ddots & &  \\
& & \,\,Q_{\sigma, \face^\ell}^{\ell,PV}\,\, & \\
& & Q_{\sigma, \face^\ell}^{\ell,NPV}  & \\
& & & \;\;\ddots\;\; & {\color{blue!30}\ddots} & & & \\
\hspace{1mm} & & & {\color{red!30}\;\;\ddots\;\;} & \ddots & & & \hspace{1mm} \\
& & & & & \ddots & & {\color{blue!0}\ddots} \\
\hspace{1mm} & & & & & & Q_{\sigma, \agg^\ell}^{\ell, B} \\
\hspace{1mm} & & & & {\color{blue!0}\ddots} & \hspace{1mm} & & \ddots \\
};
\draw[black!30,thick] (m-6-1.south west) -- (m-6-8.south east);
\draw[black!30,thick] (m-1-5.north east) -- (m-9-5.south east);
\begin{pgfonlayer}{myback}
\fhighlight[blue!30]{m-1-1}{m-1-2}
\fhighlight[red!30]{m-2-1}{m-2-2}
\fhighlight[blue!30]{m-3-3}{m-3-3}
\fhighlight[red!30]{m-4-3}{m-4-3}
\fhighlight[blue!30]{m-5-4}{m-5-5}
\fhighlight[red!30]{m-6-4}{m-6-5}
\fhighlight[green!30]{m-7-6}{m-7-6}
\fhighlight[green!30]{m-8-7}{m-8-7}
\fhighlight[green!30]{m-9-8}{m-9-8}
\node (a) at (m-1-8.north) [right=20pt]{};
\node (b) at (m-6-8.south) [right=20pt]{};
\node (c) at (m-7-8.north) [right=20pt]{};
\node (d) at (m-9-8.south) [right=20pt]{};
\draw [decorate, decoration={brace, amplitude=10pt}] (a) -- (b) 
node[midway, right=10pt] {$\displaystyle\bigoplus_{\face^\ell\in E^{\ell+1}}\bSigma^{\ell+1}(\face^\ell)$};
\draw [decorate, decoration={brace, amplitude=10pt}] (c) -- (d) 
node[midway, right=10pt] {$\displaystyle\bigoplus_{\agg^\ell\in V^{\ell+1}}\bSigma^{\ell+1}(\agg^\ell)$};
\node (e) at (m-9-1.south west) [below=10pt]{};
\node (f) at (m-9-5.south east) [below=10pt]{};
\node (g) at (m-9-6.south west) [below=10pt]{};
\node (h) at (m-9-8.south east) [below=10pt]{};
\draw [decorate, decoration={brace, amplitude=10pt, mirror}] (e) -- (f) 
node[midway, below=10pt] {$\displaystyle\bigoplus_{\face^\ell\in E^{\ell+1}}\bSigma^{\ell}(\face^\ell)$};
\draw [decorate, decoration={brace, amplitude=10pt, mirror}] (g) -- (h) 
node[midway, below=10pt] {$\displaystyle\bigoplus_{\agg^\ell\in V^{\ell+1}}\bSigma^{\ell}(\agg^\ell)$};
\node (i) at (m-8-1) [right=-8pt]{${\scriptsize \begin{bmatrix} \ddots  \\ & \hspace{-3mm} \left(P_{u, \agg^\ell}^{\ell,NPV}\right)^T \hspace{-3.5mm} \\ & & \ddots \end{bmatrix}}
D^\ell_{E^{\ell+1}}$};
\end{pgfonlayer}
\end{tikzpicture} }.
\label{eq:Qedge_structure}
\end{equation}
Now we are ready to prove Proposition~\ref{prop:projection}.\vspace{3mm} \\

\noindent{\bf Proof of Proposition~\ref{prop:projection}:} 
Recall from \eqref{eq:Pedge_structure} that $P_\sigma^\ell$ has the structure
\[
P_\sigma^\ell = {\scriptsize
\begin{tikzpicture}[baseline=-\the\dimexpr\fontdimen22\textfont2\relax ]
\matrix (m)[matrix of math nodes,left delimiter={[},right delimiter={]}]
{
\ddots & {\color{red!30}\ddots} & & & & {\color{blue!0}\ddots} & & & {\color{blue!0}\ddots}  \\
{\color{blue!30}\ddots} & \ddots & & &  \\
& & \bsigma_{\face^\ell{\color{blue!30},i}}^{\ell,PV} & P_{\sigma, \face^\ell}^{\ell,NPV} & \\
& & & & \ddots & {\color{red!30}\ddots} & & & \\
\hspace{1mm} & & & & {\color{blue!30}\ddots} & \ddots & & & \hspace{1mm} \\
& & & & & & \ddots & & {\color{blue!0}\ddots} \\
& & \hspace{1mm}  & & & & & P_{\sigma, \agg^\ell}^{\ell, B} \\
\hspace{1mm} & & & & & {\color{blue!0}\ddots} & \hspace{1mm} & & \ddots \\
};
\draw[black!30,thick] (m-5-1.south west) -- (m-5-9.south east);
\draw[black!30,thick] (m-1-6.north east) -- (m-8-6.south east);
\begin{pgfonlayer}{myback}
\fhighlight[blue!30]{m-1-1}{m-2-1}
\fhighlight[red!30]{m-1-2}{m-2-2}
\fhighlight[blue!30]{m-3-3}{m-3-3}
\fhighlight[red!30]{m-3-4}{m-3-4}
\fhighlight[blue!30]{m-4-5}{m-5-5}
\fhighlight[red!30]{m-4-6}{m-5-6}
\fhighlight[green!30]{m-6-7}{m-6-7}
\fhighlight[green!30]{m-7-8}{m-7-8}
\fhighlight[green!30]{m-8-9}{m-8-9}
\node (i) at (m-7-3) [right=0pt]{\normalsize$P_{\sigma}^{\ell,E}$};
\end{pgfonlayer}
\end{tikzpicture}. }
\]
By \eqref{eq:Qedge_structure}, Lemma~\ref{lemma:trace_projection} and \eqref{eq:bubble_projection}, we have
\begin{equation}
Q_\sigma^\ell P_\sigma^\ell = {\scriptsize
\begin{tikzpicture}[baseline=-\the\dimexpr\fontdimen22\textfont2\relax ]
\matrix (m)[matrix of math nodes,left delimiter={[},right delimiter={]}]
{
\ddots & {\color{red!0}\ddots} & & & & {\color{blue!0}\ddots} & & & {\color{blue!0}\ddots}  \\
{\color{blue!0}\ddots} & \ddots & & &  \\
& & 1 &  & \\
& & & I_{\sigma, \face^\ell}^{\ell+1,NPV} & \\
& & & & \ddots & {\color{red!0}\ddots} & & & \\
\hspace{1mm} & & & & {\color{blue!0}\ddots} & \ddots & & & \hspace{1mm} \\
& & & & & & \ddots & & {\color{blue!0}\ddots} \\
& & \hspace{1mm} & & & & & I_{\sigma, \agg^\ell}^{\ell+1, B} \\
\hspace{1mm} & & & & & {\color{blue!0}\ddots} & \hspace{1mm} & & \ddots \\
};
\draw[black!30,thick] (m-6-1.south west) -- (m-6-9.south east);
\draw[black!30,thick] (m-1-6.north east) -- (m-9-6.south east);
\begin{pgfonlayer}{myback}
\fhighlight[blue!30]{m-1-1}{m-1-1}
\fhighlight[red!30]{m-2-2}{m-2-2}
\fhighlight[blue!30]{m-3-3}{m-3-3}
\fhighlight[red!30]{m-4-4}{m-4-4}
\fhighlight[blue!30]{m-5-5}{m-5-5}
\fhighlight[red!30]{m-6-6}{m-6-6}
\fhighlight[green!30]{m-7-7}{m-7-7}
\fhighlight[green!30]{m-8-8}{m-8-8}
\fhighlight[green!30]{m-9-9}{m-9-9}
\node (i) at (m-8-3) [right=0pt]{{\normalsize $\left( Q_\sigma^\ell P_\sigma^\ell \right)_{21}$}};
\end{pgfonlayer}
\end{tikzpicture} }.
\end{equation}
Hence, it suffices to show that $\left( Q_\sigma^\ell P_\sigma^\ell \right)_{21}$, the (2,1)-block of the product $Q_\sigma^\ell P_\sigma^\ell$, is a zero matrix:
\begin{equation}
\left( Q_\sigma^\ell P_\sigma^\ell \right)_{21} = 
{\footnotesize \begin{bmatrix} \ddots  \\ & (P_{u, \agg^\ell}^{\ell,NPV})^T \\ & & \ddots \end{bmatrix}}
D^\ell_{E^{\ell+1}} 
{\footnotesize \begin{bmatrix}
\ddots \\ & P_{\sigma, \face^\ell}^{\ell, T} \\ & & \ddots
\end{bmatrix}}
+ {\footnotesize \begin{bmatrix}
\ddots \\ & Q_{\sigma, \agg^\ell}^{\ell} \\ & & \ddots
\end{bmatrix}}P_{\sigma}^{\ell,E} = 0.
\label{eq:extension_orthogonality}
\end{equation}
We will focus on the column of $\left( Q_\sigma^\ell P_\sigma^\ell \right)_{21}$ that is associated with some trace $\bsigma_{\face^\ell}$, where $\face^\ell = (\agg_i^\ell, \agg_j^\ell)$. The only possible nonzeros in this columns are the rows associated with the bubbles in $\agg^\ell_i$ and $\agg^\ell_j$. For $\agg^\ell = \agg^\ell_i$ or $\agg^\ell_j$, let $\bsigma_{\agg^\ell}^{E}$ be the extension of $\bsigma_{\face^\ell}$ into $\edgespace^\ell(\agg^\ell)$ that is defined by \eqref{eq:extension}. Then the submatrix of $\left( Q_\sigma^\ell P_\sigma^\ell \right)_{21}$ restricted to the column associated with $\bsigma_{\face^\ell}$ and the rows associated with bubbles in $\agg^\ell$ is
\[
(P_{u, \agg^\ell}^{\ell,NPV})^T D^\ell_{\agg^\ell, \face^\ell} \bsigma_{\face^\ell} \;+\; Q_{\sigma, \agg^\ell}^{\ell} \bsigma_{\agg^\ell}^{E}.
\]
By the second block of equations in \eqref{eq:extension}, we know that
\begin{equation*}
\begin{split}
 (P_{u, \agg^\ell}^{\ell,NPV})^T D^\ell_{\agg^\ell, \face^\ell} \bsigma_{\face^\ell} \;+\; Q_{\sigma, \agg^\ell}^{\ell} \bsigma_{\agg^\ell}^{E}  
& = (P_{u, \agg^\ell}^{\ell,NPV})^T \left(D^\ell_{\agg^\ell, \face^\ell} \bsigma_{\face^\ell} + D_{\agg^\ell}^{\ell} \bsigma_{\agg^\ell}^{E} \right)  \\
& = (P_{u, \agg^\ell}^{\ell,NPV})^T \left(c_\agg^{\ell} \bq^{PV}_{\agg^\ell} \right) = 0. \\
\end{split}
\end{equation*}
Thus, \eqref{eq:extension_orthogonality} and therefore the proposition follow.  {\hfill$\square$}

\vspace{5mm}

\noindent The rest of this section are some lemmas that lead to the key commutative property \eqref{eq:commutativity}.\\

\begin{lemma}
On level $\ell$, if $(D_{\bgg^\ell}^\ell)^T \bone_{\bgg^\ell}^\ell = \bzero$ for any aggregate $\bgg^\ell$, then
\begin{equation}
D^{\ell+1} Q_\sigma^\ell = Q_u^\ell D^\ell.
\label{eq:commutativity2}
\end{equation}
\label{lemma:commutativity}
\end{lemma}

{\it Proof.} Since $D^{\ell+1} = (P_u^\ell)^T D^\ell P_\sigma^\ell = Q_u^\ell D^\ell P_\sigma^\ell$, proving \eqref{eq:commutativity2} is equivalent to prove that,
\begin{equation}
Q_u^\ell \, D^\ell \left(P_\sigma^\ell Q_\sigma^\ell \bsigma \right) = Q_u^\ell \, D^\ell \bsigma, \qquad \forall\,\bsigma\in\bSigma^\ell.
\label{eq:commutativity_proof1}
\end{equation}
Our goal is to show \eqref{eq:commutativity_proof1} by considering rows corresponding to one aggregate at a time. In the aggregate $\agg^\ell$, $Q_u^\ell$ is simply (cf. \eqref{eq:vertex_basis_decomposition})
\[
(P_{u,\agg^\ell}^\ell)^T = \begin{bmatrix}
(\bq_{\agg^\ell}^{PV})^T \vspace{3mm}  \\ (P_{u, \agg^\ell}^{NPV})^T
\end{bmatrix}.
\]
Consider partitions of $P_\sigma^\ell$ and $Q_\sigma^\ell$ based on trace extensions and bubbles (cf. \eqref{eq:Pedge_structure} and \eqref{eq:Qedge_structure})
\[
P_\sigma^\ell = \begin{bmatrix}
P_\sigma^{\ell,TE} &  P_{\sigma}^{\ell, B}
\end{bmatrix},
\qquad
Q_\sigma^\ell = \begin{bmatrix}
Q_\sigma^{\ell,TE} \vspace{2mm} \\  Q_{\sigma}^{\ell, B}
\end{bmatrix}.
\]
By construction, the restriction of $D^\ell \bsigma_i$ in $\agg^\ell$ is in $\Span{\bq_{\agg^\ell}^{PV}}$ if $\bsigma_i$ is a trace extension; or it is in $\Range{P_{u, \agg}^{NPV}}$ if $\bsigma_i$ is a bubble. This observation and the orthogonality of the columns of $P_{u,\agg^\ell}^\ell$ imply that 
\begin{equation}
\begin{split}
(P_{u, \agg^\ell}^{\ell})^T \left((D^\ell P_\sigma^\ell Q_\sigma^\ell \bsigma)|_{\agg^\ell}\right) = & \begin{bmatrix}
(\bq_{\agg^\ell}^{PV})^T \vspace{3mm}  \\ (P_{u, \agg^\ell}^{NPV})^T
\end{bmatrix} \left(D^\ell P_{\sigma}^{\ell,TE} Q_{\sigma}^{\ell,TE} \bsigma + D^\ell P_{\sigma}^{\ell,B} Q_{\sigma}^{\ell,B} \bsigma\right)|_{\agg^\ell} \\
= & \begin{bmatrix}
(\bq_{\agg^\ell}^{PV})^T \left(D^\ell P_{\sigma}^{\ell,TE} Q_{\sigma}^{\ell,TE} \bsigma\right)|_{\agg^\ell} \vspace{3mm}
 \\  (P_{u, \agg^\ell}^{NPV})^T \left(D^\ell P_{\sigma}^{\ell,B} Q_{\sigma}^{\ell,B} \bsigma\right)|_{\agg^\ell}
\end{bmatrix}.
\end{split}
\label{eq:commutativity_proof_D_orthogonality}
\end{equation}

For the bubble part, we know from the structure of $Q_\sigma^\ell$ \eqref{eq:Qedge_structure} that $Q_{\sigma}^{\ell, B}$ has the form
\begin{equation}
Q_{\sigma}^{\ell, B} = \begin{bmatrix}\vspace{-2mm}\\
\;{\scriptsize \begin{bmatrix} \ddots  \\ & (P_{u, \agg^\ell}^{\ell,NPV})^T \\ & & \ddots \end{bmatrix}
 D^\ell_{E^{\ell+1}}}  & & {\scriptsize \begin{bmatrix} \ddots  \\ & Q_{\sigma, \agg^\ell}^{\ell, B} \\ & & \ddots \end{bmatrix}}
\\ 
\vspace{-2mm}\\
\end{bmatrix}
\label{eq:commutativity_proof_bubble1}
\end{equation}
By \eqref{eq:bubble_project},
\begin{equation}
\begin{bmatrix} \ddots  \\ & Q_{\sigma, \agg^\ell}^{\ell, B} \\ & & \ddots \end{bmatrix}
= \begin{bmatrix} \ddots  \\ & (P_{u, \agg^\ell}^{\ell,NPV})^T D_{\agg^\ell}^\ell \\ & & \ddots \end{bmatrix}
= \begin{bmatrix} \ddots  \\ & (P_{u, \agg^\ell}^{\ell,NPV})^T \\ & & \ddots \end{bmatrix}\begin{bmatrix} \ddots  \\ & D_{\agg^\ell}^\ell \\ & & \ddots \end{bmatrix}
\label{eq:commutativity_proof_bubble2}
\end{equation}
Substituting \eqref{eq:commutativity_proof_bubble2} into \eqref{eq:commutativity_proof_bubble1}, we can see that
\begin{equation}
Q_{\sigma}^{\ell, B} = \begin{bmatrix} \ddots  \\ & (P_{u, \agg^\ell}^{\ell,NPV})^T \\ & & \ddots \end{bmatrix}
\begin{bmatrix} 
\vspace{-2mm}\\
\;D^\ell_{E^{\ell+1}}  & {\scriptsize \begin{bmatrix} \ddots  \\ & D_{\agg^\ell}^\ell \\ & & \ddots \end{bmatrix}}
\\ 
\vspace{-2mm}\\
\end{bmatrix}
= \begin{bmatrix} \ddots  \\ & (P_{u, \agg^\ell}^{\ell,NPV})^T \\ & & \ddots \end{bmatrix}
D^\ell.
\label{eq:commutativity_proof2}
\end{equation}
Hence, by \eqref{eq:bubble_divergence} and \eqref{eq:commutativity_proof2}, 
\begin{equation}
\begin{split}
& \hspace{4.5mm} (P_{u, \agg^\ell}^{NPV})^T \left(D^\ell P_{\sigma}^{\ell,B} Q_{\sigma}^{\ell,B} \bsigma\right)|_{\agg^\ell} 
= (P_{u, \agg^\ell}^{NPV})^T \left((D^\ell P_{\sigma}^{\ell,B}) (Q_{\sigma}^{\ell,B} \bsigma) \right)|_{\agg^\ell} \\
& = (P_{u, \agg^\ell}^{NPV})^T \left.\left(\begin{bmatrix} \ddots  \\ & P_{u, \agg^\ell}^{\ell,NPV} \\ & & \ddots \end{bmatrix} \begin{bmatrix} \ddots  \\ & (P_{u, \agg^\ell}^{\ell,NPV})^T \\ & & \ddots \end{bmatrix} D^\ell \bsigma \right)\right|_{\agg^\ell} \\
& = (P_{u, \agg^\ell}^{NPV})^T P_{u, \agg^\ell}^{NPV} (P_{u, \agg^\ell}^{NPV})^T \left(D^\ell \bsigma \right)|_{\agg^\ell} \\
& = (P_{u, \agg^\ell}^{NPV})^T \left(D^\ell \bsigma \right)|_{\agg^\ell}.
\end{split}
\label{eq:commutativity_proof_bubble}
\end{equation}

For the trace extension part, first note that by Lemma~\ref{lemma:proj_pv_unique} and the hypothesis of the current lemma, we can take $Q_{\sigma, \face^\ell}^{\ell,PV}$ in the definition of PV-trace projections for faces $\face^\ell \subseteq \bdr E^{\ell+1}(\agg^\ell)$ to be (cf. Remark~\ref{rmk:proj_pv_unique})
\begin{equation}
Q_{\sigma, \face^\ell}^{\ell,PV} = \left(\bq^{PV}_{\agg^\ell})^T D^\ell_{\agg^\ell, \face^\ell} \bsigma_{\face^\ell}^{PV} \right)^{-1}(\bq^{PV}_{\agg^\ell})^T D^\ell_{\agg^\ell, \face^\ell}.
\label{eq:proj_pv_choice}
\end{equation}
Next, note also that the only relevant coarse edge dofs for the aggregate $\agg^\ell$ are the trace extensions $\bsigma_{\face^\ell}^{TE}$ associated with boundary faces of $\agg^\ell$, $\face^\ell \subseteq \bdr E^{\ell+1}(\agg^\ell)$. Moreover, by definition \eqref{eq:extension}, each $\bsigma_{\face^\ell}^{TE}$ satisfies
\begin{equation}
\left( D^\ell \bsigma_{\face^\ell}^{TE} \right)|_{\agg^\ell} = c_{\agg^\ell} \bq_{\agg^\ell}^{PV}.
\label{eq:trace_extension_divergence}
\end{equation}
Because of \eqref{eq:extension_rhs_constant} and \eqref{eq:trace_orthogonality}, the constant $c_{\agg^\ell}$ is exactly
\begin{equation}
c_{\agg^\ell} = \left\{ \begin{array} {ll}
c_{\agg^\ell, \face^\ell}^{PV} := \big(\bq^{PV}_{\agg^\ell}\big)^T D^\ell_{\agg^\ell,\face^\ell} \,\bsigma_{\face^\ell}^{PV} .
  & \quad \text{ if $\bsigma_{\face^\ell}^{TE}$ is the extension of $\bsigma_{\face^\ell}^{PV}$,} \\
0  & \quad \text{ otherwise.}
\end{array}
\right.
\label{eq:explicit_extension_rhs_constant}
\end{equation}
So, only the PV-trace projections associated with $\face^\ell \subseteq \bdr E^{\ell+1}(\agg^\ell)$ contribute to $\left(D^\ell P_{\sigma}^{\ell,TE} Q_{\sigma}^{\ell,TE} \bsigma\right)|_{\agg^\ell}$.
This observation together with \eqref{eq:trace_extension_divergence}, \eqref{eq:explicit_extension_rhs_constant}, \eqref{eq:vertex_pv_vector}, \eqref{eq:proj_pv_choice}, and the hypothesis of the lemma imply that
\begin{equation}
\begin{split}
& \hspace{4.5mm} (\bq_{\agg^\ell}^{PV})^T \left(D^\ell P_{\sigma}^{\ell,TE} Q_{\sigma}^{\ell,TE} \bsigma\right)|_{\agg^\ell} = (\bq_{\agg^\ell}^{PV})^T \left( (D^\ell P_{\sigma}^{\ell,TE}) (Q_{\sigma}^{\ell,TE} \bsigma) \right)|_{\agg^\ell} \\
& =  (\bq_{\agg^\ell}^{PV})^T \sum_{\face^\ell\subseteq\bdr E^{\ell+1}(\agg^\ell)} \left( c_{\agg^\ell, \face^\ell}^{PV}\bq_{\agg^\ell}^{PV} \right) \left( Q_{\sigma, \face^\ell}^{\ell,PV} \bsigma|_{\face^\ell} \right) \\
& = \sum_{\face^\ell\subseteq\bdr E^{\ell+1}(\agg^\ell)} \left( c_{\agg^\ell, \face^\ell}^{PV} \left\| \bq_{\agg^\ell}^{PV} \right\|^2 \right) \left( \left(\bq^{PV}_{\agg^\ell})^T D^\ell_{\agg^\ell, \face^\ell} \bsigma_{\face^\ell}^{PV} \right)^{-1}(\bq^{PV}_{\agg^\ell})^T D^\ell_{\agg^\ell, \face^\ell}\bsigma|_{\face^\ell}  \right) \\
& = \sum_{\face^\ell\subseteq\bdr E^{\ell+1}(\agg^\ell)} (\bq^{PV}_{\agg^\ell})^T D^\ell_{\agg^\ell, \face^\ell}\bsigma|_{\face^\ell} \\
& = (\bq^{PV}_{\agg^\ell})^T D^\ell_{\agg^\ell}\bsigma|_{\agg^\ell} + \sum_{\face^\ell\subseteq\bdr E^{\ell+1}(\agg^\ell)} (\bq^{PV}_{\agg^\ell})^T D^\ell_{\agg^\ell, \face^\ell}\bsigma|_{\face^\ell} \\
& = (\bq_{\agg^\ell}^{PV})^T \left(D^\ell \bsigma \right)|_{\agg^\ell}.
\end{split}
\label{eq:commutativity_proof_trace_extension}
\end{equation}
Finally, \eqref{eq:commutativity} follows from \eqref{eq:commutativity_proof1}, \eqref{eq:commutativity_proof_bubble}, and \eqref{eq:commutativity_proof_trace_extension}.              {\hfill$\square$}\\
\\

\begin{lemma}
On level $\ell\ge 0$,  if $(D_{\bgg^\ell}^\ell)^T \bone_{\bgg^\ell}^\ell = \bzero$ for any aggregate of vertices $\bgg^\ell$, then
\begin{equation}
\Span{\bone^\ell_{\agg^\ell_i}} \subseteq \Range{P_{u, \agg^\ell_i}^\ell}  \qquad \forall\, \agg^\ell_i\in V^{\ell+1}.
\label{eq:const_preserve}
\end{equation}
\label{lemma:const_preserve}
\end{lemma}
{\it Proof.} By the assumption and noticing that $\Lapl_{N(\agg_i^\ell)}$ is symmetric positive semi-definite, we know $\bone^\ell_{N(\agg_i^\ell)}$ is an eigenvector of \eqref{eq:eigen} corresponding to the smallest eigenvalue $\lambda = 0$. Recall that our algorithm picks the eigenvectors associated with the small eigenvalues, so $\bone^\ell_{N(\agg_i^\ell)}$ is always selected. After restricting $\bone^\ell_{N(\agg_i^\ell)}$ to $\agg_i^\ell$, we get $\bone^\ell_{\agg_i^\ell}$. As SVD only removes linear dependence, $\bone^\ell_{\agg_i^\ell}$ remains in the range of $P_{u, \agg^\ell_i}^\ell$.   \hfill$\square$\\


\begin{lemma}
On level $\ell\ge 0$, for any aggregate of vertices $\bgg^\ell$, we have $(D^\ell_{\bgg^\ell})^T \bone^\ell_{\bgg^\ell} = \bzero$. Furthermore, if $\bgg^\ell$ is connected, then $\Null{(D^\ell_{\bgg^\ell})^T} = \Span{\bone^\ell_{\bgg^\ell}}$.\vspace{3mm}
\label{lemma:D_null}
\end{lemma}


{\it Proof.}
We will show the lemma by induction on $\ell$. When $\ell = 0$, for any aggregate $\bgg^0$, $(D^0_{\bgg^0})^T$ is an incidence matrix of which each row has exactly 2 non-zero entries; one is 1 and the other is -1.  Clearly, the constant vector $\bone^0_{\bgg^0}$ is in $\Null{(D^0_{\bgg^0})^T}$. Moreover, if $\bgg^0$ is connected, $\bone^0_{\bgg^0}$ is the only vector spanning $\Null{(D^0_{\bgg^0})^T}$. Thus, the lemma holds when $\ell = 0$. 

Now suppose the lemma is true on level $\ell$. Our goal is to show that the lemma is also true on level $\ell+1$.
Consider an aggregate $\bgg^{\ell+1}$ on level $\ell+1$, which by definition is a set of level-$(\ell+1)$ vertices (or equivalently level-$\ell$ aggregates): $\bgg^{\ell+1} = \{ \agg^\ell_i \}$. Let $P_{u, \bgg^{\ell+1}}^\ell$ be the block diagonal matrix
\[
P_{u, \bgg^{\ell+1}}^\ell = \begin{bmatrix}
\ddots \\ & P_{u, \agg^{\ell}_i}^\ell \\ & & \ddots
\end{bmatrix}_{\agg^{\ell}_i\in \bgg^\ell}.
\]
 By the induction hypothesis and Lemma~\ref{lemma:const_preserve}, $\bone^\ell_{\agg^\ell_i} \in \Range{P_{u, \agg^\ell_i}^\ell}$ for all $\agg^\ell_i\in \bgg^{\ell+1}$. Hence, $\bone^\ell_{\bgg^{\ell+1}}$ is in $\Range{P_{u, \bgg^{\ell+1}}^\ell}$ and so $P_{u, \bgg^{\ell+1}}^\ell Q_{u, \bgg^{\ell+1}}^\ell \bone_{\bgg^{\ell+1}}^\ell = \bone_{\bgg^{\ell+1}}^\ell $. Consequently,
 \begin{equation}
(D^{\ell+1}_{\bgg^{\ell+1}})^T \bone_{\bgg^{\ell+1}}^{\ell+1} = \left( (P_{\sigma, \bgg^{\ell+1}}^\ell)^T (D^\ell_{\bgg^{\ell+1}})^T P_{u, \bgg^{\ell+1}}^\ell \right) \left( Q_{u, \bgg^{\ell+1}}^\ell \bone_{\bgg^{\ell+1}}^\ell \right) = (P_{\sigma, \bgg^{\ell+1}}^\ell)^T (D_{\bgg^{\ell+1}}^\ell)^T \bone_{\bgg^{\ell+1}}^\ell = \bzero.
\label{eq:D_null_proof}
\end{equation}
If $\bgg^{\ell+1}$ is connected, since each $\agg^\ell_i$ is a connected aggregate of level-$\ell$ vertices, $\bgg^{\ell+1}$ can also be viewed as a connected aggregate of level-$\ell$ vertices. Thus, by the induction hypothesis and Lemma~\ref{lemma:commutativity},  
\[
\Null{(D^{\ell+1}_{\bgg^{\ell+1}})^T} \subseteq \Null{(Q_{\sigma, \bgg^{\ell+1}}^\ell)^T(D^{\ell+1}_{\bgg^{\ell+1}})^T } = \Null{(D^\ell_{\bgg^{\ell+1}})^T P_{u, \bgg^{\ell+1}}^\ell} = \Span{\bone_{\bgg^{\ell+1}}^{\ell+1}}.
\]
The above set inclusion relation and \eqref{eq:D_null_proof} imply $\Null{(D^{\ell+1}_{\bgg^{\ell+1}})^T} = \Span{\bone_{\bgg^{\ell+1}}^{\ell+1}}$.  \hfill$\square$

\section{Another choice of PV trace}\label{sec:pv_trace}

In this section, we provide an alternative way to construct a PV trace $\bsigma_{\face^\ell}^{PV}$ on each face $\face^\ell = (\agg^\ell_i, \agg^\ell_j)$ that satisfies \eqref{eq:pv_condition}. Let $\agg_{\face^\ell} : = \agg^\ell_i\cup\agg^\ell_j$. Define
\[
g_{\face^\ell} = \begin{bmatrix} -\bone^\ell_{\agg^\ell_i} / \|\bone^\ell_{\agg^\ell_i}\|^2 \\ 
\bone^\ell_{\agg^\ell_j} / \|\bone^\ell_{\agg^\ell_j}\|^2 \end{bmatrix}.
\]
Then consider the following local problem defined on $\agg_{\face^\ell}$ (i.e., take $\bgg = \agg_{\face^\ell}$ in \eqref{eq:submatrices}):
\begin{equation}
 \begin{bmatrix}
    M^{\ell}_{\agg_{\face^\ell}} &  \left(D^{\ell}_{\agg_{\face^\ell}}\right)^T \\
    D^{\ell}_{\agg_{\face^\ell}} &
  \end{bmatrix} 
  \begin{bmatrix}
    \bsigma_{\agg_{\face^\ell}}^{PV} \\ \bu
  \end{bmatrix} 
  =
  \begin{bmatrix}
    0 \\ g_{\face^\ell}
  \end{bmatrix} .
  \label{eq:pv_problem}
\end{equation}
Note that by Lemma~\ref{lemma:D_null}, $\Null{(D^\ell_{\agg_{\face^\ell}})^T} = \Span{\bone^\ell_{\agg_{\face^\ell}}}$ and $\big(\bone_{\agg_{\face^\ell}}\big)^T g_{\face^\ell} = 0$, so \eqref{eq:pv_problem} is solvable. The actual PV trace is taken as the restriction of the solution $\bsigma_{\agg_{\face^\ell}}^{PV}$ to $\bSigma^\ell(\face^\ell)$. That is, $\bsigma_{\face^\ell}^{PV} := \bsigma_{\agg_{\face^\ell}}^{PV}|_{\bSigma^\ell(\face^\ell)}$.

\end{appendices}

\end{document}

%% file: dirichlet_bc_square.tex
\begin{figure}
  \begin{center}
    \includegraphics[width=0.45\textwidth]{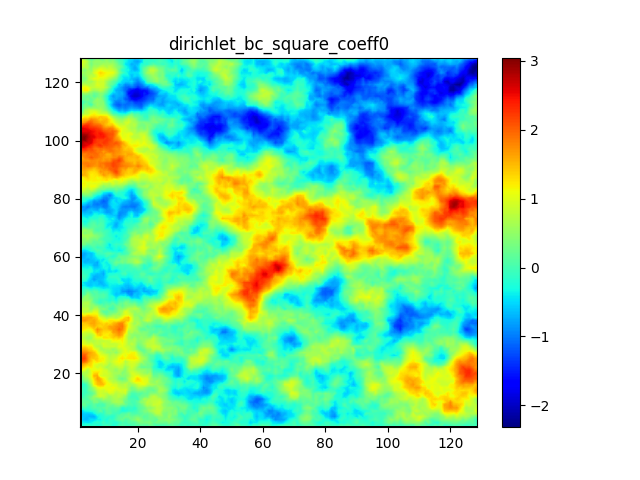}
    \includegraphics[width=0.45\textwidth]{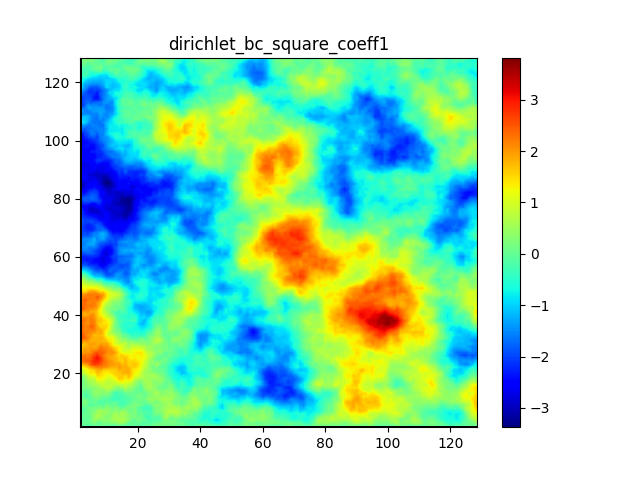}
  \end{center}
  \caption{Two typical (log) permeability field realizations for the unit square example.}
  \label{fig:dirichlet_bc_square_coeff}
\end{figure}
\begin{figure}
  \begin{center}
    \includegraphics[width=0.3\textwidth]{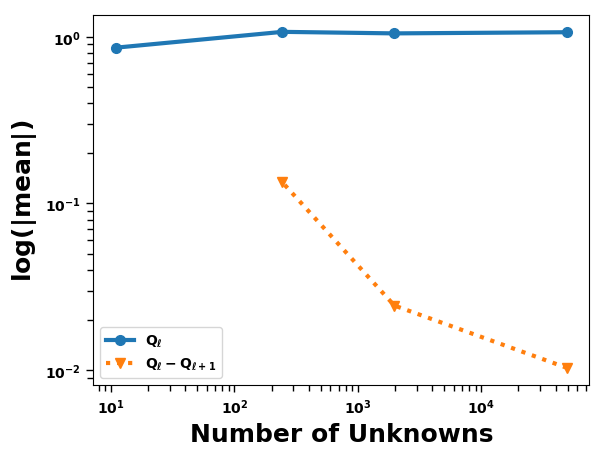}
    \includegraphics[width=0.3\textwidth]{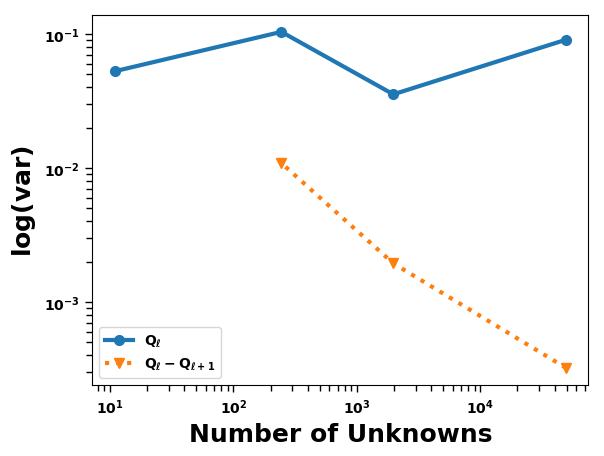}
    \includegraphics[width=0.3\textwidth]{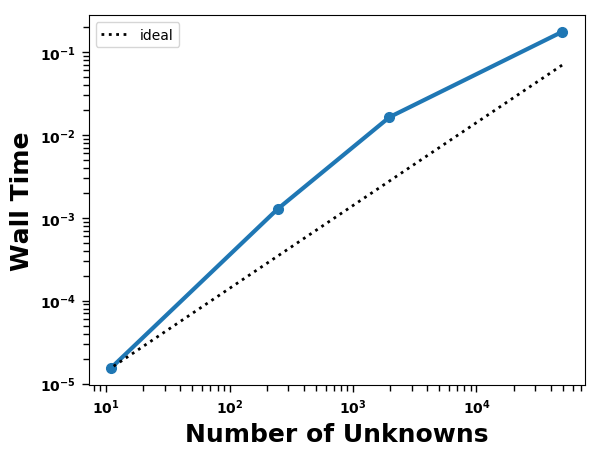}
  \end{center}
  \caption{The mean (left) and variance (middle) of the quantity of interest (or difference in quantity of interest) at several levels of the hierarchy. The scaling (right) of the algorithm with respect to problem size is seen to be close to ideal.}
  \label{fig:dirichlet_bc_square}
\end{figure}

%% file: dirichlet_square_speedup.tex
\begin{table}
\caption{Sample counts and computational cost for the unit square problem, including one-level and multilevel Monte Carlo for various tolerances on the mean-square error (MSE).
For multilevel Monte Carlo, the sample counts are listed from finest to coarsest level.}
\label{tab:dirichlet_square_speedup}
\begin{center}\begin{tabular}{llllll}
\toprule
      &  \multicolumn{2}{c}{one-level}&  \multicolumn{3}{c}{multilevel} \\
   \cmidrule(l){2-3}  \cmidrule(l){4-6}
MSE&  samples&  time(s)&  samples&  time(s)&  speedup \\
\midrule
5e-03&  19&  3.1&  [15, 15, 15, 19]&  2.9&  1.1 \\
1e-03&  72&  11.5&  [15, 15, 15, 110]&  2.9&  4.0 \\
5e-04&  129&  20.7&  [15, 15, 15, 1143]&  2.9&  7.1 \\
1e-04&  614&  98.1&  [15, 32, 282, 7739]&  3.6&  27.1 \\
1e-05&  6773&  1081.8&  [108, 785, 4525, 120036]&  39.4&  27.4 \\
\bottomrule
\end{tabular}\end{center}
\end{table}

%% file: dirichlet_bc_cube.tex
\begin{figure}
  \begin{center}
    \includegraphics[width=0.45\textwidth]{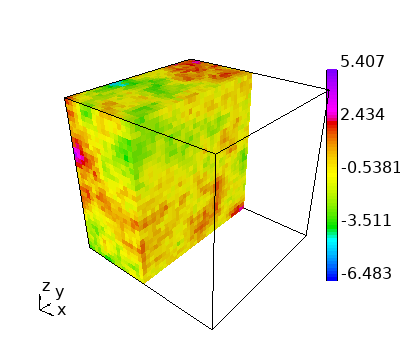}
    \includegraphics[width=0.45\textwidth]{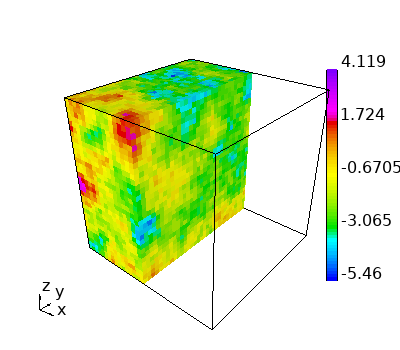}
  \end{center}
  \caption{Two typical (log) permeability field realizations for the unit cube example.}
  \label{fig:dirichlet_bc_cube_coeff}
\end{figure}
\begin{figure}
  \begin{center}
    \includegraphics[width=0.3\textwidth]{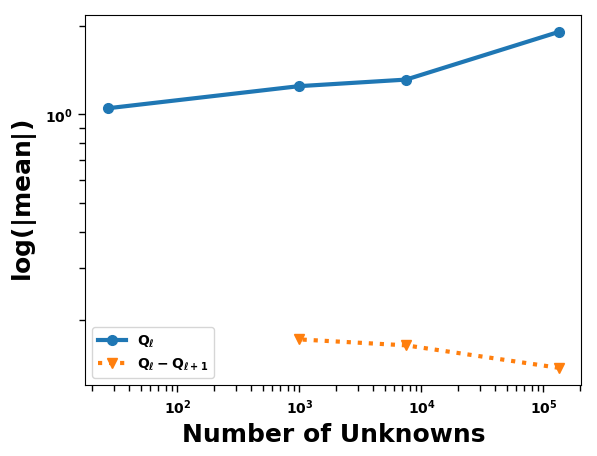}
    \includegraphics[width=0.3\textwidth]{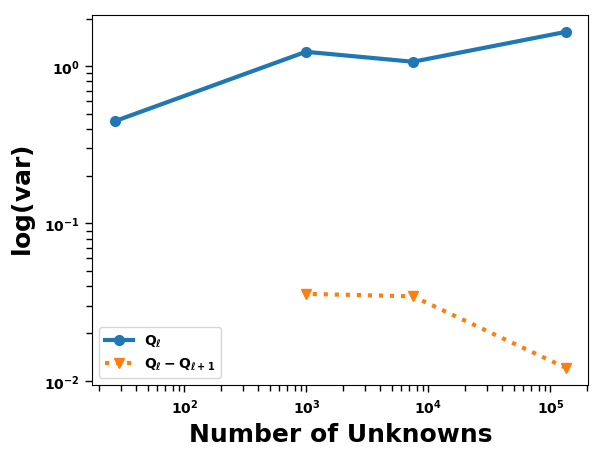}
    \includegraphics[width=0.3\textwidth]{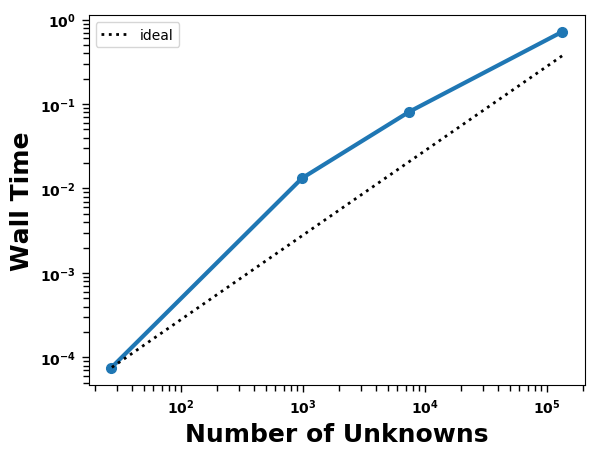}
  \end{center}
  \caption{The mean (left) and variance (middle) of the quantity of interest (or difference in quantity of interest) at several levels of the hierarchy. The scaling (right) of the algorithm with respect to problem size is seen to be close to ideal.}
  \label{fig:dirichlet_bc_cube}
\end{figure}

%% file: dirichlet_cube_speedup.tex
\begin{table}
\caption{Sample counts and computational cost for the unit cube problem, including one-level and multilevel Monte Carlo for various tolerances on the mean-square error (MSE).
For multilevel Monte Carlo, the sample counts are listed from finest to coarsest level.}
\label{tab:dirichlet_cube_speedup}
\begin{center}\begin{tabular}{llllll}
\toprule
      &  \multicolumn{2}{c}{one-level}&  \multicolumn{3}{c}{multilevel} \\
   \cmidrule(l){2-3}  \cmidrule(l){4-6}
MSE&  samples&  time(s)&  samples&  time(s)&  speedup \\
\midrule
5e-03&  390&  18.5&  [15, 66, 146, 382]&  2.4&  7.7 \\
1e-03&  1514&  70.9&  [28, 347, 547, 1858]&  9.0&  7.9 \\
5e-04&  2771&  128.1&  [76, 575, 1109, 3495]&  17.5&  7.3 \\
1e-04&  13422&  624.9&  [434, 2395, 5298, 17890]&  85.5&  7.3 \\
\bottomrule
\end{tabular}\end{center}
\end{table}

%% file: spe10.tex
\begin{figure}
  \begin{center}
    \includegraphics[width=\textwidth]{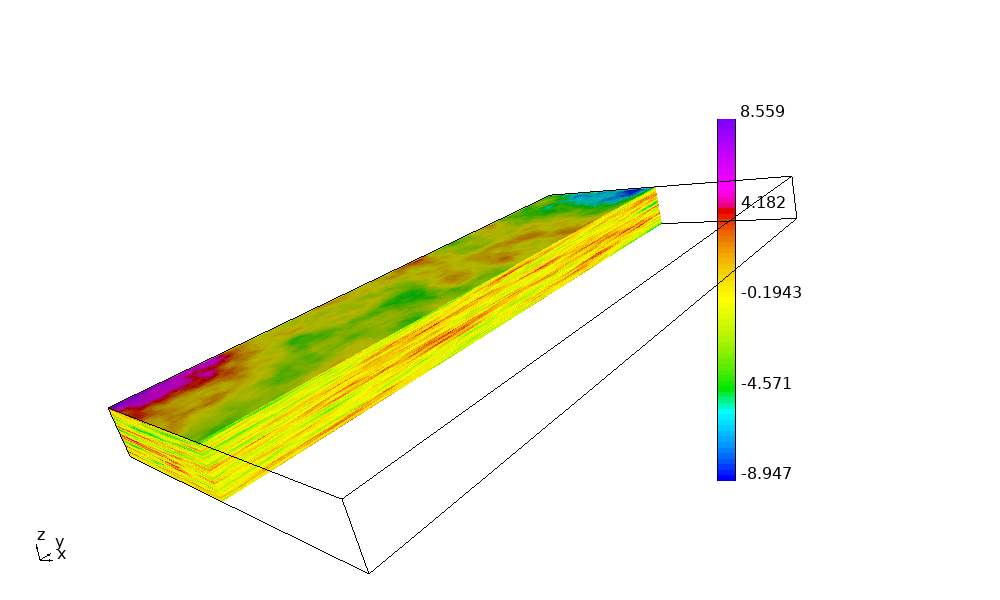}
  \end{center}
  \caption{A typical (log) permeability field realization for the SPE10 example.
The anisotropic covariance leads to a layering structure for the coefficient, where the correlation length is longer in the horizontal directions than in the vertical direction.}
  \label{fig:spe10_coeff}
\end{figure}
\begin{figure}
  \begin{center}
    \includegraphics[width=0.3\textwidth]{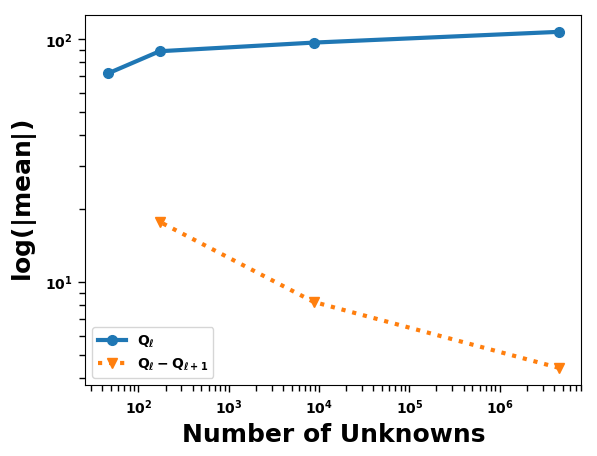}
    \includegraphics[width=0.3\textwidth]{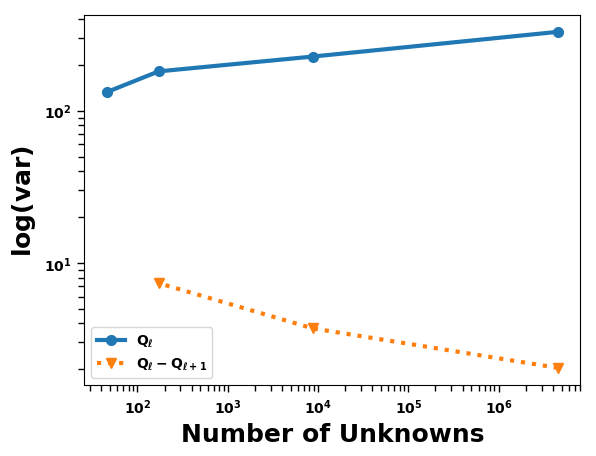}
    \includegraphics[width=0.3\textwidth]{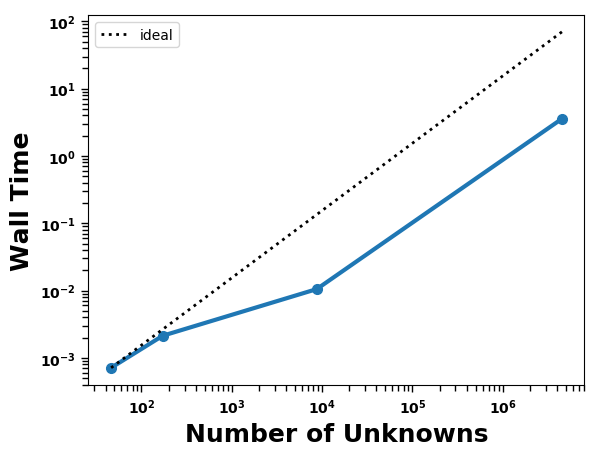}
  \end{center}
  \caption{The mean (left) and variance (middle) of the quantity of interest (or difference in quantity of interest) at several levels of the hierarchy for the SPE10 example. The scaling (right) of the algorithm with respect to problem size is seen to be close to ideal.}
  \label{fig:spe10}
\end{figure}

%% file: multilevel.bbl
\begin{thebibliography}{10}

\bibitem{glvis}
{GLVis}: {OpenGL} finite element visualization tool.
\newblock \url{glvis.org}.

\bibitem{mfem}
{MFEM}: Modular finite element methods library.
\newblock \url{mfem.org}.

\bibitem{aeev2012}
A.~Abdulle, W.~E, B.~Engquist, and E.~Vanden-Eijnden.
\newblock The heterogeneous multiscale method.
\newblock {\em Acta Numerica}, 21:1?87, 2012.

\bibitem{ArPeWY07}
T.~Arbogast, G.~Pencheva, M.~Wheeler, and I.~Yotov.
\newblock A multiscale mortar mixed finite element method.
\newblock {\em Multiscale Model. Simul.}, 6(1):319--346 (electronic), 2007.

\bibitem{ArnoldBrezzi85}
D.~N. Arnold and F.~Brezzi.
\newblock Mixed and nonconforming finite element methods : implementation,
  postprocessing and error estimates.
\newblock {\em ESAIM: M2AN}, 19(1):7--32, 1985.

\bibitem{smoothg}
A.~T. Barker, S.~A. Gelever, C.~S. Lee, C.~V. Ponce, and P.~S. Vassilevski.
\newblock {smoothG}.
\newblock \url{github.com/LLNL/smoothG}, 2017.

\bibitem{blv17}
A.~T. Barker, C.~S. Lee, and P.~S. Vassilevski.
\newblock Spectral upscaling for graph {L}aplacian problems with application to
  reservoir simulation.
\newblock {\em SIAM J. Sci. Comput.}, 39:S323--S346, 2017.

\bibitem{cfhjmmrv03}
T.~Chartier, R.~Falgout, V.~Henson, J.~Jones, T.~Manteuffel, S.~McCormick,
  J.~Ruge, and P.~Vassilevski.
\newblock Spectral {AMG}e ($\rho${AMG}e).
\newblock {\em SIAM J. Scientific Computing}, 25:1--26, 2003.

\bibitem{cps07}
G.~A. Chechkin, A.~L. Piatnitski, and A.~S. Shamaev.
\newblock {\em Homogenization: Methods and applications}, volume 234 of {\em
  Translations of Mathematical Monographs}.
\newblock American Mathematical Society, 2007.

\bibitem{spe10}
M.~A. Christie and M.~J. Blunt.
\newblock Tenth {SPE} comparative solution project: A comparison of upscaling
  techniques.
\newblock {\em Society of Petroleum Engineers}, 2001.

\bibitem{cliffe-giles-mlmc}
K.~A. Cliffe, M.~B. Giles, R.~Scheichl, and A.~L. Teckentrup.
\newblock Multilevel {M}onte {C}arlo methods and applications to elliptic
  {PDE}s with random coefficients.
\newblock {\em Computing and Visualization in Science}, 14(1):3, 2011.

\bibitem{dkltv19}
V.~Dobrev, T.~Kolev, C.~S. Lee, V.~Tomov, and P.~S. Vassilevski.
\newblock Algebraic hybridization and static condensation with application to
  scalable {$H$(div)} preconditioning.
\newblock {\em SIAM Journal on Scientific Computing}, 41(3):B425--B447, 2019.

\bibitem{ee2003}
W.~E and B.~Engquist.
\newblock The heterognous multiscale methods.
\newblock {\em Commun. Math. Sci.}, 1(1):87--132, 03 2003.

\bibitem{egh12}
Y.~Efendiev, J.~Galvis, and T.~Hou.
\newblock Generalized multiscale finite element methods.
\newblock {\em Journal of Computational Physics}, 251:116--135, 2013.

\bibitem{EfendievHouBook}
Y.~Efendiev and T.~Y. Hou.
\newblock {\em Multiscale Finite Element Methods. Theory and Applications}.
\newblock Springer Science+Business Media, LLC, 1st edition, 2009.

\bibitem{networkx}
A.~A. Hagberg, D.~A. Schult, and P.~J. Swart.
\newblock Exploring network structure, dynamics, and function using {NetworkX}.
\newblock In G.~Varoquaux, T.~Vaught, and J.~Millman, editors, {\em Proceedings
  of the 7th Python in Science Conference}, pages 11 -- 15, Pasadena, CA USA,
  2008.

\bibitem{h03}
V.~E. Henson.
\newblock Multigrid methods for nonlinear problems: an overview.
\newblock In {\em Proc. SPIE}, volume 5016, pages 36--48, 2003.

\bibitem{hornung97}
U.~Hornung.
\newblock {\em Homogenization and porous media}, volume~6 of {\em
  Interdisciplinary Applied Mathematics}.
\newblock Springer, New York, 1997.

\bibitem{hw97}
T.~Hou and X.~Wu.
\newblock A multiscale finite element method for elliptic problems in composite
  materials and porous media.
\newblock {\em J. Comput. Phys.}, 134:169--189, 1997.

\bibitem{EggModel}
J.~D. Jansen, R.~M. Fonseca, S.~Kahrobaei, M.~M. Siraj, G.~M. Van~Essen, and
  P.~M.~J. Van~den Hof.
\newblock The egg model - a geological ensemble for reservoir simulation.
\newblock {\em Geoscience Data Journal}, 1(2):192--195, 2014.

\bibitem{jennylt03}
P.~Jenny, S.~Lee, and H.~Tchelepi.
\newblock Multi-scale finite volume method for elliptic problems in subsurface
  flow simulation.
\newblock {\em J. Comput. Phys.}, 187:47--67, 2003.

\bibitem{jlt06}
P.~Jenny, S.~Lee, and H.~Tchelepi.
\newblock Adaptive fully implicit multi-scale finite-volume method for
  multi-phase flow and transport in heterogeneous porous media.
\newblock {\em Journal of Computational Physics}, 217(2):627 -- 641, 2006.

\bibitem{kalchev-lee-upscaling-mixed}
D.~Z. Kalchev, C.~S. Lee, U.~Villa, Y.~Efendiev, and P.~S. Vassilevski.
\newblock Upscaling of mixed finited element discretization problems by the
  spectral {AMG}e method.
\newblock {\em SIAM J. Sci. Comput.}, 38(5):A2912--2933, 2016.

\bibitem{lv08}
I.~V. Lashuk and P.~S. Vassilevski.
\newblock On some versions of the element agglomeration {AMG}e method.
\newblock {\em Numer. Linear Algebra Appl.}, 15:595--620, 2008.

\bibitem{fas-spectral}
C.~S. Lee, N.~Castelletto, F.~Hamon, P.~S. Vassilevski, and J.~White.
\newblock Nonlinear multigrid based on local spectral coarsening.
\newblock {\em In preparation}.

\bibitem{lv16}
C.~S. Lee and P.~S. Vassilevski.
\newblock Parallel solver for {H}(div) problems using hybridization and {AMG}.
\newblock {\em LLNL Report}, LLNL-TR-681025, 2016.

\bibitem{lrl11}
F.~Lindgren, H.~Rue, and J.~Lindstr{\"o}m.
\newblock An explicit link between {G}aussian fields and {G}aussian {M}arkov
  random fields: the stochastic differential equation approach.
\newblock {\em Journal of the Royal Statistical Society: Series B (Statistical
  Methodology)}, 73:423--498, 2011.

\bibitem{mm06}
S.~MacLachlan and J.~Moulton.
\newblock Multilevel upscaling through variational coarsening.
\newblock {\em Water Resour. Res.}, 42:W02418, doi:10.1029/2005WR003940, 2006.

\bibitem{ovv17}
S.~Osborn, P.~S. Vassilevski, and U.~Villa.
\newblock A multilevel, hierarchical sampling technique for spatially
  correlated random fields.
\newblock {\em SIAM J. Sci. Comput.}, 39:S543--S562, 2017.

\bibitem{pasciak-vassilevski}
J.~E. Pasciak and P.~S. Vassilevski.
\newblock Exact de {R}ham sequences of spaces defined on macro-elements in two
  and three spatial dimensions.
\newblock {\em SIAM J. Sci. Comput.}, 30(5):2427--2446, 2008.

\bibitem{classicAMG}
J.~W. Ruge and K.~St\"{u}ben.
\newblock Algebraic multigrid {(AMG)}.
\newblock In S.~F. McCormick, editor, {\em Multigrid Methods}, volume~3 of {\em
  Frontiers in Applied Mathematics}, pages 73--130. SIAM, Philadelphia, PA,
  1987.

\bibitem{simpson-lindgren-continuous}
D.~Simpson, F.~Lindgren, and H.~Rue.
\newblock Think continuous: {M}arkovian {G}aussian models in spatial
  statistics.
\newblock {\em Spatial Statistics}, 1:16--29, 2012.

\bibitem{tjl07}
H.~A. Tchelepi, P.~Jenny, S.~H. Lee, and C.~Wolfsteiner.
\newblock Adaptive multiscale finite-volume framework for reservoir simulation.
\newblock {\em SPE Journal}, 12(2):188--195, 2007.

\bibitem{vassilevski-upscaling}
P.~S. Vassilevski.
\newblock Coarse spaces by algebraic multigrid: Multigrid convergence and
  upscaling error estimates.
\newblock {\em Advances in Adaptive Data Analysis}, 3(1-2):229--249, 2011.

\bibitem{wdle00}
X.-H. Wen, L.~J. Durlofsky, S.~H. Lee, and M.~G. Edwards.
\newblock Full tensor upscaling of geologically complex reservoir descriptions.
\newblock In {\em SPE Annual Technical Conference and Exhibition}, page~11.
  Society of Petroleum Engineers, Dallas, Texas, 2000.

\bibitem{wg96}
X.-H. Wen and J.~G\'{o}mez-Hern\'{a}ndez.
\newblock Upscaling hydraulic conductivities in heterogeneous media: An
  overview.
\newblock {\em Journal of Hydrology}, 183(1):ix -- xxxii, 1996.

\end{thebibliography}
